\PassOptionsToPackage{dvipsnames}{xcolor} 

\documentclass[10pt,a4paper]{amsart}
\usepackage[a4paper, left=1in, right=1in, top=1.1in, bottom=1.1in]{geometry}

\usepackage[english]{babel}
\usepackage[utf8]{inputenc}
\usepackage{amsmath}
\usepackage{graphicx}
\usepackage{amssymb}
\usepackage{stmaryrd}
\usepackage{amsthm}
\usepackage{tikz-cd}
\usepackage{mathrsfs}
\usepackage{csquotes}
\usepackage[colorinlistoftodos]{todonotes}
\usepackage{yfonts}
\usepackage{ dsfont }
\usepackage{xcolor}
\usepackage{multicol}
\usepackage{indentfirst}
\usepackage[hyphens]{url}
\usepackage{hyperref}
\usepackage[hyphenbreaks]{breakurl}
\usepackage{subfiles}

 \usepackage[
    backend=biber,
    style=alphabetic, maxbibnames=99, url=false
  ]{biblatex}
\usepackage{xurl}
\hypersetup{
    colorlinks,
    linkcolor={red!70!black},
    citecolor={blue!70!black},
    urlcolor={violet}
}
\hypersetup{breaklinks=true}

\usepackage[inline]{enumitem}   
\makeatletter

\newcommand{\inlineitem}[1][]{%
\ifnum\enit@type=\tw@
    {\descriptionlabel{#1}}
  \hspace{\labelsep}%
\else
  \ifnum\enit@type=\z@
       \refstepcounter{\@listctr}\fi
    \quad\@itemlabel\hspace{\labelsep}%
\fi}
\makeatother
\allowdisplaybreaks

 \usepackage[
    backend=biber,
    style=alphabetic, maxbibnames=99, url=false
  ]{biblatex}
\usepackage[capitalise]{cleveref}
\title[Conformal Block Divisors for Discrete Series Virasoro VOA $\vir{2k+1}$]{Conformal Block Divisors for\\ Discrete Series Virasoro VOA $\vir{2k+1}$}
\subjclass[2020]{14H10, 17B69 (primary), 81R10 (secondary)}
\keywords{Moduli of curves, Vertex Operator Algebras}

\author{Daebeom Choi}
\address{Department of Mathematics\\
    University of Pennsylvania\\
    Philadelphia, PA 19104-6395}
\email{dbchoi@sas.upenn.edu}

\date{\today}
\theoremstyle{definition}
\newtheorem{thm}{Theorem}[section]

\newtheorem{lem}[thm]{Lemma}

\newtheorem{prop}[thm]{Proposition}
\newtheorem{qes}[thm]{Question}

\newtheorem{defn}[thm]{Definition}
\newtheorem{eg}[thm]{Example}

\newtheorem{conj}[thm]{Conjecture}
\newtheorem{defthm}[thm]{Definition-Theorem}

\newtheorem{cor}[thm]{Corollary}

\newtheorem{rmk}[thm]{Remark}

\newcommand{\R}{\mathbb{R}}
\newcommand{\C}{\mathbb{C}}
\newcommand{\Z}{\mathbb{Z}}
\newcommand{\N}{\mathbb{N}}
\newcommand{\Q}{\mathbb{Q}}
\newcommand{\rank}{\text{rank }}

\newcommand{\M}[2]{\overline{\rm{M}}_{#1, #2}}
\newcommand{\vir}[1]{\text{Vir}_{#1,2}}
\newcommand{\NE}[1]{\overline{\text{NE}}_1(#1)}

\DeclareSymbolFont{yhlargesymbols}{OMX}{yhex}{m}{n}
\DeclareMathAccent{\widetriangle}{\mathord}{yhlargesymbols}{"E6} 

\begin{document}

\begin{abstract}
    In this work, we study a family of vector bundles on the moduli space of curves constructed from representations of $\vir{2k+1}$, a family of vertex operator algebras derived from the Virasoro Lie algebra. Using the relationship between rank and degree, we characterize their asymptotic behavior, demonstrating that their first Chern classes are nef on $\M{g}{n}$ in many cases. This is the first nontrivial example of divisors arising from vertex operator algebras that are uniformly positive for all genera. Furthermore, for $g = 1$, these divisors form a $\mathbb{Q}$-basis of the Picard group of $\M{1}{n}$, with several desirable functorial properties. Using this basis, we characterize line bundles on certain contractions of $\M{1}{n}$. We also propose conjectures regarding the conformal blocks of Virasoro VOAs and potential generalizations. In particular, by introducing a generalization of conformal block divisors, we provide a nonlinear nef interpolation between affine and Virasoro conformal block divisors.
\end{abstract}

\maketitle

\section{Introduction}\label{sec:intro}

The moduli space of stable pointed curves, $\M{g}{n}$, is one of the most fundamental objects in algebraic geometry, yet its birational geometry remains largely mysterious. A key step toward understanding the birational geometry of the moduli space is to describe positive divisors, those that lie in the ample cone, and its closure, the nef cone. In this work, we identify the first nontrivial family of divisors, derived from representations of a vertex operator algebra, with positivity properties on $\M{g}{n}$ for all $g,n$.

There is a conjectural description of the nef cone of $\M{g}{n}$. The F-cone is an explicit closed polyhedral subcone of the Picard group, known to contain the nef cone. A line bundle is said to be \textbf{F-nef} (resp. \textbf{F-ample}) if it lies within the F-cone (resp. the interior of the F-cone). The F-conjecture asserts that the nef cone coincides with the F-cone. Various techniques have been developed for showing that F-nef divisors are nef \cite{KM13, GKM02, GF03, Gib09, Lar11, Fe15, MS19, Fe20}. Nevertheless, the F-conjecture remains open (see \cref{subsec:F} for a discussion).

One approach to the F-conjecture is to produce a large family of nef line bundles that one could hope to show collectively cover the entire F-cone. A large source of line bundles in genus $0$ arises from certain vertex operator algebras. \cite{DGT21, DGT22a, DGT22b} proves that for any $C_2$-cofinite and rational VOA $V$ and simple modules $W^1,\cdots,W^n$, there exist a corresponding vector bundle $\mathbb{V}_{g,n}(V, W^\bullet)$ on $\M{g}{n}$, referred to as the \textbf{sheaf of coinvariants}, whose first Chern class $\mathbb{D}_{g,n}(V, W^\bullet)$ we refer to as \textbf{coinvariant divisor}, and $-\mathbb{D}_{g,n}(V, W^\bullet)$ we call the \textbf{conformal block divisor}. We refer to \cref{def:CBD}  and \cref{rmk:cdiv} for terminology and \cref{subsec:cdf} for an overview of conformal block divisors.

If $V=L_k(\mathfrak{g})$, the affine VOA of level $k$ corresponds to $\mathfrak{g}$,  (cf. \cite{TK88, TUY89}) or, more generally, is generated in degree $1$, then \cite{Fak12, DG23} proved that $\mathbb{D}_{0,n}(V, W^\bullet)$ is base point free. However, this result is limited to the case $g = 0$ (cf. \cite[Section 6]{Fak12}) and has not been sufficient to resolve the F-conjecture. Here, we describe a new family of nef divisors on $\M{g}{n}$ for $g,n \ge 0$, which may be expressed as linear combinations of conformal block divisors (see \cref{thm:main1}, \cref{thm:main2}, \cref{thm:main3}). As we show, these divisors are nef even in positive genus. Moreover, as $\vir{2k+1}$ is not generated in degree $1$, the proof of their positivity is different from previous approaches (cf. \cref{thm:virdeg}). 

\begin{thm}\label{thm:main1}
    For any $g, n$, and $n$ simple modules $W^i$, $-\mathbb{D}_{g,n}(\vir{2k+1}, W^\bullet)$ is F-nef. Furthermore, if $W^i$ are nontrivial, then $-\mathbb{D}_{g,n}(\vir{2k+1}, W^\bullet)$ is either zero or F-ample. For $k < 9$, $-\mathbb{D}_{g,n}(\vir{2k+1}, W^\bullet)$ is either zero or ample.
\end{thm}

In \cref{thm:main1}, $\vir{2k+1}$ denotes a family of discrete series Virasoro VOAs, which we will review in \cref{subsec:VOA}. The divisor $-\mathbb{D}_{g,n}(\vir{2k+1}, W^\bullet)$ will be referred to as \textbf{Virasoro conformal block divisors}. Notably, this is the first example of a non-holomorphic VOA whose conformal block divisor is uniformly positive for all $g$. Furthermore, we emphasize that it is the conformal block divisors (not coinvariant divisors) for $\vir{2k+1}$ that are nef. This nefness of the conformal block divisors is unexpected since, before this work, only some coinvariant divisors and no conformal block divisors have been shown to be nef. As an application of \cref{thm:main1}, we confirm Fakhruddin's prediction \cite[\nopp 5.2.8]{Fak12} that genus $0$ coinvariant divisors for $L_1(F_4)$ and $L_1(G_2)$ are ample (see \cref{cor:f4g2}).

To prove \cref{thm:main1}, we proceed by induction on $n$. A key ingredient of the proof is the factorization of conformal blocks \cite{TUY89, DGT22a}, explained in \cref{thm:factor}. This allows us to handle the restriction of conformal block divisors to boundary divisors. We also utilize \cref{lem:effnef} and \cref{lem:Knef}, which are variants of standard tools aimed at applications to conformal block divisors. Such standard tools have appeared in previous work, such as \cite{Gib09}, but factorization makes the induction step much more feasible.

Unfortunately, we were unable to prove the ampleness part of \cref{thm:main1} for every $k$. However, we propose \cref{conj:genvireff}, which can be verified through a finite amount of computation for each $k$, and implies \cref{thm:main1} for larger values of $k$. Furthermore, we establish the positivity for another important class of Virasoro conformal block divisors. Note that any simple module over $\vir{2k+1}$ corresponds to $1 \le i \le k$, and let $W_i$ denote the corresponding modules.

\begin{thm}\label{thm:main2}
    Let $W^i=W_{a_i}$, $1\le a_i\le k$, and $\sum_{i=1}^n(a_i-1)\le 2(k-1)$. Then
    \begin{enumerate}
        \item $\mathbb{D}_{0,n}(\vir{2k+1}, W^\bullet)$ are independent of the specific value of $k$.
        \item $-\mathbb{D}_{0,n}(\vir{2k+1}, W^\bullet)$ is nef. If $a_i\ne 1$ for every $i$, this is either zero or ample.
    \end{enumerate}
\end{thm}

In particular, this theorem provides infinitely many ample line bundles on $\M{0}{n}$ for each $n\ge 4$.

\cref{thm:main2} (1) establishes the stabilization, or the ``critical level" phenomenon, of the Virasoro conformal block divisors. Namely, if $k$ is above the critical level (cf. \cref{defn:critlev}), then $-\mathbb{D}_{0,n}(\vir{2k+1}, W^\bullet)$ does not depend on $k$, and we refer to the limit of $-\mathbb{D}_{0,n}(\vir{2k+1}, W^\bullet)$ as the \textbf{stable Virasoro conformal block divisor}. This behavior parallels the critical level phenomena observed in affine VOAs \cite{Fak12,BGM15,BGM16}, where $\mathbb{D}_{0,n}(L_k(\mathfrak{g}), W^\bullet)$ converges to a trivial line bundle when $k$ is above the critical level. Like \cref{thm:main1}, this reveals that the conformal block divisors of $\vir{2k+1}$ and affine VOAs share a similar property, but for seemingly different reasons. However, in contrast to the affine VOA case, Virasoro conformal block divisors generically converge to a nontrivial line bundle if $\sum_i(a_i-1)$ is even. This stabilization phenomenon arises from the computation presented in \cref{thm:virdeg}, which demonstrates that on $\M{0}{4}$, the degree of $\mathbb{V}_{0,n}(\vir{2k+1}, W^\bullet)$ is entirely determined by its rank. This property of the sheaf for $\vir{2k+1}$ that the rank determines its degree sets it apart from other VOAs (see \cref{rmk:rank} for detail). Moreover, (2) shows that stable Virasoro conformal block divisors exhibit the predicted positivity properties.

So far in this paper, all the conformal block divisors whose positivity we have proven are either zero, ample, or pullback of such along a projection map; hence, they are not contained in the interior of an interesting face of the nef cone. However, we can also construct a wide range of F-nef line bundles that are not necessarily F-ample or zero by taking the difference of Virasoro conformal block divisors.

\begin{thm}\label{thm:main3}
    Let $W^i=W_{a_i}$, $1\le a_1\le a_2\cdots \le a_n\le k$.
    \begin{enumerate}
        \item If $\sum_{i=1}^n(a_i-1)$ is even, then $\mathbb{D}_{0,n}(\vir{2k+1}, W^\bullet)-\mathbb{D}_{0,n}(\vir{2k+3}, W^\bullet)$ is F-nef.
        \item If $\sum_{i=1}^n(a_i-1)$ is odd, $n \ge 5$, $a_i > 1$, and $l(a_1, \cdots, a_n) - 3 \le k$, then
    \[ \mathbb{D}_{0,n}(\vir{2k+3}, W^\bullet) - \mathbb{D}_{0,n}(\vir{2k+1}, W^\bullet) \]
    is F-nef, where $l(a_1, \cdots, a_n)$ is the smallest integer $l$ such that $\sum_{i=1}^n(a_i-1) \le 2(l-1)$. Moreover, if $l(a_1, \cdots, a_n) - 1 \le k$, then $\mathbb{D}_{0,n}(\vir{2k+1}, W^\bullet) = 0$, even when $n = 4$.
    \end{enumerate}
\end{thm}

Note that the signs of the differences in (1) and (2) are opposite. This difference in signs arises from the differing behavior of the fusion rules in the odd and even cases. Unfortunately, we were unable to prove the nefness of these line bundles, as their explicit representations in terms of tautological line bundles are too complicated. Although \cref{thm:main2} and \cref{thm:main3} only apply to the case $g=0$, they still have the potential to contribute to the solution of the F-conjecture, as \cite{GKM02} shows that the F-conjecture for genus $0$ implies the F-conjecture for all genera.

Furthermore, we prove another interesting property of the conformal block divisors of $\vir{5}$. In \cite[Theorem 4.3]{Fak12}, Fakhruddin proved that conformal block divisors for $L_1(\mathfrak{sl}_2)$ form a basis of the Picard group of $\M{0}{n}$. Here, we establish an analogous theorem for $\M{1}{n}$. Note that $\vir{5}$ has two simple modules, $W_1$ and $W_2$, where $W_1$ is trivial and $W_2$ is nontrivial.

\begin{thm}\label{thm:main4}
    If exactly one among $W^1, \ldots, W^n$ is $W_2$, then $\mathbb{V}_{1,n}(V, W^\bullet)$ is a trivial line bundle. For any other choices, $\mathbb{D}_{1,n}(V, W^\bullet)$ is nontrivial, and such conformal block divisors form a basis of $\text{Pic}(\overline{\mathrm{M}}_{1,n})_\Q$.
\end{thm}

Using the basis from \cref{thm:main4}, one may characterize line bundles on certain contractions of $\M{1}{n}$, as was done via Fakhruddin's basis in genus $0$ in \cite[Theorem 1.1]{Cho24}.

\begin{cor}\label{cor:main4cor}
    For nonempty subsets $S,T\subseteq [n]$, let $f:\M{1}{n}\to X_{1, S, T}$ be a contraction corresponds to $\M{1}{n}\to \M{1}{S}\times \M{1}{T}$. Then a divisor on $\M{1}{n}$ is a pullback along $f$ if and only if it intersects trivially with the F-curves contracted by $f$ if and only if $D$ is a pullback of a divisor along $\M{1}{n}\to \M{1}{S}\times \M{1}{T}$. 
\end{cor}

\begin{cor}\label{cor:main4cor2}
     If $D$ is a divisor that forms an extremal ray of the nef cone of $\M{1}{n}$, then the pullback of $D$ along any projection map $\pi:\M{1}{m}\to \M{1}{n}$ is an extremal ray of the nef cone of $\M{1}{m}$.
\end{cor}

Here, \textbf{F-curves} denote the class of curves corresponding to the irreducible components of 1-dimensional strata of $\M{g}{n}$. \cref{cor:main4cor} has two key advantages. First, the contracted curves fully determine the line bundle on contractions. Second, only F-curves, not all curves, are needed for this characterization. We refer to \cref{subsec:F} for a detailed description. 

As mentioned above, the main theorems proved here indicate that the family $\vir{2k+1}$ for $k\in \N$ exhibit highly distinctive properties, setting them apart even from other discrete series Virasoro VOAs. In \cref{subsec:vir}, we explore several remarkable features of $\vir{2k+1}$ in the context of conformal field theory and probability. Notably, unlike other discrete series Virasoro VOAs, their partition functions satisfy the KdV hierarchy. Regarding this, we raise the question of whether these features of the partition function imply certain properties of sheaves of coinvariants. We propose further conjectures concerning the conformal blocks of $\vir{2k+1}$.

Note that the proofs of positivity for conformal block divisors rely solely on the factorization theorem and explicit representations in terms of tautological line bundles. In \cref{subsec:ind}, we propose the concept of an \textbf{inductive system of line bundles}, which generalizes conformal block divisors from a VOA and extends Fedorchuk's notion of `divisors from symmetric functions on abelian groups' \cite{Fe15}. This framework retains all the desirable properties of conformal block divisors while offering greater flexibility, as it is not constrained to $C_2$-cofinite and rational VOAs, which are generally difficult to construct. As a toy example, we prove the following.

\begin{thm}\label{thm:main5}
    Let $p \ge 0$ be a nonnegative real number and
    \[R_n^p := \frac{1}{\sqrt{p^2+4}}\left(\left(\frac{p+\sqrt{p^2+4}}{2} \right)^{n-1} - \left(\frac{p-\sqrt{p^2+4}}{2} \right)^{n-1} \right). \]
    Then
    \[D_n^p := R_n^p\cdot\psi - \sum_{i=2}^{\lfloor \frac{n}{2} \rfloor} R_{i+1}^p \cdot R_{n-i+1}^p \cdot \delta_{0,i} \in \text{Pic}(\M{0}{n})_\R^{S_n} \]
It is a symmetric nef line bundle and ample if $p > 0$.
\end{thm}

Note that for $p=0$, this is the coinvariant divisor for $L_1(\mathfrak{sl}_2)$, and for $p=1$, this is the $\vir{5}$-conformal block divisor. Hence, this theorem interpolates and extends the construction of these nef divisors. Such interpolation is nontrivial as the formula that describes it is nonlinear.

Since the theory of conformal block divisors has been developed over $\C$ so far, the main theorems of this paper a priori hold over a field of characteristic zero. However, this obstruction can be circumvented without much difficulty. By \cref{thm:c1cdf}, the conformal block divisors can be expressed in terms of tautological line bundles on $\M{g}{n}$, which remain well-defined in characteristic $p$. Furthermore, since our methods do not depend on the characteristic of the base field, all of the main theorems also hold over any field. See also \cref{rmk:char}.

\subsection{Structure of the paper}
\cref{sec:pre} provides essential preliminaries for the paper, including definitions of vertex operator algebras, conformal block divisors, and a description of the F-conjecture. In \cref{subsec:fnef}, we will prove that the Virasoro conformal block divisors are F-nef. Building on this result, we will establish that they are indeed nef in \cref{subsec:nef}, thereby proving \cref{thm:main1}. \cref{subsec:rmk} provides additional material on an alternative approach to proving nefness and explores its connection to the conformal block divisors for $L_1(F_4)$ and $L_1(G_2)$. \cref{subsec:crit} focuses on the stabilization phenomena of Virasoro conformal block divisors, where \cref{thm:main2} will be proved. Using the results from this section, we will prove \cref{thm:main3} in \cref{subsec:diff}. In \cref{sec:basis}, we establish \cref{thm:main4} regarding a basis for the Picard group of $\M{1}{n}$ and prove \cref{cor:main4cor}. Finally, \cref{sec:dis} explores special properties and conjectures about $\vir{2k+1}$ and proposes an alternative approach for handling line bundles on $\M{g}{n}$ with an inductive structure. In particular, we will prove \cref{thm:main5} in \cref{subsec:ind}.

\section*{Acknowledgement}

The author would like to thank Angela Gibney for her helpful discussions, continued support, and encouragement. We are grateful to Andres Fernandez Herrero for his valuable advice with the proof of \cref{thm:gencont}. The author would also like to thank Prakash Belkale, Najmuddin Fakhruddin, and Swarnava Mukhopadhyay for useful feedback on the draft of the paper. Avik Chakravarty dedicated significant time to discussing this paper and related results, significantly improving its quality. We also appreciate Colton Griffin for his help with the proof of \cref{lem:Fibo} and for the talk that inspired the author to write this paper. Additionally, the author is indebted to Shengjing Xu for insightful comments on the KdV minimal models.

\section{Preliminary}\label{sec:pre}

\subsection{Notations and conventions}

Throughout, any Picard group $\text{Pic}(X)$, cone of curves $\NE{X}$, Chow group $A_d(X)$, curve class $[C]$, and divisors/line bundles will be considered over $\Q$. Thus, unless otherwise stated, these terms will refer to their $\Q$-versions, such as the $\Q$-Picard group, $\Q$-divisors, and so on. $\overline{\mathcal{M}}_{g,n}$ will denote the moduli stack of stable curves, while $\M{g}{n}$ will refer to its coarse moduli space. For each index $i \in [n]$, let $\pi_i : \overline{\mathcal{M}}_{g,n} \to \overline{\mathcal{M}}_{g,n-1}$ and $\pi_i : \M{g}{n} \to \M{g}{n-1}$ be the projection maps that forget the $i$th marked point. Similarly, for any subset $S \subseteq [n]$, let $\pi_S : \overline{\mathcal{M}}_{g,n} \to \overline{\mathcal{M}}_{g,|S|}$ and $\pi_S : \M{g}{n} \to \M{g}{|S|}$ be the projection maps that forget all marked points not in $S$.

The notation $\sum_{i \stackrel{\text{2}}{=} n}^{m} a_i$ denotes the summation of $a_i$ over integers $i$ between $n$ and $m$ that have the same parity as $n$.

\subsection{F-curves and the F-conjecture}\label{subsec:F}

Let $\M{g}{n}^r$ be the closed subvariety of $\M{g}{n}$ parametrizing stable curves with at least $r$ nodes. The codimension of $\M{g}{n}^r$ is $r$ and its irreducible components are called \textbf{codimension $r$ boundary strata}. When $r = 3g - 4 + n$, the space $\M{g}{n}^r$ is one-dimensional, and its irreducible components are referred to as \textbf{F-curves}.

F-curves are divided into six types based on the topological structure of the curves they parametrize. Type 1 F-curves are given by a map $f: \M{1}{1} \to \M{g}{n}$, which attaches a fixed stable genus $g-1$ curve with $n+1$ marked points. All other types of F-curves correspond to maps $f: \M{0}{4} \to \M{g}{n}$, attaching appropriate stable curves at specific points. In particular, Type 6 F-curves correspond to a partition $[n] = I_1 \sqcup I_2 \sqcup I_3 \sqcup I_4$ and nonnegative integers $g = g_1 + g_2 + g_3 + g_4$ such that $g_i > 0$ if $I_i = \emptyset$. The corresponding map $f: \M{0}{4} \to \M{g}{n}$ attaches a fixed prestable curve of genus $g_i$ with $|I_i|+1$ marked points to the $i$th point. Moreover, Type 5 F-curves correspond to a partition $[n] = I_1 \sqcup I_2 \sqcup I_3$ and nonnegative integers $g-1 = g_1 + g_2 + g_3$ such that $g_i > 0$ if $I_i = \emptyset$ for $i = 1,2$. The corresponding map $f: \M{0}{4} \to \M{g}{n}$ attaches a fixed prestable curve of genus $g_i$ with $|I_i|+1$ marked points to the first and second points, and a fixed prestable curve of genus $g_3$ with $|I_3|+2$ marked points to the third and fourth points. We will denote the corresponding F-curves of Type 5 and Type 6 by
\[ F_5(g_1,g_2,I_1,I_2),\ F_6(g_1,g_2,g_3,g_4, I_1,I_2,I_3,I_4). \]
If $g=0$, then we let $F_{I_1, I_2, I_3, I_4}:= F_6(0,0,0,0, I_1, I_2, I_3, I_4)$ for convenience. We refer to \cite[Section 2]{GKM02} for a detailed explanation of other types of F-curves.

A line bundle is \textbf{F-nef} (resp. \textbf{F-ample}) if its intersection with any F-curve is nonnegative (resp. positive). By definition, F-nefness (resp. F-ampleness) is weaker than nefness (resp. ampleness). However, it has the advantage of being straightforward to verify; see \cite[Theorem 2.1]{GKM02}. The following conjecture asserts that they are equivalent, thereby completely determining the F-cone of $\M{g}{n}$.

\begin{conj}\label{conj:F}[F-conjecture, \cite{KM13, GKM02}]
A line bundle on $\M{g}{n}$ is nef if and only if it is F-nef. Equivalently, it is ample if and only if it is F-ample.
\end{conj}

As mentioned in the introduction, there are many previous works on \cref{conj:F}, including \cite{KM13, GKM02, GF03, Gib09, Lar11, Fe15, MS19, Fe20}. In particular, \cite{KM13} proves the F-conjecture for $\M{0}{n}$ with $n \leq 7$ when the characteristic of the base field is $0$, while \cite{Lar11} establishes the same result for characteristic $p$. Furthermore, \cite{GKM02} proves the following result, thereby reducing the F-conjecture to the genus $0$ case.

\begin{thm}(\cite[Theorem 0.3]{GKM02})\label{thm:fconj0}
    Let $F:\M{0}{g+n}\to \M{g}{n}$ be the flag map defined in \cite[Section 0]{GKM02}. If $L$ is a F-nef divisor on $\M{g}{n}$, then  $L$ is nef if and only if $F^\ast L$ is nef.
\end{thm} 

\begin{thm}(\cite[Corollary 0.4]{GKM02})\label{cor:fconj}
    If $g+n\le 7$, then the F-conjecture holds for $\M{g}{n}$.
\end{thm} 

Note that even the F-conjecture for $\M{0}{8}$ remains open.

\subsection{Vertex Operator Algebras and Virasoro VOAs}\label{subsec:VOA}

The vertex operator algebra (VOA) is a mathematical framework used for the rigorous study of $2$-dimensional conformal field theory, first introduced by Borcherds \cite{Bor86}. Since the precise definitions and concepts in the theory of VOAs are somewhat intricate, we will provide only a heuristic description here. We refer the reader to \cite{FLM88, LL04, MT10} for a detailed and formal treatment.

\begin{defthm}\label{defn:VOA}
    A \textbf{vertex operator algebra} (VOA for short) is a graded $\C$-algebra $V=\oplus_{n\in \N}V_n$ with infinitely many products, indexed by $\Z$, and an action of the Virasoro Lie algebra of the central charge $c$. These products admit corresponding associativity and commutativity relations.

    A \textbf{module} over a vertex operator algebra $V$ is a graded $\C$-vector space $M=\oplus_{n\in \N}M_n$ with an action of $V$ and a compatible action of the Virasoro Lie algebra. If $L_0$ of the Virasoro Lie algebra acts on $M_n$ by $h+n$ for some $h\in \C$, then $h$ is called the \textbf{conformal weight} of $M$. If $M$ is a simple $V$-module, then $M$ admits a conformal weight $h$. 

    $V$ is \textbf{of CFT type} if $V_0=\C\textbf{1}$. $V$ is \textbf{rational} if any $V$-module is semisimple, i.e., a sum of simple $V$-modules, and there are only finitely many isomorphism classes of simple $V$-modules. $V$ is $\mathbf{C_2}$\textbf{-cofinite} if $C_2(V):=\text{Span} \left\{A_{(-2)}B\ \mid\ A,B\in V \right\}$ is of finite codimension inside $V$.
\end{defthm}

Unless otherwise stated, every VOA in this paper is assumed to be of CFT type. The preceding definition of $V$-module is called the \textbf{admissible modules}. We refer the reader to \cite{MT10} for a detailed treatment. Moreover, there is a notion of tensor product of $V$-modules, which we will denote by $\boxtimes$ for historical reasons.

Now, we will consider two examples of VOAs: discrete series Virasoro VOAs and affine VOAs, which are important examples and play a significant role in this paper. We refer to \cite{FZ92, LL04} for more information. For any coprime integers $2\le p,q$, there exists a $C_2$-cofinite, rational VOA $\text{Vir}_{p,q}$, called a \textbf{discrete series Virasoro VOA}, whose central charge is
\[ c_{p,q}=1-\frac{6(p-q)^2}{pq}. \]

For any pair of integers $(r,s)$ such that $1\le r\le q-1$, $1\le s\le p-1$, there exists a simple module $W_{r,s}$ over $\text{Vir}_{p,q}$, whose conformal weight is
\[ h_{r,s}=\frac{(pr-qs)^2-(p-q)^2}{4pq}.\]

Any simple module over $\text{Vir}_{p,q}$ is isomorphic to one of them, and they are not isomorphic to each other except $W_{r,s}\simeq W_{p-r,q-s}$. Hence, there are $\frac{(p-1)(q-1)}{2}$ simple modules. If $q=2$ and $p=2k+1$, then $r$ is always $1$; hence, we will denote $W_{r,s}$ by $W_s$. Hence, we have $2k$ modules, $W_1,\cdots, W_{2k}$ such that $W_i\simeq W_{2k+1-i}$ and non-isomorphic otherwise. Note that $W_1=W_{2k}=\vir{2k+1}$ as $\vir{2k+1}$-modules. Hence, $\left\{W_1,\cdots,W_k\right\}, \left\{W_i\ |\ i \text{ even}\right\}$ and $\left\{W_i\ |\ i \text{ odd}\right\}$ are representatives of non-isomorphic irreducible modules. We will denote $c_{2k+1,2}$ by $c_k$ and $h_{1, i}$ by $h_i$. Explicitly, we have
\[ c_k=-\frac{2(k-1)(6k-1)}{2k+1}\text{ and } h_i=-\frac{(2k-i)(i-1)}{2(2k+1)}. \]

The fusion rules are given by
\[ W_{i_1}\boxtimes W_{i_2}=\sum_{i\stackrel{\text{2}}{=}|i_1-i_2|+1}^{\min(i_1+i_2, 4k+2-(i_1+i_2))-1}  W_{i}  \]

where $\stackrel{\text{2}}{=}$ denotes sums with an increment of $2$. Note that every simple module over $\vir{2k+1}$ is self-dual. If $k=1$, this construction gives a holomorphic vertex operator algebra with $c=0$, which is not our interest since its corresponding conformal block divisors are trivial (cf. \cref{thm:c1cdf}). Hence, we will only consider the case $k > 1$.

Another important class of VOAs relevant to this paper is the \textbf{affine vertex operator algebras} associated with a simple Lie algebra $\mathfrak{g}$. For any positive integer $k$, there exists a $C_2$-cofinite, rational VOA $L_k(\mathfrak{g})$, known as the level $k$ affine vertex operator algebra for $\mathfrak{g}$. Furthermore, the representation theory of $L_k(\mathfrak{g})$ is closely related to the representation theory of $\mathfrak{g}$.  It is worth noting that for any complex number $k$ not equal to the negation of the dual Coxeter number of $\mathfrak{g}$, we can also define the level $k$ affine vertex operator algebra for $\mathfrak{g}$. However, this general case is not $C_2$-cofinite or rational in general. In this paper, we restrict our attention to the case where $k$ is a positive integer.

We refer to \cite{FZ92, LL04, MT10} for a detailed description of affine vertex operator algebras and their representations and focus only on the representation theory of $L_1(\mathfrak{sl}_{r+1})$. Its representation theory is elementary because it is a lattice VOA. Specifically, the set of simple modules over $L_1(\mathfrak{sl}_{r+1})$ corresponds to elements of $\mathbb{Z}/(r+1)$ and the tensor product structure is given by the group law of $\mathbb{Z}/(r+1)$. More precisely, for any $0 \leq i < r+1$, if we let $U_i$ denote the simple module corresponding to $i$ over $L_1(\mathfrak{sl}_{r+1})$, then $U_i \boxtimes U_j = U_{i+j}$. Note that the central charge of $L_1(\mathfrak{sl}_{r+1})$ and the conformal weight of $U_i$ are given by
\[ c_r'=r\text{ and }h'_i=\frac{i(r+1-i)}{2(r+1)}. \]
 
\subsection{Conformal Block Divisors}\label{subsec:cdf}

The \textbf{sheaf of coinvariants} was introduced in \cite{TK88, TUY89} in the context of conformal field theory. This concept was further developed in \cite{BFM91, FBZ04, NT05, Uen08} and eventually extended to any $C_2$-cofinite and rational VOA over $\overline{\mathcal{M}}_{g,n}$ for any genus $g$ by \cite{DGT21, DGT22a, DGT22b}.

\begin{defthm}[\cite{DGT21, DGT22a, DGT22b}]\label{def:CBD}
Let $V$ be a $C_2$-cofinite and rational VOA, and let $W^1, \ldots, W^n$ be simple $V$-modules. There exists a vector bundle $\mathbb{V}_{g,n}(V, W^{\bullet})$ on $\overline{\mathcal{M}}_{g,n}$, referred to as the \textbf{sheaf of coinvariants}. The dual of this vector bundle, denoted by $\mathbb{V}_{g,n}(V, W^{\bullet})^\vee$, is called the \textbf{sheaf of conformal blocks}. These bundles' first Chern classes $\mathbb{D}_{g,n}(V, W^{\bullet})$, $-\mathbb{D}_{g,n}(V, W^{\bullet})$ on $\overline{\mathcal{M}}_{g,n}$ are referred to as the \textbf{coinvariant divisor} and the \textbf{conformal block divisor}.
\end{defthm}

Note that, as the $\Q$-Picard group of $\overline{\mathcal{M}}_{g,n}$ and $\M{g}{n}$ coincide, we will consider coinvariant divisors and conformal block divisors as $\Q$-divisors on $\M{g}{n}$.

\begin{rmk}\label{rmk:cdiv}
    In earlier documents, coinvariant divisors associated with affine VOAs were often referred to as conformal block divisors. For example, see \cite{TUY89, Fak12, BGM15, BGM16}. However, this terminology may be confusing, as they are, in fact, the first Chern class of the sheaf of coinvariants, not of the sheaf of conformal blocks, which is the dual of the sheaf of coinvariants. We propose adopting a new and more accurate term for these divisors to address this issue.
\end{rmk}

The two most interesting properties of the sheaf of coinvariants are the propagation of vacua and the factorization theorem.

\begin{thm}\label{thm:propvac}(Propagation of Vacua, \cite[Theorem 3.6]{Cod20}, \cite[Theorem 4.3.1]{DGT22a})
     Let $V$ be a $C_2$-cofinite and rational VOA. Then the following holds:
    \[ \pi_{n+1}^\ast\mathbb{V}_{g,n}\left(V,W^\bullet \right)=\mathbb{V}_{g,n+1}\left(V,W^\bullet\otimes V \right).  \]
\end{thm}

\begin{thm}\label{thm:factor}(Factorization Theorem, \cite[Theorem 6.2.6]{TUY89}, \cite[Theorem 7.0.1]{DGT22a})
    Let $V$ be a $C_2$-cofinite and rational VOA and $S$ be the set of isomorphism classes of simple $V$-modules. Also, let $\xi:\overline{\mathcal{M}}_{g_1,n_1}\times \overline{\mathcal{M}}_{g_2,n_2}\to \overline{\mathcal{M}}_{g,n}$ be the gluing map. Then, there exists a canonical isomorphism
    \[ \xi^\ast \mathbb{V}_{g,n}\left(V, W^{\bullet} \right) \simeq \bigoplus_{W \in S} \pi_1^\ast \mathbb{V}_{g_1,n_1}\left(V, W^{\bullet}\otimes W \right) \otimes \pi_2^\ast\mathbb{V}_{g_2,n_2}\left(V, W^{\bullet}\otimes W' \right). \]
    Similarly, for the gluing map $\xi_{irr}:\overline{\mathcal{M}}_{g-1,n+2}\to \overline{\mathcal{M}}_{g,n}$, we have a canonical isomorphism
    \[ \xi_{irr}^\ast \mathbb{V}_{g,n}\left(V, W^{\bullet} \right) \simeq \bigoplus_{W \in S}  \mathbb{V}_{g-1,n+2}\left(V, W^{\bullet}\otimes W \otimes W' \right). \]
    Here, $W'$ denotes the dual module. 
\end{thm}

Fortunately, we can explicitly express the coinvariant divisor in terms of tautological line bundles.

\begin{thm}\cite[Corollary 2]{DGT22b}\label{thm:c1cdf}
    Let $V$ be a rational, $C_2$-cofinite VOA with central charge $c$, and let $W^1,\cdots, W^n$ be simple modules with conformal weights $h_1,\cdots, h_n$. Let $S$ be the set of simple modules over $V$. Then
    \[ \mathbb{D}_{g,n}\left(V, W^\bullet\right)=\rank \mathbb{V}_{g,n}(V, W^\bullet)\left( \frac{c}{2}\lambda+\sum_{i=1}^n h_i\psi_i\right)-b_{\text{irr}}\delta_{irr}-\sum b_{i,I}\delta_{i,I} \]
    where
    \[ b_{\text{irr}}=\sum_{W\in S} h_W \cdot \rank \mathbb{V}_{g,n}(V, W^\bullet\otimes W\otimes W'), \]
    \[ b_{i,I}=\sum_{W\in S} h_W \cdot \rank \mathbb{V}_{i,|I|+1}(V, W^I\otimes W) \cdot \rank \mathbb{V}_{g-i,|I^c|+1}(V, W^{I^c}\otimes W').\]
\end{thm}

\section{Positivity of Conformal Block Divisors}\label{sec:pos}

\subsection{F-nefness of Conformal Block Divisors}\label{subsec:fnef}

The main result of this section is the following.

\begin{thm}\label{thm:virfnef}
    Any conformal block divisor of $V=\vir{2k+1}$ is F-nef. Furthermore, if $W^1\cdots, W^n$ be nontrivial simple $\vir{2k+1}$-modules, then the conformal block divisor $-\mathbb{D}_{g,n}\left(V, W^{\bullet}\right)$ is either zero or F-ample. More precisely, $-\mathbb{D}_{g,n}(V, W^\bullet)$ is trivial if and only if $\text{rank }\mathbb{V}_{g,n}(V, W^\bullet) \leq 1$, and it is F-ample otherwise. In particular, if $W^i$'s are nontrivial,
    \begin{enumerate}
        \item If $g\ge 2$, then $-\mathbb{D}_{g,n}\left(V, W^{\bullet}\right)$ is F-ample.
        \item If $g=1$, then $-\mathbb{D}_{g,n}\left(V, W^{\bullet}\right)$ is F-ample except in the special case where $n=1$ and $W^1=W_2$.
        \item If $g=0$ and $n>2k-1$, then $-\mathbb{D}_{g,n}\left(V, W^{\bullet}\right)$ is F-ample.
    \end{enumerate}
\end{thm}

To prove this theorem, we need the following lemma. We will also prove an analogous lemma for nefness together, which will be used in \cref{subsec:nef}.

\begin{lem}\label{lem:VOAGKM}
    Let $V$ be a $C_2$-cofinite and rational VOA. 
    \begin{enumerate}
        \item If conformal block divisors (resp. coinvariant divisor) of $V$ are nef  on $\overline{\rm{M}}_{0,4}$ and $\overline{\rm{M}}_{1,1}$, then any conformal block divisor (resp. coinvariant divisor) of $V$ on $\overline{\rm{M}}_{g,n}$ is F-nef.
        \item If conformal block divisors (resp. coinvariant divisor) of $V$ are nef on $\overline{\rm{M}}_{0,n}$ for all $n$ and $\overline{\rm{M}}_{1,1}$, then any conformal block divisor (resp. coinvariant divisor) of $V$ on $\overline{\rm{M}}_{g,n}$ is nef.
    \end{enumerate}
\end{lem}

\begin{proof}
    (1) The core of the proof is that, by factorization \cref{thm:factor}, the intersection of an F-curve and a conformal block divisor can be expressed in terms of the rank of the conformal blocks and the degree of the conformal block divisors on $\overline{\rm{M}}_{0,4}$. We have six cases corresponding to the six types of F-curves described in \cite[Section 2]{GKM02}. Since F-curves of type 2-5 are similar to F-curves of type 6, we will only describe two cases in detail: type 1 and type 6.

    The type 1 F-curve corresponds to a map $F:\M{1}{1}\to \M{g}{n}$, attaching a fixed stable curve of genus $g-1$ with $n+1$ points. Then, by \cref{thm:factor},
    \[ F^\ast \mathbb{D}_{g,n}(V, W^\bullet)=\sum_W \text{rank }\mathbb{V}_{g-1,n+1}(V, W^\bullet\otimes W) \cdot \mathbb{D}_{1,1}(V, W') \]
    where $W$ is runs through simple $V$-modules. Note that $\text{rank }\mathbb{V}_{g-1,n+1}(V, W^\bullet\otimes W)$ is nonnegative, hence if the degree of $\mathbb{D}_{1,1}(V, W')$ are nonnegative (resp. nonpositive), then the degree of $F^\ast \mathbb{D}_{g,n}(V, W^\bullet)$ is also nonnegative (resp. nonpositive).

    As mentioned in \cref{subsec:F}, the type 6 F-curves correspond to a pair consisting of a partition $[n]=I_1\sqcup I_2\sqcup I_3\sqcup I_4$ and nonnegative integers $g=g_1+g_2+g_3+g_4$. There exists an associated map $F: \M{0}{4} \to \M{g}{n}$, attaching a fixed stable curve of genus $g_i$ with $|I_i| + 1$ points for each $i = 1, 2, 3, 4$. Then, by \cref{thm:factor},
    \[ F^\ast \mathbb{D}_{g,n}(V, W^\bullet)=\sum_{W_1, W_2, W_3, W_4} \left( \prod_{i=1}^{4}\text{rank }\mathbb{V}_{g_i,|I_i|+1}(V, W^{I_i}\otimes W_i)\right) \cdot \mathbb{D}_{0,4}(V, \otimes_{i=1}^4 W_i'). \]
    The rest of the proof is the same as the type 1 curve.

    (2) By \cref{thm:fconj0} and (1), it is enough to show that $F^\ast \mathbb{D}_{g,n}(V, W^\bullet)$ is nef. By \cref{thm:factor},
    \[ F^\ast\mathbb{D}_{g,n}(V, W^\bullet)=\sum_{W_1,\cdots, W_n} \left( \prod_{i=1}^{n}\text{rank }\mathbb{V}_{1,1}(V, W_i)\right) \cdot \mathbb{D}_{0,g+n}(V, W^\bullet\otimes\left(\otimes_{i=1}^n W_i')\right). \]
    Since $\text{rank }\mathbb{V}_{1,1}(V, W_i)$ is nonnegative for each $i$, $F^\ast\mathbb{D}_{g,n}(V, W^\bullet)$ is nef by the assumption.
\end{proof}

\begin{thm}\label{thm:virdeg}
\begin{enumerate}
    \item Let $r=\rank \mathbb{V}_{0,4}\left(\vir{2k+1}, W^\bullet \right)$. Then
    \[ \deg \mathbb{V}_{0,4}\left(\vir{2k+1}, W^\bullet  \right) =-\frac{r(r-1)}{2}. \]
    \item For $0\le i\le k-1$, let  $r=\rank \mathbb{V}_{1,1}\left(\vir{2k+1}, W_{2i+1}  \right)$. Then $r=k-i$ and 
    \[ \deg \mathbb{V}_{1,1}\left(\vir{2k+1}, W_{2i+1} \right) =-r(r-1). \]
\end{enumerate}
\end{thm}

\begin{proof}
    (1) Let $W^i=W_{a_i}$, where $1\le a_i\le k$. We may assume that $a_1\le a_2\le a_3\le a_4$. Following the notation of \cref{thm:c1cdf}, it is enough to compute $r$ and $b_{0, \left\{1,2\right\}}, b_{0, \left\{1,3\right\}}, b_{0, \left\{1,4\right\}}$. 

    \textbf{Case 1. } $2\mid a_1+a_2+a_3+a_4$ and $a_1+a_4\ge a_2+a_3$.
    By the fusion rules,
    \[ W_{a_1}\boxtimes W_{a_2} = \sum_{i\stackrel{\text{2}}{=}a_2-a_1+1}^{a_1+a_2-1}  W_{i},\ W_{a_3}\boxtimes W_{a_4} = \sum_{i\stackrel{\text{2}}{=}a_4-a_3+1}^{a_3+a_4-1}  W_{i}.\]
    By the factorization theorem, to compute the rank, it is enough to compute the number of common summands of right-hand sides. Since both RHSs consist solely of even or odd modules, and $a_1+a_4\ge a_2+a_3$, the overlapping modules are 
    \[ \left\{W_i\ |\ a_4-a_3+1\le i\le a_1+a_2-1, i\equiv a_4-a_3+1 \text{ mod }2 \right\}. \]
    If $a_1+a_2-1<a_4-a_3+1$, the rank is zero, so (1) is trivially true. Hence, we may assume that $a_1+a_2-1\ge a_4-a_3+1$. In this case,
    \[ r=\frac{a_1+a_2+a_3-a_4}{2}. \]
    By definition, 
    \begin{align*}
        -b_{0, \left\{1,2\right\}}&=-\sum_{i\stackrel{\text{2}}{=}a_4-a_3+1}^{a_1+a_2-1}h_i=\sum_{i\stackrel{\text{2}}{=}a_4-a_3+1}^{a_1+a_2-1}\frac{(2k-i)(i-1)}{2(2k+1)}=\sum_{i\stackrel{\text{2}}{=}a_4-a_3+1}^{a_1+a_2-1}\frac{i}{2}-\frac{1}{2(2k+1)}\sum_{i\stackrel{\text{2}}{=}a_4-a_3+1}^{a_1+a_2-1}i^2-\frac{rk}{2k+1}\\
        &=\frac{r(r+a_4-a_3)}{2}-\frac{(a_4-a_3+1)^2r}{2(2k+1)}-\frac{(a_4-a_3+1)r(r-1)}{2k+1}-\frac{(r-1)r(2r-1)}{3(2k+1)}-\frac{rk}{2k+1}.
    \end{align*}
    By the same process,
    \[ -b_{0, \left\{1,3\right\}}=\frac{r(r+a_4-a_2)}{2}-\frac{(a_4-a_2+1)^2r}{2(2k+1)}-\frac{(a_4-a_2+1)r(r-1)}{2k+1}-\frac{(r-1)r(2r-1)}{3(2k+1)}-\frac{rk}{2k+1} \]
    and
    \[  -b_{0, \left\{1,4\right\}}=\frac{r(r+a_4-a_1)}{2}-\frac{(a_4-a_1+1)^2r}{2(2k+1)}-\frac{(a_4-a_1+1)r(r-1)}{2k+1}-\frac{(r-1)r(2r-1)}{3(2k+1)}-\frac{rk}{2k+1}. \]
    Each representation consists of five terms. We will compute $-b_{0, \left\{1,2\right\}}-b_{0, \left\{1,3\right\}}-b_{0, \left\{1,4\right\}}$. The sum of the first terms is
    \[ \frac{r(3r+3a_4-a_1-a_2-a_3)}{2}=\frac{r(r+2a_4)}{2}. \]
    The sum of the second terms is the most complicated one. It is
    \begin{align*}
        \sum_{i=1}^{3}\frac{(a_4-a_i+1)^2r}{2(2k+1)}&=\frac{r}{2(2k+1)}\left(3+2(3a_4-a_1-a_2-a_3)+3a_4^2+a_1^2+a_2^2+a_3^2-2a_4(a_1+a_2+a_3)\right)\\
        &=\frac{r(3+4(a_4-r)-4ra_4)}{2(2k+1)}+\frac{r(a_1^2+a_2^2+a_3^2+a_4^2)}{2(2k+1)}.
    \end{align*}
    The sum of the third terms is
    \[ \frac{r(r-1)}{2k+1}\left(3a_4-a_1-a_2-a_3+3\right)=\frac{r(r-1)}{2k+1}\left(2a_4-2r+3\right) \]
    The sum of the fourth and fifth terms is, by inspection
    \[ \frac{r(r-1)(2r-1)}{2k+1}+\frac{3rk}{2k+1}. \]
    By summing the third, fourth and fifth terms, we obtain
    \[ \frac{r(r-1)(2a_4+2)}{2k+1}+\frac{3rk}{2k+1}. \]
    Adding this with the second term gives
    \[ -\frac{r}{2(2k+1)}+\frac{3rk}{2k+1}+\frac{r(a_1^2+a_2^2+a_3^2+a_4^2)}{2(2k+1)}. \]
    All in all, we have
    \[ -b_{0, \left\{1,2\right\}}-b_{0, \left\{1,3\right\}}-b_{0, \left\{1,4\right\}}=\frac{r(r+2a_4)}{2}+\frac{r}{2(2k+1)}-\frac{3rk}{2k+1}-\frac{r(a_1^2+a_2^2+a_3^2+a_4^2)}{2(2k+1)}. \]
    On the other hand,
    \[r\left(h_{a_1}+h_{a_2}+h_{a_3}+h_{a_4} \right)=-r\frac{a_1+a_2+a_3+a_4}{2}+\frac{r(a_1^2+a_2^2+a_3^2+a_4^2)}{2(2k+1)}+\frac{4rk}{2k+1}.  \]
    Hence, by \cref{thm:c1cdf}, the degree is
    \[ -\frac{r^2}{2}+\frac{r}{2(2k+1)}+\frac{rk}{2k+1}=\frac{r(1-r)}{2}. \]

    \textbf{Case 2. } $2\mid a_1+a_2+a_3+a_4$ and $a_1+a_4\le a_2+a_3$.
    By the same approach as above, we can obtain that the set of common simple modules in the decompositions for $W_{a_1}\boxtimes W_{a_2}$ and $W_{a_3}\boxtimes W_{a_4}$ is 
    \[ \left\{W_i\ |\ a_2-a_1+1\le i\le a_1+a_2-1, i\equiv a_2-a_1+1 \text{ mod }2. \right\}. \]
    Hence, in particular, $r=a_1$. Moreover,
    \begin{align*}
        -b_{0, \left\{1,2\right\}}&=-\sum_{i\stackrel{\text{2}}{=}a_2-a_1+1}^{a_1+a_2-1}h_i=\sum_{i\stackrel{\text{2}}{=}a_2-a_1+1}^{a_1+a_2-1}\frac{(2k-i)(i-1)}{2(2k+1)}=\sum_{i\stackrel{\text{2}}{=}a_2-a_1+1}^{a_1+a_2-1}\frac{i}{2}-\frac{1}{2(2k+1)}\sum_{i\stackrel{\text{2}}{=}a_2-a_1+1}^{a_1+a_2-1}i^2-\frac{rk}{2k+1}\\
        &=\frac{r(r+a_2-a_1)}{2}-\frac{(a_2-a_1+1)^2r}{2(2k+1)}-\frac{(a_2-a_1+1)r(r-1)}{2k+1}-\frac{(r-1)r(2r-1)}{3(2k+1)}-\frac{rk}{2k+1}.
    \end{align*}
    Similarly,
    \[  -b_{0, \left\{1,j\right\}}=\frac{r(r+a_j-a_1)}{2}-\frac{(a_j-a_1+1)^2r}{2(2k+1)}-\frac{(a_j-a_1+1)r(r-1)}{2k+1}-\frac{(r-1)r(2r-1)}{3(2k+1)}-\frac{rk}{2k+1} \]
    for $j=2,3,4$. Now, it is evident that the remainder of the proof is the same as in Case 1.

    \textbf{Case 3. }$2\nmid a_1+a_2+a_3+a_4$ and $2k+1\ge a_4+a_3+a_2-a_1$. By the fusion rules,
    \[ W_{a_1}\boxtimes W_{a_2} = \sum_{i\stackrel{\text{2}}{=}a_2-a_1+1}^{a_1+a_2-1}  W_{i},\ W_{a_3}\boxtimes W_{a_4} = \sum_{i\stackrel{\text{2}}{=}a_4-a_3+1}^{a_3+a_4-1}  W_{i}=\sum_{i\stackrel{\text{2}}{=}2k+2-(a_3+a_4)}^{2k-a_4+a_3}  W_{i}.\]
    Hence, by the same argument as in Case 1, $r=0$ if $a_1+a_2+a_3+a_4<2k+3$. Therefore, we may assume that $a_1+a_2+a_3+a_4\ge 2k+3$. In this case,
    \[ r=\frac{a_1+a_2+a_3+a_4-(2k+1)}{2}. \]
    Also, by the definition of $b_{0, \left\{1,2\right\}}$,
    \begin{align*}
        -b_{0, \left\{1,2\right\}}&=-\sum_{i\stackrel{\text{2}}{=}2k+2-(a_3+a_4)}^{a_1+a_2-1}h_i=\sum_{i\stackrel{\text{2}}{=}2k+2-(a_3+a_4)}^{a_1+a_2-1}\frac{i}{2}-\frac{1}{2(2k+1)}\sum_{i\stackrel{\text{2}}{=}2k+2-(a_3+a_4)}^{a_1+a_2-1}i^2-\frac{rk}{2k+1}\\
        &=\frac{r(a_1+a_2-r)}{2}-\frac{(a_1+a_2-1)^2r}{2(2k+1)}+\frac{(a_1+a_2-1)r(r-1)}{2k+1}-\frac{(r-1)r(2r-1)}{3(2k+1)}-\frac{rk}{2k+1}.
    \end{align*}
    Hence, by the same computation, 
    \[  -b_{0, \left\{1,j\right\}}=\frac{r(a_1+a_j-r)}{2}-\frac{(a_1+a_j-1)^2r}{2(2k+1)}+\frac{(a_1+a_j-1)r(r-1)}{2k+1}-\frac{(r-1)r(2r-1)}{3(2k+1)}-\frac{rk}{2k+1} \]
    for $j=2,3,4$.
    
    The sum of the third terms is
    \[ \frac{r(r-1)}{2k+1}\left(3a_1+a_2+a_3+a_4-3\right)=\frac{r(r-1)}{2k+1}\left(2a_1+2r+2k-2\right). \]
    By summing the third, fourth, and fifth terms, we obtain
    \[ \frac{r(r-1)}{2k+1}\left(2a_1+2k-1\right)-\frac{3rk}{2k+1}. \]
    The sum of the second terms is
    \begin{align*}
        &-\frac{3r}{2(2k+1)}+\frac{r}{2k+1}(3a_1+a_2+a_3+a_4)-\frac{ra_1}{2k+1}(a_2+a_3+a_4)-\frac{r}{2(2k+1)}(3a_1^2+a_2^2+a_3^2+a_4^2)\\
        &=-\frac{3r}{2(2k+1)}+\frac{2r}{2k+1}(r+a_1)+r-\frac{ra_1}{2k+1}(2r-a_1)-ra_1-\frac{r}{2(2k+1)}(a_1^2+a_2^2+a_3^2+a_4^2)-\frac{ra_1^2}{2k+1}.
    \end{align*}
    Hence, the summation of everything except the first term is
    \[ r(r-1)-\frac{r}{2}+\frac{2r}{2k+1}-ra_1-\frac{r}{2(2k+1)}(a_1^2+a_2^2+a_3^2+a_4^2).\]
    Consequently, the total summation is
    \[ r(r-1)-\frac{r^2+r}{2}+\frac{2r}{2k+1}-\frac{r}{2(2k+1)}(a_1^2+a_2^2+a_3^2+a_4^2)+\frac{r(2k+1)}{2}.\]
    On the other hand,
    \[r\left(h_{a_1}+h_{a_2}+h_{a_3}+h_{a_4} \right)=-r\frac{a_1+a_2+a_3+a_4}{2}+\frac{r(a_1^2+a_2^2+a_3^2+a_4^2)}{2(2k+1)}+\frac{4rk}{2k+1}.  \]
    Hence, the degree is
    \begin{align*}
         &r(r-1)-\frac{r^2+r}{2}+\frac{2r}{2k+1}+\frac{r(2k+1)}{2}-r\frac{a_1+a_2+a_3+a_4}{2}+\frac{4rk}{2k+1}\\
         &=\frac{r^2+r}{2}+\frac{r(2k+1)}{2}-r\frac{2r+2k+1}{2}=-\frac{r(r-1)}{2}.
    \end{align*}

    \textbf{Case 4. }$2\nmid a_1+a_2+a_3+a_4$ and $2k+1<a_4+a_3+a_2-a_1$. Since $a_1-a_2+a_3+a_4\le a_3+a_4 \le 2k$, as in Step 3, we have $r=a_1$. Moreover,
    \begin{align*}
        -b_{0, \left\{1,2\right\}}&=-\sum_{i\stackrel{\text{2}}{=}a_2-a_1+1}^{a_1+a_2-1}h_i=\sum_{i\stackrel{\text{2}}{=}a_2-a_1+1}^{a_1+a_2-1}\frac{(2k-i)(i-1)}{2(2k+1)}=\sum_{i\stackrel{\text{2}}{=}a_2-a_1+1}^{a_1+a_2-1}\frac{i}{2}-\frac{1}{2(2k+1)}\sum_{i\stackrel{\text{2}}{=}a_2-a_1+1}^{a_1+a_2-1}i^2-\frac{rk}{2k+1}\\
        &=\frac{r(r+a_2-a_1)}{2}-\frac{(a_2-a_1+1)^2r}{2(2k+1)}-\frac{(a_2-a_1+1)r(r-1)}{2k+1}-\frac{(r-1)r(2r-1)}{3(2k+1)}-\frac{rk}{2k+1}
    \end{align*}
     and
    \[  -b_{0, \left\{1,j\right\}}=\frac{r(r+a_j-a_1)}{2}-\frac{(a_j-a_1+1)^2r}{2(2k+1)}-\frac{(a_j-a_1+1)r(r-1)}{2k+1}-\frac{(r-1)r(2r-1)}{3(2k+1)}-\frac{rk}{2k+1} \]
    for $j=2,3,4$. Hence, the computation for this case is exactly the same as that of Case 2.
    
    (2) It is enough to prove this for $W=W_{2i+1}$. $0\le i\le k-1$. By the fusion rules, for $1\le j\le k$, $\rank \mathbb{V}_{0,3}\left(\vir{2k+1}, W_j \otimes W_j'\otimes W_{2i+1} \right)=1$ if $j>i$ and $0$ otherwise. Therefore, by the factorization theorem, $\rank\mathbb{V}_{1,1}\left(\vir{2k+1},  W_{2i+1} \right)=k-i$. Moreover, by \cref{thm:c1cdf},
    \[ -\mathbb{D}_{1,1}\left(V, W_{2i+1}\right)=\left[(k-i)\left(\frac{(k-1)(6k-1)}{2k+1}+\frac{i(2k-2i-1)}{2k+1}\right)-12\sum_{j=i+1}^{k}\frac{(2k-j)(j-1)}{2(2k+1)}\right]\lambda \]
    since $\delta_{\text{irr}}=12\lambda$ and $\lambda=\psi_1$. We can manually compute the second term to obtain
    \[ 12\sum_{j=i+1}^{k}\frac{(2k-j)(j-1)}{2(2k+1)}=\frac{2(k-i)(2k^2+2ki-i^2-3k+1)}{2k+1}. \]
    Since
    \[ (k-1)(6k-1)+i(2k-2i-1)-2(2k^2+2ki-i^2-3k+1)=(2k+1)(k-i-1), \]
    we have
    \[ -\mathbb{D}_{1,1}\left(V, W_{2i+1}\right)=(k-i)(k-i-1)\lambda. \]
    This proves the theorem.
\end{proof}

\begin{rmk}\label{rmk:rank}
    Note that we have demonstrated that the degree of the sheaf of conformal blocks of $\vir{2k+1}$ on $\M{0}{4}$ and $\M{1}{1}$ is uniquely determined by its rank. The author is unaware of any other non-holomorphic vertex operator algebras with this property. For instance, affine VOAs \cite[Proposition 4.2, Lemma 5.1]{Fak12}, lattice VOAs \cite[Section 9]{DG23}, and other discrete series Virasoro VOAs \cite[Section 8]{DG23} do not exhibit this behavior. Using \cite[Virasoro general]{Choigit24}, this can be explicitly verified for other discrete series Virasoro VOAs.

    In the remainder of this paper, particularly in \cref{subsec:crit}, we will leverage this unique property to uncover interesting characteristics of the conformal block divisors associated with $\vir{2k+1}$.
\end{rmk}

Using \cref{lem:VOAGKM}, \cref{thm:virdeg} proves the first part of \cref{thm:virfnef}. Establishing the latter part, however, requires a detailed analysis of the fusion rules of $\vir{2k+1}$. From the fusion rules, we observe the following: (1) The tensor product of two simple modules is the sum of $W_i$'s, where $i$ ranges over an arithmetic progression with a common difference of $2$. By further experimentation, e.g.over $\vir{9}$,
\[ W_2\boxtimes W_3\boxtimes W_5=W_2+2W_4+2W_6+W_8, \]
we note that (1) also applies to the tensor product of more than two modules. Additionally, we observe that (2) the coefficients of modules, except for the smallest and largest ones, are $\ge 2$, and (3) the length of the expression increases. These observations are rigorously summarized in the following lemma.

\begin{lem}\label{lem:virfus}
    Let $W^1\cdots, W^n$ be nontrivial simple $\vir{2k+1}$-modules. Then there exist integers $1\le a<b\le 2k$ of the same parity such that
    \[ \boxtimes_{i=1}^n W^i=\sum_{j\stackrel{\text{2}}{=}a}^b n_jW_j \]
    where
    \begin{enumerate}
        \item $n_j\ge 1$.
        \item If $n\ge 3$, then $n_j\ge 2$ for $j\ne a,b$.
        \item If $k\ge 3$, then $b-a\ge \min\left(2k-2, 2\lceil \frac{n-1}{2}\rceil \right)$.
    \end{enumerate}
\end{lem}

\begin{proof}
    Use induction. For $n=2$, this is evident by the fusion rules of $\vir{2k+1}$. We will prove the following: If $W=\sum_{j\stackrel{\text{2}}{=}a}^b n_jW_j$, $b\ne a$ are of the same parity and $n_j\ge 1$,
    
    \begin{enumerate}[label=(\alph*)] 
    \item There exist $1\le c<d\le 2k$ of the same parity such that $W\boxtimes W_i=\sum_{j\stackrel{\text{2}}{=}c}^d m_jW_j$ and $m_j\ge 1$. 
    \item If $W_i\ne V$, then $m_j\ge 2$ for $j\ne c,d$. 
    \item If $W_i\ne V, W_2$ and $b-a<2k-2$, then $d-c\ge b-a+2$. 
    \end{enumerate}

    By the induction hypothesis, (a) and (b) imply (1) and (2). For (3), we use the following form of induction: assume (3) holds for $n$ and $n+1$ and prove it for $n+2$ and $n+3$. If $k\ge 3$, for any nontrivial simple modules $W^1$ and $W^2$, $W^1\boxtimes W^2$ contains a nontrivial simple module other than $W_2$ as a summand. Therefore, by (c) and the induction hypothesis, (3) follows.

    Without loss of generality, assume $a,b$ are odd, and $i$ is even.
    
    \textbf{Case 1. }$i<a$. Then
    \[ W\boxtimes W_i=\sum_{j\stackrel{\text{2}}{=}a}^b n_j \sum_{r=j-i+1}^{\min(i+j, 4k+2-i-j)-1} W_r\ \]
    
    \textbf{Subcase 1.1. }$i+b<2k+1$. In this case, $W\boxtimes W_i=\sum_{j\stackrel{\text{2}}{=}a-i+1}^{b+i-1}m_jW_j$, and $m_j\ge 1$. Moreover, in the preceding expression, all terms except $c=a-i+1$ and $d=b+i-1$ appear more than twice. Hence, $m_j\ge 2$ for $j\ne c,d$. Finally, $d-c=b-a+2i-2\ge b-a+2$ if $i\ge 2$.

    \textbf{Subcase 1.2. }$i+b\ge 2k+1$. In this case, $W\boxtimes W_i=\sum_{j\stackrel{\text{2}}{=}a-i+1}^{2k}m_jW_j$, and $m_j\ge 1$. Moreover, in the preceding expression, all terms except $c=a-i+1$ and $d=2k$ appear more than twice. Hence, $m_j\ge 2$ for $j\ne c,d$. Finally, $d-c=2k-(a-i+1)\ge b-a+2$ if $i\ge 2$ (note that $b$ is odd, so $b\le 2k-1$).

    \textbf{Case 2. }$i>b$. Then
    \[ W\boxtimes W_i=\sum_{j\stackrel{\text{2}}{=}a}^b n_j \sum_{r=i-j+1}^{\min(i+j, 4k+2-i-j)-1} W_r\ \]

    \textbf{Subcase 2.1. } $i+b<2k+1$. In this case, $W\boxtimes W_i=\sum_{j\stackrel{\text{2}}{=}i-b+1}^{i+b-1}m_jW_j$, and $m_j\ge 1$. Moreover, in the preceding expression, all terms except $c=i-b+1$ and $d=i+b-1$ appear more than twice. Hence, $m_j\ge 2$ for $j\ne c,d$. Finally, $d-c=2b-2\ge b-a+2$. Note that $a,b$ are distinct odd natural numbers, so $a+b\ge 4$.

    \textbf{Subcase 2.2. }$i+b\ge 2k+1$. In this case, $W\boxtimes W_i=\sum_{j\stackrel{\text{2}}{=}i-b+1}^{2k}m_jW_j$, and $m_j\ge 1$. Moreover, in the preceding expression, all terms except $c=i-b+1$ and $d=2k$ appear more than twice. Hence, $m_j\ge 2$ for $j\ne c,d$. Finally, $d-c=2k-(i-b+1)\ge b-a+2$. Note that $W_i\ne V$ and $i$ is even, so $i\le 2k-2$.

    \textbf{Case 3.} $a<i<b$. Then
    \[  W\boxtimes W_i =\left(\sum_{j\stackrel{\text{2}}{=}a}^{i-1}\sum_{r=i-j+1}^{\min(i+j, 4k+2-i-j)-1} W_r\right)+\left(\sum_{j\stackrel{\text{2}}{=}i+1}^{b}\sum_{r=j-i+1}^{\min(i+j, 4k+2-i-j)-1} W_r \right). \]

    \textbf{Subcase 3.1. } $b+i\le 2k+1$. In this case, $W\boxtimes W_i=\sum_{j\stackrel{\text{2}}{=}2}^{i+b-1}m_jW_j$, and $m_j\ge 1$. Then $c=2$ and $d=i+b-1$. Moreover, all terms except $d$ appear more than twice in the preceding expression. Hence, $m_j\ge 2$ for $j\ne d$. Finally, $d-c=i+b-3\ge b-a+2$. Note that $W_i\ne V, W_2$, and $i$ is even, so $i\ge 4$.
    
    \textbf{Subcase 3.2. }$i+a\ge 2k+1$. In this case, $W\boxtimes W_i=\sum_{j\stackrel{\text{2}}{=}2}^{4k+1-i-a}m_jW_j$, and $m_j\ge 1$. Then $c=2$ and $d=4k+1-i-a$. Moreover, all terms except $d$ appear more than twice in the preceding expression. Hence, $m_j\ge 2$ for $j\ne d$. Finally, $d-c=4k-1-i-a\ge b-a+2$. Note that since $b$ is an odd number less than $2k$, $i+b\le 2b-1\le 4k-3$.

    \textbf{Subcase 3.3. }$i+a<2k+1<i+b$ and $i\ge k+1$. In this case, $W\boxtimes W_i=\sum_{j\stackrel{\text{2}}{=}2}^{2k}m_jW_j$, and $m_j\ge 1$. Then $c=2$ and $d=2k$. Moreover, all terms except $d$ appear more than twice in the preceding expression. Hence, $m_j\ge 2$ for $j\ne d$. Finally, $d-c=2k-2\ge b-a+2$ since $b-a<2k-2$.

    \textbf{Subcase 3.4. }$i+a<2k+1<i+b$ and $i<k+1$. In this case, $W\boxtimes W_i=\sum_{j\stackrel{\text{2}}{=}2}^{2k}m_jW_j$, and $m_j\ge 1$. Then $c=2$ and $d=2k$. Moreover, all terms except $d$ appear more than twice in the preceding expression. Hence, $m_j\ge 2$ for $j\ne d$. Finally, $d-c=2k-2\ge b-a+2$ since $b-a<2k$.
\end{proof}

\begin{lem}\label{lem:virdim}
    Let $V=\vir{2k+1}$ and $W^1,\cdots, W^n$ be simple $V$-modules.
    \begin{enumerate}
    \item $\rank \mathbb{V}_{g,n}(V, W^\bullet) \ge 2$ for $g\ge 2$. 
    \item $\rank \mathbb{V}_{1,n}(V, W^\bullet) \ge 2$ except in the special case where all but one module is trivial and that one module is $W_2$. In the latter case, $\rank \mathbb{V}_{1,n}(V, W^\bullet) = 1$.  Indeed, $\rank \mathbb{V}_{1,1}(V, W_{2i}) = i$.
    \item $\rank \mathbb{V}_{0,n}(V, W^\bullet \otimes W_k \otimes W_k) \ge 1$ for $n \ge 3$. If the modules $W^i$'s are nontrivial, $\rank \mathbb{V}_{0,n}(V, W^\bullet \otimes W_k \otimes W_k) \ge 2$ for $n \ge 4$. 
    \end{enumerate} 
\end{lem}

\begin{proof}
    The core of the proof is the identity $W_k \boxtimes W_k = \sum_{i=1}^{k} W_i$. Since all the modules $W^i$ are nontrivial, (3) holds by the fusion rules. By the \cref{thm:propvac}, we may also assume that $W^i$'s are nontrivial for (1). Hence, (1) is a consequence of (3). For (2), note that $\rank \mathbb{V}_{1,1}(V, W_{2i}) = i$ for $1 \le i \le k$ by \cref{thm:factor} and $W_r \boxtimes W_r = \sum_{i=1}^{r} W_{2i}$. Hence, by \cref{thm:propvac}, it is enough to consider the case that at least two of $W^1, \cdots, W^n$ are nontrivial. Again, this case follows from (3).
\end{proof}

\begin{proof}[proof of \cref{thm:virfnef}]
    If $\text{rank }\mathbb{V}_{g,n}(V, W^\bullet) \leq 1$, then by \cref{thm:virdeg} and \cref{thm:factor}, it is straightforward to conclude that $\mathbb{D}_{g,n}(V, W^\bullet) = 0$. Indeed, it suffices to show that the intersection of any F-curve with $\mathbb{D}_{g,n}(V, W^\bullet)$ is trivial. Using the same method as in the proof of \cref{lem:VOAGKM}, we can represent this intersection as a sum of degrees of conformal block divisors on $\M{0}{4}$ and $\M{1}{1}$. Since $\text{rank }\mathbb{V}_{g,n}(V, W^\bullet) \leq 1$, by \cref{thm:virdeg}, these degrees are $0$. Therefore, $\mathbb{D}_{g,n}(V, W^\bullet) = 0$.

    First, consider the case of $g \ge 2$. By \cref{lem:virdim} (1), it suffices to prove that $-\mathbb{D}_{g,n}(V, W^\bullet)$ is F-ample for $g \ge 2$. Specifically, we need to show that the intersection of F-curves with $\mathbb{D}_{g,n}(V, W^\bullet)$ is negative. Following the method used in the proof of \cref{lem:VOAGKM}, we will demonstrate this only for Type 1, Type 5, and Type 6 F-curves, as the proofs for other F-curves are similar. Note that each degree appearing in expressions like those in the proof of \cref{lem:VOAGKM} is nonpositive by \cref{thm:virdeg}. Therefore, it suffices to show that a negative component of the intersection exists.
    
    Consider the Type 1 F-curve. By \cref{thm:factor}, the intersection of $\mathbb{D}_{g,n}(V, W^\bullet)$ with the Type 1 F-curve is given by $\sum_W \text{rank }\mathbb{V}_{g-1,n+1}(V, W^\bullet\otimes W) \cdot \mathbb{D}_{1,1}(V, W)$. By \cref{lem:virdim} (1) and (2), $\mathbb{D}_{1,1}(V, V)$ is a summand of this expression. Since $\text{rank }\mathbb{V}_{1,1}(V, V) = k \ge 2$, by \cref{thm:virdeg}, $\deg \mathbb{D}_{1,1}(V, V) < 0$. Therefore, the intersection is negative.

    Next, consider a Type 5 F-curve. By \cref{thm:factor}, the intersection of $\mathbb{D}_{g,n}(V, W^\bullet)$ with the Type 5 F-curve is given by
    \begin{align*}
        \sum_{M_1, M_2, M_3, M_4} &\rank \mathbb{V}_{g_1,|I|+1}(V, W^{I}\otimes M_1)\cdot \rank \mathbb{V}_{g_2,|J|+1}(V, W^{J}\otimes M_2)\cdot \\&\rank \mathbb{V}_{g_3,|(I\cup J)^c|+1}(V, W^{(I\cup J)^c}\otimes M_3\otimes M_4) \cdot \mathbb{D}_{0,4}(V, \otimes_{i=1}^4 M_i). 
    \end{align*}
    By \cref{lem:virdim}, there is a summand where $M_1$ and $M_2$ are nontrivial simple modules and $M_3 = M_4 = W_k$. Hence, by \cref{lem:virdim} (3) and \cref{thm:virdeg}, this summand is negative. Therefore, the intersection is negative. The intersections for Type 2, Type 3, and Type 4 F-curves are similar.

    Finally, consider the case of Type 6 F-curves. By \cref{thm:factor}, the intersection is given by
    \[ \sum_{M_1, M_2, M_3, M_4} \left( \prod_{i=1}^{4}\text{rank }\mathbb{V}_{g_i,|I_i|+1}(V, W^{I_i}\otimes M_i)\right) \cdot \mathbb{D}_{0,4}(V, \otimes_{i=1}^4 M_i). \]
    At least one of the $g_i$'s is positive. Without loss of generality, assume $g_1 \ge 1$. By \cref{lem:virdim}, regardless of the choices of $M_2, M_3, M_4$, we can select $M_1$ to be any simple module. By the fusion rules, a summand exists where $M_2, M_3, M_4$ are nontrivial simple modules. If $k \ge 3$, then by \cref{lem:virfus}, we can choose $M_1$ such that $\text{rank } \mathbb{V}_{0,4}(V, \otimes_{i=1}^4 M_i) \ge 2$, which implies $\deg \mathbb{D}_{0,4}(V, \otimes_{i=1}^4 M_i) < 0$. If $k = 2$, there exists a summand where all $M_i$'s are $W_2$, and hence $\deg \mathbb{D}_{0,4}(V, \otimes_{i=1}^4 M_i) < 0$. Therefore, the intersection is negative.

    Now consider the case of $g=1$. The case of $n=1$ follows easily from \cref{thm:virdeg}. Assume $n \ge 2$. We need to prove that the intersection of $\mathbb{D}_{1,n}(V, W^\bullet)$ with F-curves of Type 1, Type 5, and Type 6 is negative. The proofs for the cases of Type 5 and Type 6 F-curves are very similar to those for $g \ge 2$, so we omit them here. For Type 1 F-curves, the intersection number is given by $\sum_W \text{rank }\mathbb{V}_{0,n+1}(V, W^\bullet\otimes W) \cdot \mathbb{D}_{1,1}(V, W)$. Since $n \ge 2$, by the fusion rules, there is a contribution from some module $W$ not equal to $W_2$. Hence, the intersection is also negative by \cref{thm:virdeg} and \cref{lem:virdim}.
     
    Finally, consider the case of $g=0$. If $k=2$, only one nontrivial $V$-module exists, $W=W_2$. By \cref{lem:linem1n} (4) or a direct argument using \cref{thm:factor}, $\text{rank } \mathbb{V}_{0,n}(V, W^\bullet) = F_{n-1}$, where $F$ denotes the Fibonacci sequence. Hence, it suffices to prove that $-\mathbb{D}_{0,n}(V, W^\bullet)$ is F-ample. By \cref{thm:factor}, the intersection of $\mathbb{D}_{0,n}(V, W^\bullet)$ with any F-curve is a positive multiple of $\deg \mathbb{D}_{0,4}(V, W^\bullet) = -1$, and is therefore negative.    

    For $k \ge 3$, we proceed by induction on $n$. The case of $n=4$ follows directly from \cref{thm:virdeg}. We will prove that, if $n \ge 5$ and $\text{rank } \mathbb{V}_{0,n}(V, W^\bullet) \ge 2$, then for every subset $I \subseteq [n]$ with $|I| = n-2$, there exists a nontrivial simple $V$-module $W$ such that $\text{rank } \mathbb{V}_{0,n-1}(V, W^I \otimes W) \ge 2$ and $\text{rank } \mathbb{V}_{0,3}(V, W^{I^c} \otimes W) \ge 1$. For any $F$-curve $[C]$ on $\M{0}{n}$, there exists a subset $I \subseteq [n]$ and an $F$-curve $D$ on $\M{0}{|I|+1}$ such that $|I| = n-2$ and $[C] = f_\ast [D]$, where $f:\M{0}{I\cup \left\{p\right\}}\to \M{0}{n}$ is a map attaching a stable rational curve with $|I^c|+1$ points to $p$. By the induction hypothesis, $\text{rank } \mathbb{V}_{0,n-1}(V, W^I \otimes W) \ge 2$ implies that $-\mathbb{D}_{0,n}(V, W^\bullet)$ is F-ample.

    First, consider the case of $n=5$. Assume $I=\left\{1,2,3\right\}$. By the fusion rules, $\boxtimes_{i=1}^3 W^i=\sum_{j\stackrel{\text{2}}{=}a}^b n_jW_j$ where (1) $b-a>4$ or (2) $b=a+2$ and $n_a$ or $n_b$ is $>1$. In any case, since $\text{rank } \mathbb{V}_{0,n}(V, W^\bullet) \ge 2$, there exists $c$ such that $n_c\ge 2$ and $\text{rank } \mathbb{V}_{0,3}(V, W^{I^c} \otimes W_c) \ge 1$, by \cref{lem:virfus}. Hence, this $W_c$ satisfies the requirement.
 
    Now consider the case of $n\ge 6$. Since $|I| \ge 4$ and $k \ge 3$, by \cref{lem:virfus}, $\boxtimes_{i\in I} W^i=\sum_{j\stackrel{\text{2}}{=}a}^b n_jW_j$, where $n_j \ge 2$ for $j \ne a, b$, and $b-a \ge 4$. Without loss of generality, assume $I = [n-2]$. Then, by the fusion rules, we may write $W^{n-1}\boxtimes W^n=\sum_{i\stackrel{\text{2}}{=}c}^d W_j$ where $c$ and $d$ are integers of the same parity as $a$ and $b$. By comparing these two expressions and using the condition $\text{rank } \mathbb{V}_{0,n}(V, W^\bullet) \ge 2$, we deduce that there exists an integer $j$ such that $c \le j \le d$, $j$ has the same parity as $c$ and $d$, and $n_j \ge 2$ (here we use the condition $b-a \ge 4$). Hence, $W = W_j$ satisfies the required condition.

    To complete the proof, it suffices to show that if $n > 2k-1$, then $\text{rank } \mathbb{V}_{0,n}(V, W^\bullet) \ge 2$. If $k=2$, this follows trivially from the earlier observation that $\text{rank } \mathbb{V}_{0,n}(V, W^\bullet) = F_{n-1}$, where $F$ is the Fibonacci sequence. Now assume $k \ge 3$. By \cref{lem:virfus}, since $n-2 \ge 2k-2$, $W^{[n-2]}$ contains every simple $V$-module as a summand. Furthermore, since any product of two nontrivial simple modules has at least two simple modules as components, it follows that $\text{rank } \mathbb{V}_{0,n}(V, W^\bullet) \ge 2$. This completes the proof.
\end{proof}

\begin{rmk}
    The bound in \cref{thm:virfnef}(3) is the best possible, since $\mathbb{D}_{0,2k-1}(\vir{2k+1}, W_2^{\otimes(2k-1)})=0$.
\end{rmk}

\begin{rmk}
    As a consequence, we can trivially verify \cite[Question 2]{DG23} for $V = \vir{2k+1}$. Since the Zhu algebra of discrete series Virasoro algebras is commutative \cite{FZ92}, the degree zero part of any simple module over $\vir{2k+1}$ is $1$-dimensional. Hence, the only case where \cite[Question 2]{DG23} applies is when $\text{rank } \mathbb{V}_{0,n}(V, W^\bullet) \le 1$. \cref{thm:virfnef} implies that $\mathbb{V}_{0,n}(V, W^\bullet)$ is trivial in this case.
\end{rmk}

\subsection{Nefness and Ampleness of Conformal Block Divisors}\label{subsec:nef}
This section is devoted to proving the following theorem and related discussions.

\begin{thm}\label{thm:virnef}
    Let $k$ be an integer such that $1 < k < 9$. Then, any conformal block divisor $-\mathbb{D}_{g,n}(V, W^\bullet)$ for $V = \vir{2k+1}$ is nef. If every simple module $W^i$ is nontrivial, it is either trivial or ample. Moreover, under the same conditions, if $g = 0$ and $n > 2k+1$, or $g = 1$ and $n > 1$, or $g \geq 2$, then $-\mathbb{D}_{g,n}(V, W^\bullet)$ is ample.
\end{thm}

We begin with two criteria for nefness and ampleness. Only \cref{lem:effnef} will be used in the proof of \cref{thm:virnef}. However, \cref{lem:Knef} provides a simpler proof of \cref{thm:virnef} if $k\le 4$. We will discuss this in \cref{subsec:rmk} and another application of \cref{lem:Knef}.

\begin{lem}\label{lem:effnef}
    Let $X$ be a smooth projective variety and $D$ be an effective divisor on $X$ of the form $D=\sum_i a_i E_i$ where $E_i$ are prime divisors and $a_i>0$.

    \begin{enumerate}
        \item If $D|_{E_i}$ is nef for every $i$, then $D$ is nef.
        \item If $D|_{E_i}$ is ample for every $i$ and $E_i$'s span $\text{Pic}_\Q (X)$ as a $\Q$-vector space, then $D$ is also ample.
    \end{enumerate}
\end{lem}

\begin{proof}
    (1) Suppose not. Let $C$ be an integral curve such that $D \cdot C < 0$. Since $D = \sum_i a_i E_i$, $C$ must be contained in the support of $D$, so $C \subseteq E_i$ for some $i$. However, this contradicts the assumption.
    
    (2) By (1), $D$ is nef. Let $L$ be an ample divisor on $X$. Since $a_i > 0$ for each $i$, there exists $\epsilon > 0$ and $b_i > 0$ such that $D = \epsilon L + \sum_i b_i E_i$. In particular, $D$ is big. Let $Z$ be an irreducible subvariety of $D$ such that $Z \not\subseteq E_i$ for every $i$. Then, $E_i|_Z$ is an effective divisor on $Z$. Hence, $D|_Z$ is a sum of ample ($L|_Z$) and effective ($E_i|_Z$) divisors, so $D|_Z$ is big. Since $D|_Z$ is nef, the top dimensional intersection $D^{\dim Z}\cdot Z$ is positive. If $Z \subseteq E_i$ for some $i$, since $D|_{E_i}$ is ample for every $i$, this implies $D^{\dim Z}\cdot Z>0$. Hence, by the Nakai-Moishezon criterion, $D$ is ample.
\end{proof}

\begin{lem}\label{lem:Knef}[cf. \cite[Theorem 1.2 (2)]{KM13}, \cite[Theorem 4]{GF03}. char $k=0$]
    Let $D$ be a $\mathbb{Q}$-divisor on $\M{0}{n}$ of the form $D=K_{\M{0}{n}}+E$, where $E=\sum_I a_I\delta_{0,I}$ an effective sum of boundary divisors (i.e. $a_I\ge 0$). If $D$ is F-nef (resp. F-ample) and the restrictions of $D$ to boundary divisors are nef (resp. ample), then $D$ is nef (resp. ample).  
\end{lem}

\begin{rmk}\label{rmk:compare}
Note that \cite[Theorem 1.2 (2)]{KM13} and \cite[Theorem 4]{GF03} state a slightly different form of \cref{lem:Knef}. Specifically, they assume that $a_I \leq 1$ instead of the nefness of boundary restrictions of $D$. The following proof of \cref{lem:Knef} adopts the approach of \cite[Theorem 3.1]{Gib09}.
\end{rmk}

\begin{proof}
    We will prove the nef case of the lemma. The ample case is almost identical since $D$ is ample if and only if $D \cdot R > 0$ for any extremal ray $R$ of the closed cone of curves. Assume $D$ is not a nef divisor. Then, there exists a class $C \in \NE{\M{0}{n}}$ which generates an extremal ray such that $D \cdot C < 0$. Since $C$ is contained in the closed cone of curves, $C = \lim_{n \to \infty} C_n$, where $C_n$'s are effective sums of curves. Choose integral curves $C^{1}_{n,i}$ and $C^{2}_{n,j}$ so that
\begin{enumerate}
        \item $[C_n]=\sum_i x_n^i [C^{1}_{n,i}]+\sum_j y_n^j [C^{2}_{n,j}]$ where $x_n^i, y_n^j>0$.
        \item $E\cdot C^{1}_{n,i}\ge 0$ and $E\cdot C^{2}_{n,j}<0$.
    \end{enumerate}
    Define $C_n^1:=\sum_i x_n^i [C^{1}_{n,i}]$ and $C_n^2:=\sum_j y_n^j [C^{2}_{n,j}]$. Let $L_1,\cdots, L_r\in \text{Pic}(\M{0}{n})$ be ample divisors on $\M{0}{n}$ which form a basis of the $\Q$-Picard group. Then
    \[ L:\text{A}_1(\M{0}{n})\otimes \R\to \R^r \]
    defined by $X\mapsto (L_i\cdot X)_{i=1}^{r}$ is an isomorphism of $\R$-vector spaces. In particular, this is a homeomorphism.
    Let $b_i=L_i\cdot C$. Then, by definition, for sufficiently large $n$, $L\left(C_n^1 \right), L\left(C_n^2 \right)\in\prod_i [0, b_i+1]$. Hence, by the Bolzano-Weierstrass theorem, after taking a subsequence, we may assume that $C_n^1, C_n^2$ converges to $C^1, C^2\in \text{A}_1(\M{0}{n})\otimes \R$. Since they are limits of $C_n^1, C_n^2$, they are contained in $\NE{\M{0}{n}}$. Moreover, since $C$ is an extremal ray of $\NE{\M{0}{n}}$ and $C=C^1+C^2$, $C^1$ and $C^2$ are both a multiple of $C$. In particular, $C$ is a nonzero multiple of $C^1$ or $C^2$.

    Assume $C$ is a nonzero multiple of $C^1$. By the definition of $C_n^1$, $E \cdot C^1 \geq 0$. Hence, $K_{\M{0}{n}} \cdot C^1 \leq 0$, which implies $K_{\M{0}{n}} \cdot C \leq 0$. By \cite[Theorem 4]{GF03}, this implies $C$ is equivalent to an F-curve. However, this contradicts the assumption that $D$ is F-nef.

    Now, assume $C$ is a nonzero multiple of $C^2$. Since $E \cdot C_{n,j}^2 < 0$ and $C_{n,j}^2$ is an integral curve, $C_{n,j}^2$ must be contained in the support of $E$. Therefore, each $C_n^2$ is contained in the image of $\NE{\Delta}$ under the map $i_\ast : \text{A}_1(\Delta) \otimes \R \to \text{A}_1(\M{0}{n}) \otimes \R$, whre $\Delta$ is the boundary. Note that, by \cite[Theorem 9.1]{Roc70}, $i\ast\left(\NE{\Delta}\right)$ is a closed subset of $\text{A}_1(\M{0}{n}) \otimes \R$ since $i: \Delta \to \M{0}{n}$ is a closed embedding. Hence, $C^2 = i_\ast C_0$ for some $C_0 \in \NE{\Delta}$. By the assumption, $D \cdot C^2 < 0$, which implies $i^\ast D \cdot C_0 < 0$. This contradicts the assumption that the boundary restrictions of $D$ are nef. All in all, $D$ is a nef divisor.
\end{proof} 

Note that \cref{lem:effnef} and \cref{lem:Knef} are variants of standard tools for proving nefness on $\M{0}{n}$. See \cite{KM13, GF03, Gib09} for earlier discussions and \cite{Pix13} for the limitations of this method.

We will prove \cref{thm:virnef} using \cref{lem:effnef}. Note that, by \cref{lem:VOAGKM}, it suffices to consider the case $g=0$. To do this, we need to establish the following: (1) the nefness (or ampleness, respectively) of the line bundle when restricted to divisors, and (2) the representation of the line bundle as a sum of sufficiently many effective divisors. Our proof proceeds by induction on the number of marked points. In the inductive step, condition (1) is ensured by the induction hypothesis and factorization (\cref{thm:factor}).

(2) is the hardest and most technical part of the proof. For this, we will use \cref{thm:c1cdf} and the following `big average' identity from \cite[Lemma 2.8]{Gib09}:
    \[ \psi_i=\sum_{Y\subseteq [n]\setminus\left\{i\right\}}\frac{(n-1-|Y|)(n-2-|Y|)}{(n-1)(n-2)}\delta_{0, Y\cup\left\{i\right\} }. \]
\cref{thm:c1cdf} provides the following expression:
\[ \mathbb{D}_{0,n}\left(V, W^\bullet\right)=\rank \mathbb{V}_{0,n}(V, W^\bullet)\cdot \sum_{i=1}^n h_i\psi_i-\sum_I b_{0,I}\delta_{0,I} \]
Here, we can substitute $\psi_i$ by its big average and let $\sum_{I} c_I \delta_{0,I} $ be the resulting sum. We will call this the \textbf{standard form} of the conformal block divisor for simplicity. Note that this is the original form of the first Chern class formula \cref{thm:c1cdf}, proved by Fakhruddin \cite[Corollary 3.4]{Fak12}.

\begin{thm}\label{thm:genvireff}
    Let $V=\text{Vir}_{2,2k+1}$, $n\ge 4k-4$ and $W_{a_i}$ be a set of nontrivial simple $V$-modules. Then we have an expression
    \[ \mathbb{D}_{0,n}\left(V, \otimes_{i=1}^n W_{a_i} \right)=\sum_{I} c_I \delta_{0,I} \]
    where $c_I<0$. More precisely, for any choice of $n$ and $k$, if we let $\sum_{I} c_I \delta_{0,I}$ be the standard form of the conformal block divisor, then we have $c_I<0$ if$|I|\le |I^c|$ and
    \[ |I|\ge k\text{ or }\frac{1}{2k-1-|I|}\ge \frac{1}{n-|I|}+\frac{1}{n-2} \]
\end{thm}

This theorem requires a detailed analysis of conformal weights of $\vir{2k+1}$. Define
\[ w_i:=-2(2k+1)h_i=(i-1)(2k-i) \]
be the normalized conformal weights.

\begin{lem}\label{lem:fusbd}
    Let $1\le a_1,\cdots, a_i\le k$ be integers and 
    \[ \boxtimes_{j=1}^{i}W_{a_j}=\sum_{r=1}^k n_r W_r. \]
    Then $n_r=0$ for $r\ge \sum_j a_j-(i-1)$.
\end{lem}

\begin{proof}
    Use induction. The case $i = 2$ is straightforward from the fusion rules. The induction step is also a direct consequence of the fusion rules.
\end{proof}

\begin{lem}\label{lem:cwlem}
    Let $2\le a_1,\cdots, a_i\le k$ be integers, $i<k$ and $m:=\min\left(k, \sum_j a_j-(i-1) \right)$. Then
    \[ \frac{\sum_{j=1}^{i}w_{a_j}}{w_m}\ge \frac{2k-2}{2k-1-i}. \]
\end{lem}

\begin{proof}
    As a function of $i$, $w_i$ is an increasing function in the interval $[1,k]$. Since $i < k$, we may assume that $k \ge \sum_j a_j - (i-1)$; hence $m = \sum_j a_j - (i-1)$. Note that, since $w_i$ is a concave function in $i$, for any $2 \le a < b \le k-1$,
    \[ w_a+w_b\ge w_{a-1}+w_{b+1}. \]
    Therefore, by applying this iteratively, we may assume that $a_1 = \cdots = a_{i-1} = 2$ and $a_i = r+1$. Then $m = r + i$, and it is enough to show that
    \[ \frac{2(k-1)(i-1)+r(2k-r-1)}{(r+i-1)(2k-r-i)}\ge \frac{2k-2}{2k-1-i}.  \]
    By simplifying the inequality, it is enough to show that
    \[ (i-1)r(r+2k-3)\ge 2(k-1)(i-1) \]
    which is evident since $r\ge 1$.
\end{proof}

\begin{proof}[proof of \cref{thm:genvireff}]
    Let $W^i=W_{a_i}$, $2\le a_i\le k$. Then
    \begin{align*}
        -\frac{2(2k+1)}{\rank \mathbb{V}_{0,n}(V, W^\bullet)}\mathbb{D}_{0,n}\left(V, W^\bullet\right)=\sum_{i=1}^n w_i\psi_i-\sum_I \left(-\frac{2(2k+1)b_{0,I}}{\rank \mathbb{V}_{0,n}(V, W^\bullet)}\right)\delta_{0,I}=\sum_I d_I \delta_{0,I}
    \end{align*}
    where
    \[ d_I=\sum_{p\in I} w_{a_p}\frac{(n-|I|)(n-|I|-1)}{(n-1)(n-2)}+\sum_{p\not\in I} w_{a_p}\frac{|I|(|I|-1)}{(n-1)(n-2)}+\frac{2(2k+1)b_{0,I}}{\rank \mathbb{V}_{0,n}(V, W^\bullet)}. \]
    It is enough to show that $d_I>0$ under the assumption. Note that
    \[ \sum_{1\le a\le k} \rank \mathbb{V}_{0,|I|+1}(V, W^I\otimes W_a) \cdot \rank \mathbb{V}_{0,|I^c|+1}(V, W^{I^c}\otimes W_a)=\rank \mathbb{V}_{0,n}(V, W^\bullet) \]
    by \cref{thm:factor}. Hence, 
    \[ -\frac{2(2k+1)b_{0,I}}{\rank \mathbb{V}_{0,n}(V, W^\bullet)}=\sum_{1\le a\le k} w_a \cdot \frac{\rank \mathbb{V}_{0,|I|+1}(V, W^I\otimes W_a) \cdot \rank \mathbb{V}_{0,|I^c|+1}(V, W^{I^c}\otimes W_a)}{\rank \mathbb{V}_{0,n}(V, W^\bullet)}\le w_m \]
    where $m = \min\left(k, \sum_j a_j - (i-1)\right)$, since $m$ is the largest value $a$ can take in this expression by \cref{lem:fusbd}. Hence, it is enough to show that
    \[  \sum_{p\in I} w_{a_p}\frac{(n-|I|)(n-|I|-1)}{(n-1)(n-2)}+\sum_{p\not\in I} w_{a_p}\frac{|I|(|I|-1)}{(n-1)(n-2)}> w_m. \]
    Note that we choose $I$ so that $2|I|\le n$.
 
    \textbf{Case 1. }$|I|\ge k$. Since $a_i\ge 2$, $w_{a_i}\ge 2(k-1)$ for every $i$. Hence
    \begin{align*}
        &\sum_{p\in I} w_{a_p}\frac{(n-|I|)(n-|I|-1)}{(n-1)(n-2)}+\sum_{p\not\in I} w_{a_p}\frac{|I|(|I|-1)}{(n-1)(n-2)}\\
        &\ge 2(k-1)\frac{|I|(n-|I|)(n-|I|-1)+|I|(|I|-1)(n-|I|)}{(n-1)(n-2)}\\& =2(k-1)\frac{|I|(n-|I|)}{n-1}> k(k-1)\ge w_m.
    \end{align*}
    In the last part, we used $\frac{n - |I|}{n - 1} >\frac{1}{2}$, which follows from $2|I| \le n$. 

    \textbf{Case 2. }$|I|<k$ and $\sum_{p\in I}w_{a_p}-w_m\ge |I|(2k-2)$. Then 
    \begin{align*}
        &\sum_{p\in I} w_{a_p}\frac{(n-|I|)(n-|I|-1)}{(n-1)(n-2)}+\sum_{p\not\in I} w_{a_p}\frac{|I|(|I|-1)}{(n-1)(n-2)}\\
        &\ge w_m \frac{(n-|I|)(n-|I|-1)}{(n-1)(n-2)}+\sum_{p\in I} 2(k-1)\frac{(n-|I|)(n-|I|-1)}{(n-1)(n-2)}+\sum_{p\not\in I} 2(k-1)\frac{|I|(|I|-1)}{(n-1)(n-2)}\\
        &\ge w_m \frac{(n-|I|)(n-|I|-1)}{(n-1)(n-2)}+2(k-1)\frac{|I|(n-|I|)}{n-1}
    \end{align*}
    Hence, it is enough to show
    \[ 2(k-1)\frac{|I|(n-|I|)}{n-1}> w_m\frac{(|I|-1)(2n-|I|-2)}{(n-1)(n-2)}. \]
    By using $|I|> |I|-1$ and $k(k-1)\ge w_m$, this simplifies to
    \[ 2(n-|I|)\ge k\cdot \frac{2n-|I|-2}{n-2}\text{, or equivalently, }\frac{2}{k}\ge \frac{1}{n-|I|}+\frac{1}{n-2}. \]
    Since $\frac{1}{2k-1-|I|} \ge \frac{1}{n-|I|} + \frac{1}{n-2}$ and $|I| < k$, this is straightforward to see. Note that the preceding inequality can also be proved under the $n \ge 2k$ assumption.

    \textbf{Case 3. }$|I|<k$ and $\sum_{p\in I}w_{a_p}-w_m<|I|(2k-2)$. Let $\Delta:=\sum_{p\in I}w_{a_p}-w_m$. Then 
    \begin{align*}
        &\sum_{p\in I} w_{a_p}\frac{(n-|I|)(n-|I|-1)}{(n-1)(n-2)}+\sum_{p\not\in I} w_{a_p}\frac{|I|(|I|-1)}{(n-1)(n-2)}\\
        &\ge w_m\frac{(n-|I|)(n-|I|-1)}{(n-1)(n-2)}+\Delta\frac{(n-|I|)(n-|I|-1)}{(n-1)(n-2)}+(2k-2)|I|\frac{(n-|I|)(|I|-1)}{(n-1)(n-2)}\\
        &> w_m\frac{(n-|I|)(n-|I|-1)}{(n-1)(n-2)}+\Delta \frac{n-|I|}{n-1}.
    \end{align*}
     Hence, it is enough to show
    \[ \Delta\frac{n-|I|}{n-1}\ge w_m\frac{(|I|-1)(2n-|I|-2)}{(n-1)(n-2)}. \]
    By \cref{lem:cwlem}, $\Delta\ge \frac{|I|-1}{2k-|I|-1}w_m $. Therefore, it is enough to show
    \[ \frac{n-|I|}{2k-|I|-1}\ge \frac{(2n-|I|-2)}{(n-2)} \]
    which is another form of $\frac{1}{2k-1-|I|}\ge \frac{1}{n-|I|}+\frac{1}{n-2}$.

    Finally, to prove the first assertion of \cref{thm:genvireff}, we need to show that $n \ge 4k - 4$ implies $\frac{1}{2k-1-|I|} \ge \frac{1}{n-|I|} + \frac{1}{n-2}$. We will omit the proof, as this part is rather straightforward.
\end{proof}

We will prove nefness of $\mathbb{D}_{0,n}\left(\text{Vir}_{2,2k+1}, \otimes_{i=1}^n W_{a_i} \right)$ by verifying the following conjecture for specific values of $k$. 

\begin{conj}\label{conj:genvireff}
    Let $V=\text{Vir}_{2,2k+1}$ and $\left\{W_{a_i}\right\}_{i=1}^n$ be a set of nontrivial simple $V$-modules, i.e., $a_i \neq 1$. Then $\mathbb{D}_{0,n}\left(V, \otimes_{i=1}^n W_{a_i} \right) = 0$ or $-\mathbb{D}_{0,n}\left(V, \otimes_{i=1}^n W_{a_i} \right)$ lies in the interior of the cone generated by the boundary divisors.
\end{conj}

\begin{rmk}\label{rmk:eff}
    By \cref{thm:genvireff}, to prove \cref{conj:genvireff}, it is enough to consider finitely many cases where $n < 4k - 4$. Moreover, for each case, verifying whether $-\mathbb{D}_{0,n}\left(\vir{2k+1}, \otimes_{i=1}^n W_{a_i} \right)$ is contained in the cone generated by the boundary divisors is computable, since it is a linear programming problem. Hence, for each $k$, \cref{conj:genvireff} is verifiable within a finite amount of time.
\end{rmk}

\begin{prop}\label{prop:sublem}
    If \cref{conj:genvireff} holds for a fixed value of $k$ and for every choice of simple modules, then $-\mathbb{D}_{0,n}\left(\vir{2k+1}, \otimes_{i=1}^n W_{a_i} \right)$ is nef. Furthermore, if all simple modules are nontrivial, it is either zero or ample.
\end{prop}

\begin{proof}
    First, by \cref{thm:propvac}, \cref{lem:effnef} (1), and the assumption, the first assertion follows directly. We will prove the second assertion using induction on $n$. For $n \leq 7$, this follows from \cref{cor:fconj} and \cref{thm:virfnef}. Now assume that $D = -\mathbb{D}_{0,n}\left(\vir{2k+1}, \otimes_{i=1}^n W_{a_i} \right)$ is a nontrivial Virasoro conformal block divisor on $\M{0}{n}$ for nontrivial simple modules. The induction step relies on \cref{lem:effnef} (2). The effectivity condition follows from the assumption about \cref{conj:genvireff}, so it suffices to show that the restriction of $D$ to the boundary divisors is ample.

    Consider $D|_{\Delta_{0,I}}$. Without loss of generality, assume $|I| \geq \lfloor \frac{n}{2} \rfloor$. If $|I| = n-2$, then $\Delta_{0,I}$ corresponds to the natural map $\iota: \M{0}{n-1} \to \M{0}{n}$ obtained by attaching a projective line with three marked points. Since $n-2 \geq 6$, by \cref{lem:virfus}, there exists a nontrivial simple module $W_j$ such that
    \[\rank \mathbb{V}_{0,n-1}\left(\vir{2k+1}, \otimes_{i \in I} W_{a_i} \otimes W_j \right) \geq 2 \quad \text{and} \quad \rank \mathbb{V}_{0,3}\left(\vir{2k+1}, \otimes_{i \in I^c} W_{a_i} \otimes W_j \right) \geq 1.\]
    Then, by \cref{thm:factor}, $\iota^\ast D$ is the sum of a nef line bundle and
    \[-\mathbb{D}_{0,n-1}\left(\vir{2k+1}, \otimes_{i \in I} W_{a_i} \otimes W_j \right).\]
    By \cref{thm:virfnef} and the induction hypothesis, this is ample. Therefore, $D|_{\Delta_{0,I}}$ is ample.

    Now, assume $|I| \neq n-2$. Then $\Delta_{0,I}$ corresponds to the natural map $\iota: \M{0}{|I|+1} \times \M{0}{|I^c|+1} \to \M{0}{n}$. Since $D$ is nontrivial, 
    \[\rank \mathbb{V}_{0,n}\left(\vir{2k+1}, \otimes_{i=1}^n W_{a_i} \right) \geq 2.\]
    Hence, by \cref{lem:virfus}, there exist $W_l$ and $W_j$ such that
    \begin{multline*}
        \rank \mathbb{V}_{0,|I|+1}\left(\vir{2k+1}, \otimes_{i \in I} W_{a_i} \otimes W_l \right) \geq 2, \quad \rank \mathbb{V}_{0,|I^c|+1}\left(\vir{2k+1}, \otimes_{i \in I^c} W_{a_i} \otimes W_l \right) \geq 1,\\
        \rank \mathbb{V}_{0,|I^c|+1}\left(\vir{2k+1}, \otimes_{i \in I^c} W_{a_i} \otimes W_j \right) \geq 2, \quad \text{and} \quad \rank \mathbb{V}_{0,|I|+1}\left(\vir{2k+1}, \otimes_{i \in I} W_{a_i} \otimes W_j \right) \geq 1.
    \end{multline*}
    Then, by \cref{thm:factor}, $\iota^\ast D$ is the sum of a nef line bundle and
    \[-\pi_1^\ast \mathbb{D}_{0,|I|+1}\left(\vir{2k+1}, \otimes_{i \in I} W_{a_i} \otimes W_l \right) - \pi_2^\ast \mathbb{D}_{0,|I^c|+1}\left(\vir{2k+1}, \otimes_{i \in I^c} W_{a_i} \otimes W_j \right),\]
    which is ample by the induction hypothesis. This completes the proof.

\end{proof}

\begin{proof}[proof of \cref{thm:virnef}]
    First, we verify \cref{conj:genvireff} for $k < 9$. This can be checked through brute force using a computer. A Python implementation for this verification is available in \cite[effectivity]{Choigit24}. Below, we provide a brief explanation of the code. In most cases, the standard form suffices to produce the desired effective sum as described in \cref{thm:genvireff}. In such cases, the code routinely checks that the coefficient of the standard form is negative. However, in certain exceptional cases, the standard form does not work (e.g., $-\mathbb{D}_{0,8}\left(\vir{17}, W_2^{\otimes 5} \otimes W_4^{\otimes 2} \otimes W_7 \right)$). For these cases, explicitly solving the corresponding linear programming problem is necessary. We use the Z3 solver to handle linear programming to obtain definitive results. We need to verify this for the cases not covered by \cref{thm:genvireff}. There are only finitely many such cases, as \cref{thm:genvireff} verifies the conjecture for $n \geq 4k - 4$. Moreover, we use \cref{thm:steff} to reduce the runtime.
    
    The code successfully verifies \cref{conj:genvireff} for $k < 9$. On the laptop of the author, the runtime of the code is as follows: less than 2 seconds for $k \leq 5$, 79 seconds for $k = 6$, 4158 seconds for $k = 7$, and 48477 seconds for $k = 8$. 

    By \cref{prop:sublem}, this completes the proof for the case of $g = 0$. Moreover, by \cref{lem:VOAGKM} (2), this implies that any Virasoro conformal block divisor is nef for $k < 9$. Hence, it only remains to prove the ampleness assertion. We will use \cref{thm:fconj0}. Note that, since the ample cone is the interior of the nef cone, the ample version of \cref{thm:fconj0} also holds. In other words, a line bundle $L$ on $\M{g}{n}$ is ample if and only if its intersection with F-curves is positive and $F^\ast L$ is ample.

    Let $L = -\mathbb{D}_{g,n}\left(\vir{2k+1}, \otimes_{i=1}^n W_{a_i} \right)$. By \cref{thm:virfnef}, it is enough to prove that $F^\ast L$ is ample. First, consider the case of $g \geq 2$. By \cref{thm:factor} and \cref{thm:virdeg} (2), $F^\ast L$ is the sum of a nef line bundle and
    \[ -c\cdot\mathbb{D}_{0,n+g}\left(\vir{2k+1}, \otimes_{i=1}^n W_{a_i} \otimes W_{k}^{\otimes g} \right) \]
    for some $c > 0$. Since $g \geq 2$, by the fusion rules, 
    \[ \rank \mathbb{V}_{0,n+g}\left(\vir{2k+1}, \otimes_{i=1}^n W_{a_i} \otimes W_{k}^{\otimes g} \right) \geq 2. \]
    Therefore, this is ample by \cref{thm:virfnef} and the genus $0$ case of \cref{thm:virnef}. Hence, $F^\ast L$ is ample.

    Now, consider the case of $g = 1$. If $n \leq 6$, then by \cref{cor:fconj}, the assertion holds. Now, assume $n > 6$. By \cref{lem:virfus}, there exists a nontrivial simple module $W_i$ such that
    \[\rank \mathbb{V}_{0,n+1}\left(\vir{2k+1}, \otimes_{i=1}^n W_{a_i} \otimes W_{i} \right) \geq 2. \]
    Hence, by \cref{thm:virnef} and \cref{thm:factor}, $F^\ast L$ is ample. This completes the proof.

\end{proof}

\subsection{Digression on \texorpdfstring{\cref{subsec:nef}}{}}\label{subsec:rmk}
Here, we provide some results that allow for a more efficient proof of \cref{thm:main1} for $2 \le k \le 4$ and establish a connection to Fakhruddin's work \cite{Fak12}. \cref{lem:Knef} and the induction process allow us to prove a general nefness theorem for conformal block divisors.

\begin{thm}\label{thm:VOAnefcw}[char $k=0$]
    Let $V$ be a $C_2$ and rational VOA such that 
    \begin{enumerate}
        \item $V$ is a simple module over itself
        \item The conformal block divisors of $V$ are nef on $\M{0}{4}$.
        \item The conformal weights of nontrivial simple modules over $V$ are negative. If $h_{\text{min}}$ and $h_{\text{max}}$ denote the conformal weights with minimum and maximum absolute value, then $h_{\text{max}} \leq 2h_{\text{min}}$.
    \end{enumerate}
    Then $-\mathbb{D}_{0,n}(V, W^\bullet)$'s are nef. Moreover, if $-\mathbb{D}_{0,n}(V, W^\bullet)$ is F-ample, then it is ample.
\end{thm}

\begin{rmk}\label{rmk:char}
    \cref{lem:Knef} and \cref{thm:VOAnefcw} are the only results outside \cref{sec:pre} that require the base field to have characteristic $0$. This restriction arises from \cite{KM13}, which relies on standard cone and contraction theorems. All other positivity statements, including \cref{cor:f4g2} (see \cref{thm:virnef}), can be established using only \cref{lem:effnef}, which is independent of the characteristic.
\end{rmk}

\begin{proof}
    By (2) and \cref{lem:VOAGKM} (1), $-\mathbb{D}_{0,n}(V, W^\bullet)$'s are F-nef. We will use induction on $n$ and \cref{lem:Knef}. Since the F-conjecture is true for $n\le 7$, the base case directly follows. By \cref{thm:factor}, from the induction hypothesis, we deduce that the last assumption of \cref{lem:Knef} holds. By the propagation of vacua (\cref{thm:propvac}), we may assume that every $W_i$'s are nontrivial. If $\rank \mathbb{V}_{0,n}(V, W^\bullet)=0$, then the theorem is trivial. Assume not.

    By \cref{thm:c1cdf}, using $\lambda=0$ if $g=0$,
    \begin{align*}
         &\mathbb{D}_{0,n}\left(V, W^\bullet\right)=\rank \mathbb{V}_{0,n}(V, W^\bullet)\left(\sum_{i=1}^n h_i\psi_i\right)-\sum b_{i,I}\delta_{i,I}\\&
         =\rank \mathbb{V}_{0,n}(V, W^\bullet)\cdot h_{\text{min}}\cdot\left(\psi-\sum_{I}\frac{b_{0,I}}{\rank \mathbb{V}_{0,n}(V, W^\bullet)\cdot h_{min}}\delta_{0,I}+\sum_{i=1}^{n}\frac{h_i-h_{\text{min}}}{h_{\text{min}}}\psi_i  \right)\\&
         =\rank \mathbb{V}_{0,n}(V, W^\bullet)\cdot h_{\text{min}}\cdot\left(K_{\M{0}{n}}+\sum_{I}\left(2-\frac{b_{0,I}}{\rank \mathbb{V}_{0,n}(V, W^\bullet)\cdot h_{min}}\right)\delta_{0,I}+\sum_{i=1}^{n}\frac{h_i-h_{\text{min}}}{h_{\text{min}}}\psi_i  \right)
    \end{align*}
    since $K_{\M{0}{n}}=\psi-2\delta$. By \cref{thm:factor},
    \[ b_{0,I}=\sum_{W\in S} h_W \cdot \rank \mathbb{V}_{0,|I|+1}(V, W^I\otimes W) \cdot \rank \mathbb{V}_{0,|I^c|+1}(V, W^{I^c}\otimes W')\ge h_{\text{max}}\cdot \rank \mathbb{V}_{0,n}(V, W^\bullet), \]
   and hence $2-\frac{b_{0,I}}{\rank \mathbb{V}_{0,n}(V, W^\bullet)\cdot h_{min}}\ge 0$. Since $\psi_i$'s are also a sum of boundary divisors, we can apply \cref{lem:Knef} here and deduce $\mathbb{D}_{0,n}\left(V, W^\bullet\right)$ is nef. 

   For the last assertion, we prove a slightly stronger statement: if an effective sum $D$ of $-\mathbb{D}_{0,n}(V, W^\bullet)$'s is F-ample, then it is ample. With this statement, the overall argument remains the same. The F-conjecture for $n \leq 7$ establishes the base case. By the preceding computation, $D$ is of the form $c(K_{\M{0}{n}} + E)$ for some $c < 0$ and an effective sum of boundaries $E$.  By \cref{thm:factor}, any boundary restriction of $D$ is also an effective sum of $\mathbb{D}_{0,n}(V, W^\bullet)$'s. Moreover, any boundary restriction of an F-ample divisor is also an F-ample divisor. From this and the induction hypothesis, the boundary restrictions of $D$ are ample. Therefore, by \cref{lem:Knef}, $D$ is ample.
\end{proof}

\begin{rmk}\label{rmk:VOA1}
    In particular, if $V$ is a VOA with exactly two simple modules, one of them is itself, then (3) is immediate. Also, one can easily state and prove the nef and ample version of this theorem.
\end{rmk}

\begin{cor}\label{cor:virnefeasy}
    If $V=\vir{2k+1}$ for $2\le k\le 4$, then the conformal block divisors are nef, and it is ample if it is F-ample.
\end{cor}

\begin{proof}
    By \cref{subsec:vir}, $V$ satisfies the assumptions of \cref{thm:VOAnefcw}. Hence $-\mathbb{D}_{0,n}(V, W^\bullet)$'s are nef. By \cref{lem:VOAGKM} (2), $-\mathbb{D}_{0,n}(V, W^\bullet)$'s are nef. The last assertion also follows from \cref{thm:VOAnefcw}.
\end{proof}

\begin{lem}\label{lem:fusion}
    Let $V_1, V_2$ be two $C_2$-cofinite and rational VOAs such that they
    \begin{enumerate}
        \item Are simple modules over themselves,
        \item Have exactly two simple modules, and
        \item Share the same fusion rules.
    \end{enumerate}
    Then their conformal block divisors on $\overline{\rm{M}}_{0,n}$ are proportional.
\end{lem}

\begin{proof}
    Let $W_i$ be the nontrivial simple module of $V_i$ with conformal weight $h_i$ for $i=1,2$. It is enough to prove that
    \[ h_2\mathbb{D}_{0,n}(V_1, W_1^{\otimes n})=h_1\mathbb{D}_{0,n}(V_1, W_1^{\otimes n}). \]
    Since the rank of the conformal blocks is determined by the fusion rules (cf. \cref{thm:factor}), the ranks for $V_1, V_2$ are the same. Hence, this follows from \cref{thm:c1cdf}. 
\end{proof}

Note that this argument is essentially the same as the one in \cite{Fak12}. There, Fakhruddin observed that if $\mathfrak{g}$ is $F_2$ or $G_4$, then $L_1(\mathfrak{g})$ has a unique nontrivial module $W$, and $\mathbb{D}_{0,n}(L_1(\mathfrak{g}), W^{\otimes n})$ is F-ample. The following corollary, a direct consequence of \cref{lem:fusion}, proves that this is ample without assuming the F-conjecture.

\begin{cor}\label{cor:f4g2}
    Let $\mathfrak{g}$ be $F_2$ or $G_4$ and $W$ be the nontrivial simple module of $L_1(\mathfrak{g})$. Then
    \[ \mathbb{D}_{0,n}(L_1(\mathfrak{g}), W^{\otimes n})=c\cdot\mathbb{D}_{0,n}(\vir{5}, W_2^{\otimes n}) \]
    for some $c<0$. In particular, they are ample.
\end{cor}

Note that the ampleness assertion can be proved by the same method as \cref{thm:VOAnefcw} by using \cite[Theorem 1.5]{Muk16} instead of $\vir{5}$.

\subsection{Critical level of Virasoro conformal block divisors}\label{subsec:crit}

As we mentioned in the remark after \cref{thm:virdeg}, a distinguishing property of conformal block divisors of $\vir{2k+1}$ is that their degree on $\M{0}{4}$ and $\M{1}{1}$ is entirely determined by their rank. In this section, by utilizing this property, we will show that on genus $0$, Virasoro conformal block divisors stabilize when we fix the modules and increase $k$. This stabilization is analogous to the critical level phenomena of conformal block divisors of affine VOAs discussed in \cite{Fak12, BGM15, BGM16}.

\begin{defn}\label{defn:critlev}
    Let $a_1, \cdots, a_n$ be a sequence of positive integers. The \textbf{critical level} $l(a_1, \cdots, a_n)$ of this tuple (with respect to Virasoro algebras) is the smallest integer $l$ which satisfies    
    \[ 2(l-1)\ge \sum_{i=1}^n (a_i-1). \]
\end{defn}

\begin{prop}\label{prop:rkstab}
    Let $a_1, \cdots, a_n$ be a sequence of positive integers, and let $W^i := W_{a_i}$. If $k \geq l(a_1, \cdots, a_n)$, then $\rank \mathbb{V}_{0,n}(\vir{2k+1}, W^\bullet)$ is independent of $k$. Furthermore, if $\sum_{i=1}^n (a_i - 1)$ is odd, then $\mathbb{V}_{0,n}(\vir{2k+1}, W^\bullet) = 0$.
\end{prop}

\begin{proof}
    We will prove a stronger statement. As a $\vir{2k+1}$-module,
    \[ \boxtimes_{i=1}^n W_{a_i}=\sum_{j}n_j W_j \]
    where $n_j \neq 0$ only if $j$ has the same parity as $\sum_{i=1}^n (a_i - 1) + 1$. Note that such an expression exists uniquely since $\left\{W_a \ |\ 1 \leq a \leq 2k, a \text{ odd}\right\}$ and $\left\{W_a \ |\ 1 \leq a \leq 2k, a \text{ even}\right\}$ are both representatives of simple $\vir{2k+1}$-modules. We will prove that if $2k-1\ge \sum_{i=1}^n (a_i-1)$ then this expression does not depend on $k$ and $n_j=0$ for $j>\sum_{i=1}^n (a_i - 1) + 1$. Indeed, this follows directly from the fusion rules of $\vir{2k+1}$ and induction. If $n=2$, then
    \[ W_{a_1}\boxtimes W_{a_2}=\sum_{i\stackrel{\text{2}}{=}a_2-a_1+1}^{a_1+a_2-1}W_i \]
    and this expression does not depend on $k$. In the induction step, the condition $2k-1\ge \sum_{i=1}^n (a_i-1)$ ensures that each of the tensor products satisfies this condition. Thus the resulting expression is independent of $k$, by the induction hypothesis. Since $\rank \mathbb{V}_{0,n}(\vir{2k+1}, W^\bullet)$ corresponds to the coefficient of the trivial module in this expression, it also does not depend on $k$. Note that we should be more careful in the odd case since in the even case, $W_1$ is the trivial module regardless of $k$, but in the odd case, $W_{2k}$ is the trivial module, and this varies as $k$ varies. This is why we need $2k-2\ge \sum_{i=1}^n (a_i-1)$ instead of $2k-1\ge \sum_{i=1}^n (a_i-1)$ in the statement.

    By the same reasoning as in \cref{lem:fusbd}, we have $n_j = 0$ for $j > \sum_{i=1}^n (a_i - 1) + 1$. Moreover, $j$ has the same parity as $\sum_{i=1}^n (a_i - 1) + 1 \leq 2k - 1$. Therefore, if $\sum_{i=1}^n (a_i - 1)$ is odd, there is no trivial module in the expression. This implies $\mathbb{V}_{0,n}(\vir{2k+1}, W^\bullet) = 0$.
\end{proof}

In summary, the fusion rules stabilize for a fixed set of modules as $k$ increases. Consequently, $\rank \mathbb{V}_{0,n}(\vir{2k+1}, W^\bullet)$ also stabilizes. Combined with \cref{thm:virdeg}, this further implies that the conformal block divisor $\mathbb{D}_{0,n}(\vir{2k+1}, W^\bullet)$ stabilizes as well.

\begin{thm}\label{thm:cbdstab}
    Under the same assumption, $\mathbb{D}_{0,n}(\vir{2k+1}, W^\bullet)$ is independent of $k$. Furthermore, if $\sum_{i=1}^n (a_i - 1)$ is odd, then $\mathbb{D}_{0,n}(\vir{2k+1}, W^\bullet) = 0$.
\end{thm}

\begin{proof}
    To prove this, it suffices to show that its intersection with $F$-curves is independent of $k$. Since any $F$-curve on $\M{0}{n}, $ is of Type 6, it is enough to show that
    \[ \sum_{M_1, M_2, M_3, M_4} \left( \prod_{i=1}^{4}\text{rank }\mathbb{V}_{g_i,|I_i|+1}(\vir{2k+1}, W^{I_i}\otimes M_i)\right) \cdot \mathbb{D}_{0,4}(\vir{2k+1}, \otimes_{i=1}^4 M_i) \]
    does not depend on $k$, where $I_1 \sqcup I_2 \sqcup I_3 \sqcup I_4 = [n]$ is a nonempty partition of $[n]$.
    
    First, the fusion rules stabilize by the proof of \cref{prop:rkstab}. Therefore, the list of modules $M_i$ that can appear in this expression does not depend on $k$, as it corresponds to the modules appearing in the expression $\boxtimes_{r \in I_i} W_{a_r} = \sum_{j} n_j W_j$, where $j$ has the same parity as $\sum_{r \in I_i} (a_r - 1) + 1$. This implies that $\text{rank } \mathbb{V}_{g_i,|I_i|+1}(\vir{2k+1}, W^{I_i} \otimes M_i)$ is also independent of $k$.

    Moreover, if $M_i = W_{b_i}$ is a module appearing in this expression, then $b_i \leq \sum_{r \in I_i} (a_r - 1) + 1$, and it has the same parity as $\sum_{r \in I_i} (a_r - 1) + 1$. Therefore, $W_{b_1}, \cdots, W_{b_4}$ also satisfy the conditions of \cref{prop:rkstab}. Hence, $\mathbb{D}_{0,4}(\vir{2k+1}, \otimes_{i=1}^4 M_i)$ does not depend on $k$. This completes the proof of the first assertion. The last assertion is a direct consequence of \cref{prop:rkstab}.
\end{proof}

\begin{rmk}
    Note that \cref{thm:cbdstab} is not a direct consequence of the explicit representation \cref{thm:c1cdf}. Therefore, \cref{thm:cbdstab} implies a lot of relations between tautological line bundles on $\M{0}{n}$.
\end{rmk}

\cref{thm:cbdstab} motivates the following definition.

\begin{defn}\label{defn:stcdb}
    Let $a_1, \cdots, a_n$ be a sequence of positive integers. The associated \textbf{stable Virasoro conformal block divisor} is defined as
    \[ \mathbb{D}_{0,n}(a_1,\cdots, a_n):=\mathbb{D}_{0,n}(\vir{2k+1}, \otimes_{i=1}^n W_{a_i} ) \]
    for any $k \geq l(a_1, \cdots, a_n)$. Note that this is zero if $\sum_{i=1}^n (a_i - 1)$ is odd.
\end{defn}

\begin{thm}\label{thm:stnef}
    $-\mathbb{D}_{0,n}(a_1,\cdots, a_n)$ is either a trivial or an ample line bundle.
\end{thm}

This theorem provides an infinite set of ample line bundles on $\M{0}{n}$. Note that there are many examples of nontrivial line bundles among $\mathbb{D}_{0,n}(a_1, \cdots, a_n)$ by \cref{thm:virfnef}. For example, if $n = 2m$ and $(a_1, \cdots, a_n) = (b_1, b_1, \cdots, b_m, b_m)$, then $-\mathbb{D}_{0,n}(a_1, \cdots, a_n)$ is nontrivial and hence ample.

\begin{prop}\label{thm:steff}
    Assume $a_i\ge 2$, $\sum_{i=1}^n (a_i - 1)$ is even and $2k - 1 = \sum_{i=1}^n (a_i - 1) + 1$, i.e., $k$ is exactly the critical level. Then $-\mathbb{D}_{0,n}(\vir{2k+1}, \otimes_{i=1}^n W_{a_i})$ is either zero or equal to $\sum_I c_I\delta_{0,I}$, where $c_I>0$.
\end{prop}

\begin{proof}
    In this proof, we use the notation and results of \cref{subsec:VOA}, \cref{subsec:cdf}, and \cref{subsec:nef}. If $\rank \mathbb{V}_{0,n}(\vir{2k+1}, \otimes_{i=1}^n W_{a_i})\leq 1$, then $\mathbb{D}_{0,n}(\vir{2k+1}, \otimes_{i=1}^n W_{a_i})=0$. Hence, we may assume that $r=\rank \mathbb{V}_{0,n}(\vir{2k+1}, \otimes_{i=1}^n W_{a_i})\geq 2$. By \cref{thm:c1cdf},
    \[ -\frac{1}{r}\mathbb{D}_{0,n}(\vir{2k+1}, \otimes_{i=1}^n W_{a_i})=\sum_{i=1}^n(- h_{a_i}) \psi_i-\sum_I b_I \delta_{0,I} \]
    where
    \[ b_I=\sum_{i=1}^{k} -h_i \cdot \frac{\rank \mathbb{V}_{0,|I|+1}(\vir{2k+1}, W_I\otimes W_i) \cdot \rank \mathbb{V}_{0,|I^c|+1}(\vir{2k+1}, W_{I^c}\otimes W_I)}{\rank \mathbb{V}_{0,n}(\vir{2k+1}, \otimes_{i=1}^n W_{a_i})}. \]
    Take the limit $k\to \infty$. Note that since $\otimes_{i=1}^n W_{a_i}$ is exactly at the critical level, by \cref{lem:fusbd}, for $W_i$'s such that $\rank \mathbb{V}_{0,|I|+1}(V, W_I\otimes W_i)$ and $\rank \mathbb{V}_{0,|I^c|+1}(V, W_{I^c}\otimes W_I)$ are both nonzero, $W_I\otimes W_i$ and $W_{I^c}\otimes W_I$ are also exactly at or above the critical level. Hence, by \cref{prop:rkstab}, all terms except the conformal weights stabilize after the critical level. Since $-h_i=\frac{(i-1)(2k-i)}{2(2k+1)}$, $-h_i$ converges to $\frac{i-1}{2}$. Therefore,
    \[-\frac{2}{r}\mathbb{D}_{0,n}(\vir{2k+1}, \otimes_{i=1}^n W_{a_i})=\sum_{i=1}^n (a_i-1) \psi_i-\sum_I b_I' \delta_{0,I} \]
    where
    \[ b_I'=\sum_{i=1}^{k} (i-1) \cdot \frac{\rank \mathbb{V}_{0,|I|+1}(\vir{2k+1}, W_I\otimes W_i) \cdot \rank \mathbb{V}_{0,|I^c|+1}(\vir{2k+1}, W_{I^c}\otimes W_I)}{\rank \mathbb{V}_{0,n}(\vir{2k+1}, \otimes_{i=1}^n W_{a_i})}. \]
    By taking a difference,
     \[-\frac{4k}{r}\mathbb{D}_{0,n}(\vir{2k+1}, \otimes_{i=1}^n W_{a_i})=\sum_{i=1}^n (a_i-1)(2k-a_i-1) \psi_i-\sum_I b_I'' \delta_{0,I} \]
    where
    \[ b_I''=\sum_{i=1}^{k} (i-1)(2k-i-1)\cdot \frac{\rank \mathbb{V}_{0,|I|+1}(\vir{2k+1}, W_I\otimes W_i) \cdot \rank \mathbb{V}_{0,|I^c|+1}(\vir{2k+1}, W_{I^c}\otimes W_I)}{\rank \mathbb{V}_{0,n}(\vir{2k+1}, \otimes_{i=1}^n W_{a_i})}. \]
    Now consider $\mathbb{D}_{0,n}(L_1(\mathfrak{sl}_{2k-2}), \otimes_{i=1}^n U_{a_i-1})$. The critical level of $(U_{a_1-1},\cdots, U_{a_n-1})$, in the sense of \cite{BGM15}, is $0$. Therefore, by \cite[Proposition 1.3]{BGM15}, this divisor is $0$. Hence, by \cref{thm:c1cdf},
    \begin{align*}
        &(4k-4)\cdot \mathbb{D}_{0,n}(L_1(\mathfrak{sl}_{2k-2}), \otimes_{i=1}^n U_{a_i-1})\\
        &=\sum_{i=1}^n (a_i-1)(2k-1-a_i)\psi_i-\sum_{I} \left(\sum_{j\in I}(a_j-1)\right)\left(2k-2-\sum_{j\in I}(a_j-1) \right) \delta_{0,I}=0.
    \end{align*}
    By comparing this with the preceding expression of $\mathbb{D}_{0,n}(\vir{2k+1}, \otimes_{i=1}^n W_{a_i})$, to complete the proof, it is enough to show that
    \begin{align*}
    &\sum_{i=1}^{k} (i-1)(2k-i-1)\cdot \frac{\rank \mathbb{V}_{0,|I|+1}(\vir{2k+1}, W_I\otimes W_i) \cdot \rank \mathbb{V}_{0,|I^c|+1}(\vir{2k+1}, W_{I^c}\otimes W_I)}{\rank \mathbb{V}_{0,n}(\vir{2k+1}, \otimes_{i=1}^n W_{a_i})}\\&<\left(\sum_{j\in I}(a_j-1)\right)\left(2k-2-\sum_{j\in I}(a_j-1) \right).     
    \end{align*}
    Since $\sum_{j\in I}(a_j-1) = 2k-2-\sum_{j\in I^c}(a_j-1)$, we may replace $I$ with $I^c$ and assume $\sum_{j\in I}(a_j-1) \leq k-1$. By \cref{lem:fusbd}, if $\rank \mathbb{V}_{0,|I|+1}(\vir{2k+1}, W_I\otimes W_i) \ne 0$, then $i \leq \sum_{j\in I}(a_j-1)+1$. Hence, $(i-1)(2k-i-1)\le \left(\sum_{j\in I}(a_j-1)\right)\left(2k-2-\sum_{j\in I}(a_j-1)\right)$. Since the first term is a weighted average of $h_i$'s, it should be less than or equal to the second term. Moreover, since $\sum_{j\in I}(a_j-1) \leq k-1$, by the fusion rules, it is straightforward to see that the multiplicity of $W_{\sum_{j\in I}(a_j-1)+1}$ is $1$ in $\otimes_{j\in I}W_{a_j}$. Since we assumed that $\rank \mathbb{V}_{0,n}(\vir{2k+1}, \otimes_{i=1}^n W_{a_i}) \geq 2$, there exists an $i < \sum_{j\in I}(a_j-1)$ that contributes to the first term. Hence, the second term is strictly larger than the first term.
\end{proof}

\begin{proof}[proof of \cref{thm:stnef}]
    We will use induction on $n$. The base case is evident by \cref{cor:fconj}. We will use \cref{lem:effnef} for the induction step. By \cref{thm:steff}, it suffices to show that if $\mathbb{D}_{0,n}(a_1, \dots, a_n)$ is nontrivial, then its restriction to boundary divisors is ample. We apply \cref{thm:factor} to handle the boundary restrictions. Note that, by \cref{lem:fusbd}, the restriction of stable Virasoro conformal block divisors is also an effective sum of stable Virasoro conformal block divisors. Therefore, the remainder of the proof follows the same argument as in the proof of \cref{thm:VOAnefcw}.
\end{proof}

\subsection{Differences of Virasoro conformal block divisors}\label{subsec:diff}

Unfortunately, the Virasoro conformal block divisors are pullbacks of F-ample divisors along projection maps. Hence, they are not contained in the interior of an interesting face of the nef cone. Here, we will examine the differences between Virasoro conformal block divisors and demonstrate that many are F-nef, and moreover, some of them are non-F-ample, nonzero F-nef divisors. In particular, we will prove \cref{thm:main3}. We begin with the first assertion, which is relatively straightforward to prove.

\begin{thm}\label{thm:criteven}
    Let $a_1, \cdots, a_n$ be integers such that $1 \leq a_i \leq 2k$ and $\sum_{i=1}^{n}(a_i - 1)$ is even. Then
    \[
    \mathbb{D}_{0,n}\left(\vir{2k+1}, \otimes_{i=1}^{n}W_{a_i} \right) - \mathbb{D}_{0,n}\left(\vir{2k+3}, \otimes_{i=1}^{n}W_{a_i} \right)
    \]
    is F-nef. 
\end{thm}

To prove this, we introduce a notation. Let $1 \leq a_1, \cdots, a_n \leq 2k$ be integers between $1$ and $2k$. Define
\begin{align*}
     &M_k(a_1,\cdots, a_n):=\\&\left\{b\ |\ 1\le b\le 2k,\ \sum_{i=1}^n (a_i-1)+(b-1)\text{ is even and } \rank \mathbb{V}_{0,n+1}\left(\vir{2k+1}, \otimes_{i=1}^{n}W_{a_i}\otimes W_b \right)\ge 1 \right\}. 
\end{align*}

\begin{prop}\label{prop:evenrank}
    Let $a_1, \cdots, a_n$ be integers such that $1 \leq a_i \leq 2k$ and $\sum_{i=1}^{4}(a_i - 1)$ is even. Then 
    \[\rank \mathbb{V}_{0,n}\left(\vir{2k+1}, \otimes_{i=1}^{n}W_{a_i} \right) \leq \rank \mathbb{V}_{0,n}\left(\vir{2k+3}, \otimes_{i=1}^{n}W_{a_i} \right)\]
    for any $n \geq 3$. In particular, for any $1 \leq a_1, \cdots, a_n \leq 2k$ (without the assumption of parity),
    \[M_k(a_1, \cdots, a_n) \subseteq M_{k+1}(a_1, \cdots, a_n).\]
\end{prop}

\begin{proof}
   Note that the second assertion directly follows from the first assertion. We will proceed by induction on $n$. For the base case, when $n = 3$, the statement directly follows from the fusion rules of $\vir{2k+1}$ and $\vir{2k+3}$. In particular, the second assertion holds for $n = 2$.

    Assume the statement holds for $n$, and consider the case of $n+1$. By \cref{thm:factor}, we have:
    \begin{align*}
    &\rank \mathbb{V}_{0,n+1}\left(\vir{2k+1}, \otimes_{i=1}^{n+1}W_{a_i} \right) = \\
    &\sum_{b \in M_k(a_n, a_{n+1})} \rank \mathbb{V}_{0,n}\left(\vir{2k+1}, \otimes_{i=1}^{n-1}W_{a_i} \otimes W_b \right) \cdot \rank \mathbb{V}_{0,3}\left(\vir{2k+1}, W_{a_{n}} \otimes W_{a_{n+1}} \otimes W_b \right), \\
    &\rank \mathbb{V}_{0,n+1}\left(\vir{2k+3}, \otimes_{i=1}^{n+1}W_{a_i} \right) = \\
    &\sum_{b \in M_{k+1}(a_n, a_{n+1})} \rank \mathbb{V}_{0,n}\left(\vir{2k+3}, \otimes_{i=1}^{n-1}W_{a_i} \otimes W_b \right) \cdot \rank \mathbb{V}_{0,3}\left(\vir{2k+3}, W_{a_{n}} \otimes W_{a_{n+1}} \otimes W_b \right).
    \end{align*}

    Since $M_k(a_n, a_{n+1}) \subseteq M_{k+1}(a_n, a_{n+1})$, it follows by induction that 
    \[\rank \mathbb{V}_{0,n+1}\left(\vir{2k+1}, \otimes_{i=1}^{n+1}W_{a_i} \right) \leq \rank \mathbb{V}_{0,n+1}\left(\vir{2k+3}, \otimes_{i=1}^{n+1}W_{a_i} \right).\]
\end{proof}

\begin{proof}[proof of \cref{thm:criteven}]
    By the same method as in the proof of \cref{lem:VOAGKM}, 
    \begin{align*}
    &\mathbb{D}_{0,n}\left(\vir{2k+1}, \otimes_{i=1}^{n}W_{a_i} \right)\cdot F_{I_1, I_2, I_3, I_4} = \\
    &\sum_{b_i \in M_k(a_{I_i})} \left( \prod_{i=1}^{4}\text{rank }\mathbb{V}_{0,|I_i|+1}(\vir{2k+1}, W^{I_i}\otimes W_{b_i})\right) \cdot \mathbb{D}_{0,4}(\vir{2k+1}, \otimes_{i=1}^4 W_{b_i}),
    \end{align*}
    and we have a similar expression for $\mathbb{D}_{0,n}\left(\vir{2k+3}, \otimes_{i=1}^{n}W_{a_i} \right)\cdot F_{I_1, I_2, I_3, I_4}$. By \cref{prop:evenrank} and \cref{thm:virdeg},
    \begin{enumerate}
    \item $M_k(a_{I_i}) \subseteq M_{k+1}(a_{I_i})$,
    \item $\rank \mathbb{V}_{0,n}\left(\vir{2k+1}, \otimes_{i=1}^{n}W_{a_i} \right) \leq \rank \mathbb{V}_{0,n}\left(\vir{2k+3}, \otimes_{i=1}^{n}W_{a_i} \right)$, and
    \item $\mathbb{D}_{0,4}(\vir{2k+1}, \otimes_{i=1}^4 W_{b_i}) \ge \mathbb{D}_{0,4}(\vir{2k+3}, \otimes_{i=1}^4 W_{b_i})$.        
    \end{enumerate}
    Note that $\sum_{i=1}^{4}(b_i-1)$ and $\sum_{i\in I_j}(a_i-1) + (b_j-1)$ are even. This implies
    \[\mathbb{D}_{0,n}\left(\vir{2k+1}, \otimes_{i=1}^{n}W_{a_i} \right)\cdot F_{I_1, I_2, I_3, I_4} \geq \mathbb{D}_{0,n}\left(\vir{2k+3}, \otimes_{i=1}^{n}W_{a_i} \right)\cdot F_{I_1, I_2, I_3, I_4}.\]
    Hence, $\mathbb{D}_{0,n}\left(\vir{2k+1}, \otimes_{i=1}^{n}W_{a_i} \right) - \mathbb{D}_{0,n}\left(\vir{2k+3}, \otimes_{i=1}^{n}W_{a_i} \right)$ is F-nef.
\end{proof}

Now we will see some examples. All of them are computed by \cite[Difference of Virasoro CBD]{Choigit24}.

\begin{eg}\label{eg:even1}
   If $(a_1, \cdots, a_6) = (3, 4, 5, 6, 6, 6)$, then the critical level is $l = 13$ and 
    \[D = \mathbb{D}_{0,6}(\text{Vir}_{2,13}, W_3 \otimes W_4 \otimes W_5 \otimes W_6^{\otimes 3}) - \mathbb{D}_{0,6}(\text{Vir}_{2,15}, W_3 \otimes W_4 \otimes W_5 \otimes W_6^{\otimes 3})\]
    is F-nef. Indeed, since the F-conjecture holds for $\M{0}{6}$, $D$ is nef. Moreover,
    \[D \cdot F_{\left\{1,5,6\right\}, 2,3,4} = 0,\text{ and } D \cdot F_{\left\{1,2,3\right\}, 4,5,6} \ne 0,\]
    $D$ is a non-ample, non-zero nef divisor.
\end{eg}

\begin{eg}\label{eg:even2}
    If $(a_1, \cdots, a_8) = (2, 3, 3, 4, 4, 5, 5, 6)$, then the critical level is $l = 13$ and 
    \[D = \mathbb{D}_{0,8}(\text{Vir}_{2,13}, W_2\otimes W_3^{\otimes 2} \otimes W_4^{\otimes 2} \otimes W_5^{\otimes 2} \otimes W_6) - \mathbb{D}_{0,8}(\text{Vir}_{2,15}, W_2\otimes W_3^{\otimes 2} \otimes W_4^{\otimes 2} \otimes W_5^{\otimes 2} \otimes W_6)\]
    is F-nef. By explicit computer verification \cite[effectivity]{Choigit24}, we can show that $D$ is an effective sum of boundaries. Hence, by \cref{lem:effnef}, $D$ is nef. Moreover,
    \[D \cdot F_{\left\{1,2,3,4,6\right\}, 5, 7, 8} = 0, \text{ and } D\cdot F_{\left\{1,2,3,4,5\right\}, 6, 7, 8} \ne 0,\]
    hence $D$ is a non-ample, non-zero nef divisor. This gives a desired new example of a nef divisor on $\M{0}{n}$ using the difference of Virasoro conformal block divisors.

    Note that $D$ is not a pullback of a divisor along the projection map since the intersection of $D$ with
    \[ F_{1,\left\{2,3,4,5\right\}, \left\{6,7\right\}, 8}, F_{\left\{1,3,4\right\},2, \left\{5,6,7\right\}, 8}, F_{\left\{1,2,4\right\},3, \left\{5,6,7\right\}, 8}, F_{\left\{1,2,3,7\right\}, 4, \left\{5,6\right\}, 8}, F_{\left\{1,2,3,7\right\}, \left\{4,6\right\}, 5, 8}  \]
    are also nonzero.

\end{eg}

Unfortunately, it seems very hard to prove that the divisors in \cref{thm:criteven} are nef using our current methodology. The problem is that we need to prove that these divisors are a sum of an effective divisor and the canonical divisor. However, since we are using two Virasoro VOAs, it is very difficult to use \cref{thm:c1cdf}. In particular, the fusion rules and conformal weights vary as the VOA varies, which makes the estimation of coefficients very hard.

Now we consider the case where $\sum_{i=1}^{n}(a_i-1)$ is odd. This case is more complicated.

\begin{thm}\label{thm:critodd}
    Let $a_1, \cdots, a_n$ be integers such that $1 \leq a_i \leq 2k$ and $\sum_{i=1}^{n}(a_i - 1)$ is even. 
    \begin{enumerate}
        \item If $\sum_{i=1}^{n}(a_i - 1) \leq 2k-1$, or equivalently, $l(a_1,\cdots, a_n)-1\le k$ then
        \[\mathbb{D}_{0,n}\left(\vir{2k+1}, \otimes_{i=1}^{n}W_{a_i} \right) = 0.\]
        Consequently, if $a_i\le 2k-2$ for every $i$,
        \[\mathbb{D}_{0,n}\left(\vir{2k+1}, \otimes_{i=1}^{n}W_{a_i} \right) - \mathbb{D}_{0,n}\left(\vir{2k-1}, \otimes_{i=1}^{n}W_{a_i} \right)\]
        is F-ample or zero.
        \item If $\sum_{i=1}^{n}(a_i - 1) \leq 2k+3$, or equivalently, $l(a_1,\cdots, a_n)-3\le k$, $a_i>1$ and $n\ge 5$, then
        \[\mathbb{D}_{0,n}\left(\vir{2k+3}, \otimes_{i=1}^{n}W_{a_i} \right) - \mathbb{D}_{0,n}\left(\vir{2k+1}, \otimes_{i=1}^{n}W_{a_i} \right)\]
        is F-nef.
    \end{enumerate}
\end{thm}

\cref{thm:critodd} is more subtle to prove than \cref{thm:criteven}. Note that the sign of \cref{thm:criteven} and \cref{thm:critodd} are opposite. The reason for this is the following. Let $1 \leq a_1 \leq a_2 \leq a_3 \leq a_4 \leq 2k$ be integers such that their sum is odd. Assume that $a_2 + a_3 + a_4 - a_1 \leq 2k+1$. Then $a_1+a_4, a_2+a_3\le 2k+1$, so
\[W_{a_1}\boxtimes W_{a_4} = \sum_{i \stackrel{\text{2}}{=} a_4-a_1+1}^{a_1+a_4-1}W_i, \quad W_{a_2}\boxtimes W_{a_3} = \sum_{i \stackrel{\text{2}}{=} a_3-a_2+1}^{a_2+a_3-1}W_i = \sum_{i \stackrel{\text{2}}{=} 2k+2-(a_2+a_3)}^{2k+a_2-a_3}W_i.\]
By computing the number of common summands, we observe that $a_2 + a_3 + a_4 - a_1 \leq 2k+1$ implies $\rank \mathbb{V}_{0,4}(\vir{2k+1}, \otimes_{i=1}^4 W_{a_i})$ decreases as $k$ increases. Therefore, as opposed to \cref{prop:evenrank}, the rank decreases in this case. However, this requires the assumption $a_2 + a_3 + a_4 - a_1 \leq 2k+1$, and since $M_k(a_1, \cdots, a_n) \subseteq M_{k+1}(a_1, \cdots, a_n)$ holds regardless of parity, we need further assumptions to ensure $M_k(a_1, \cdots, a_n) = M_{k+1}(a_1, \cdots, a_n)$. This makes the proof of \cref{thm:critodd} more involved and the conclusion much weaker. The analysis we need is summarized in the following proposition.

\begin{prop}\label{prop:odddeg}
    Let $1\le a_1\le a_2\le a_3\le a_4\le 2k$ be integers such that $\sum_{i=1}^4(a_i-1)$ is odd.
    \begin{enumerate}
        \item If $\sum_{i=1}^4(a_i-1)\le 2k-1$, then $\mathbb{D}_{0,4}(\vir{2k+1}, \otimes_{i=1}^4 W_{a_i})=0$.
        \item If $\sum_{i=1}^4(a_i-1)\le 2k+3$, then $\mathbb{D}_{0,4}(\vir{2k+3}, \otimes_{i=1}^4 W_{a_i})\ge \mathbb{D}_{0,4}(\vir{2k+1}, \otimes_{i=1}^4 W_{a_i})$ except the case of $a_1=a_2=2$ and $a_3+a_4=2k+3$. In this case,
        \[ \mathbb{D}_{0,4}(\vir{2k+3}, \otimes_{i=1}^4 W_{a_i})=-1, \ \mathbb{D}_{0,4}(\vir{2k+1},\ \otimes_{i=1}^4 W_{a_i})=0. \]
    \end{enumerate}
\end{prop}

\begin{proof}
    (1) $\sum_{i=1}^4(a_i-1)\le 2k-1$ implies $a_i+a_j\le 2k+1$ for any $i\ne j$. Therefore, as above, 
    \[W_{a_1}\boxtimes W_{a_4} = \sum_{i \stackrel{\text{2}}{=} a_4-a_1+1}^{a_1+a_4-1}W_i \text{  and  } W_{a_2}\boxtimes W_{a_3} = \sum_{i \stackrel{\text{2}}{=} 2k+2-(a_2+a_3)}^{2k+a_2-a_3}W_i.\]
    Since $\sum_{i=1}^4(a_i-1)\le 2k-1$, $\rank \mathbb{V}_{0,4}(\vir{2k+1}, \otimes_{i=1}^4 W_{a_i})\le 1$. Therefore, 
    \[\mathbb{D}_{0,4}(\vir{2k+1}, \otimes_{i=1}^4 W_{a_i})=0.\]

    (2) We may assume that $1<a_i$ for each $i$, since otherwise this holds trivially by \cref{thm:propvac}. Moreover, if $\sum_{i=1}^4(a_i-1)\le 2k+1$, then this is evident by (1) and \cref{thm:virdeg}. Hence, we may assume that $\sum_{i=1}^4(a_i-1)= 2k+3$.  

    First, consider the case of $(a_1,a_2)=(2,2)$ or $(2,3)$. If $(a_1,a_2)=(2,2)$, then $a_3+a_4=2k+3$. Hence, 
    \[ W_{a_3}\boxtimes W_{a_4}=\sum_{i \stackrel{\text{2}}{=} a_4-a_3+1}^{2k-2}W_i \text{ for }\vir{2k+1},\ W_{a_3}\boxtimes W_{a_4}=\sum_{i \stackrel{\text{2}}{=} a_4-a_3+1}^{2k+2}W_i \text{ for }\vir{2k+3},  \] 
    so $\mathbb{D}_{0,4}(\vir{2k+3}, \otimes_{i=1}^4 W_{a_i})=-1$, $\mathbb{D}_{0,4}(\vir{2k+1},\ \otimes_{i=1}^4 W_{a_i})=0$.
    If $(a_1,a_2)=(2,3)$, then $a_3+a_4=2k+2$, so 
    \[ W_{a_3}\boxtimes W_{a_4}=\sum_{i \stackrel{\text{2}}{=} a_4-a_3+1}^{2k-1}W_i \text{ for }\vir{2k+1},\ W_{a_3}\boxtimes W_{a_4}=\sum_{i \stackrel{\text{2}}{=} a_4-a_3+1}^{2k+1}W_i \text{ for }\vir{2k+3}. \]
    This implies $\mathbb{D}_{0,4}(\vir{2k+3}, \otimes_{i=1}^4 W_{a_i})=\mathbb{D}_{0,4}(\vir{2k+1}, \otimes_{i=1}^4 W_{a_i})=-1$.
    
    The remaining cases satisfy $a_i+a_j\le 2k+1$ for any $i\ne j$. Observe that $a_4-a_2+1\le 2k+2-(a_1+a_3)$. Since we have
    \[W_{a_2}\boxtimes W_{a_4} = \sum_{i \stackrel{\text{2}}{=} a_4-a_2+1}^{a_2+a_4-1}W_i\text{  and  } W_{a_1}\boxtimes W_{a_3} = \sum_{i \stackrel{\text{2}}{=} 2k+2-(a_1+a_3)}^{2k+a_1-a_3}W_i,\]
    as $k$ increases, $\rank \mathbb{V}_{0,4}(\vir{2k+1}, \otimes_{i=1}^4 W_{a_i})$ decreases, so $\mathbb{D}_{0,4}(\vir{2k+1}, \otimes_{i=1}^4 W_{a_i})$ decreases.
\end{proof}

\begin{proof}[proof of \cref{thm:critodd}]
    (1) As before, 
    \begin{align*}
    &\mathbb{D}_{0,n}\left(\vir{2k+1}, \otimes_{i=1}^{n}W_{a_i} \right)\cdot F_{I_1, I_2, I_3, I_4} = \\
    &\sum_{b_i \in M_k(a_{I_i})} \left( \prod_{i=1}^{4}\text{rank }\mathbb{V}_{0,|I_i|+1}(\vir{2k+1}, W^{I_i}\otimes W_{b_i})\right) \cdot \mathbb{D}_{0,4}(\vir{2k+1}, \otimes_{i=1}^4 W_{b_i}).
    \end{align*}
    For $b_i \in M_k(a_{I_i})$, by \cref{lem:fusbd}, $\sum_{i=1}^4(b_i-1)\le 2k-1$. Hence we can apply \cref{prop:odddeg} and obtain that any such intersection is zero. Since the F-curves generate $\text{A}_1(\M{0}{n})$, this implies the conclusion.

    (2) By the same argument as \cref{prop:odddeg}, we may assume that $\sum_{i=1}^{n}(a_i-1)=2k+3$. We will use the same expression as above. First, consider the case that there exists $b_i' \in M_{k+1}(a_{I_i})$ for each $i$ such that
    \[(b_1', b_2', b_3', b_4') = (2, 2, c, d)\text{ up to permutation and }c+d=2k+3.\]
    We may assume that $(b_1', b_2', b_3', b_4') = (2, 2, c, d)$. Since $\sum_{i=1}^4(b_i'-1)=2k+3$, we know that, by \cref{lem:fusbd}, $b_i'$ is the maximum of $M_{k+1}(a_{I_i})$. In particular, $I_1, I_2$ are singletons, and we may assume that $I_4$ is not a singleton since $n\ge 5$.

    First, consider the case of $c, d < 2k+1$. Then by the proof of \cref{prop:rkstab}, $M_k(a_{I_i}) = M_{k+1}(a_{I_i})$ for each $i$ and
    \[\mathbb{V}_{0,|I_i|+1}(\vir{2k+1}, W^{I_i}\otimes W_{b_i}) = \mathbb{V}_{0,|I_i|+1}(\vir{2k+3}, W^{I_i}\otimes W_{b_i})\]
    for each $i$ and $b_i \in M_k(a_{I_i})$. Therefore, it is enough to prove that 
    \begin{align*}
    &\sum_{b_i \in M_k(a_{I_i})} \left( \prod_{i=1}^{4}\text{rank }\mathbb{V}_{0,|I_i|+1}(\vir{2k+1}, W^{I_i}\otimes W_{b_i})\right) \cdot \mathbb{D}_{0,4}(\vir{2k+1}, \otimes_{i=1}^4 W_{b_i})\\
    &\le \sum_{b_i \in M_k(a_{I_i})} \left( \prod_{i=1}^{4}\text{rank }\mathbb{V}_{0,|I_i|+1}(\vir{2k+1}, W^{I_i}\otimes W_{b_i})\right) \cdot \mathbb{D}_{0,4}(\vir{2k+3}, \otimes_{i=1}^4 W_{b_i}).
    \end{align*}
    By \cref{prop:odddeg}, apart from the case of $(b_1, b_2, b_3, b_4) = (2, 2, c, d)$, each term with fixed $(b_1, b_2, b_3, b_4)$ of the LHS is less than or equal to the corresponding term of the RHS. Hence, the only problematic term is $(b_1, b_2, b_3, b_4) = (2, 2, c, d)$. Since $b_i'$ is the maximum of $M_{k+1}(a_{I_i})$ for each $i$, $\text{rank }\mathbb{V}_{0,|I_i|+1}(\vir{2k+1}, W^{I_i}\otimes W_{b_i'}) = 1$ for each $i$, by the fusion rules. Hence, the difference between the LHS and RHS of the terms corresponding to $(b_1, b_2, b_3, b_4) = (2, 2, c, d)$ is $1$. Therefore, it is enough to show that $(b_1, b_2, b_3, b_4)$ exists such that the corresponding term of the RHS is strictly larger than that of the LHS. Since $I_4$ is not a singleton, by \cref{lem:virfus}, $(b_1, b_2, b_3, b_4) = (2, 2, c, d-2)$ also contributes to the summation, and this term satisfies the requirement. This completes the proof for the first case.
    
    Now we will consider the case of $(b_1', b_2', b_3', b_4') = (2, 2, 2, 2k+1)$. In this case, $I_1, I_2, I_3$ are singletons. Hence, we need to prove that
    \begin{align*}
    &\sum_{b_4 \in M_k(a_{I_4})} \text{rank }\mathbb{V}_{0,|I_4|+1}(\vir{2k+1}, W^{I_4}\otimes W_{b_4}) \cdot \mathbb{D}_{0,4}(\vir{2k+1}, W_2^{\otimes 3}\otimes W_{b_4})\\
    &\le \sum_{b_4 \in M_{k+1}(a_{I_4})} \text{rank }\mathbb{V}_{0,|I_4|+1}(\vir{2k+3}, W^{I_4}\otimes W_{b_4}) \cdot \mathbb{D}_{0,4}(\vir{2k+3}, W_2^{\otimes 3}\otimes W_{b_4}).
    \end{align*}    
    Let $j \in I_4$. Then $\sum_{i \in I_4 \setminus \left\{j\right\}} (a_i-1) \le 2k-1$. Hence, by the proof of \cref{prop:rkstab}, the representation
    \[\boxtimes_{i \in I_4 \setminus \left\{j\right\}} W_{a_i} = \sum_{s} n_s W_s\] 
    where $n_s \neq 0$ only if $s$ has the same parity as $\sum_{i \in I_4 \setminus \left\{j\right\}} (a_i - 1) + 1$, does not depend on whether we consider them as a module over $\vir{2k+1}$ or $\vir{2k+3}$. Also, if $m$ is the maximum index $s$ such that $n_s \ne 0$, then only $W_m \boxtimes W_{a_j}$ depends on the choice of algebra between $\vir{2k+1}$ and $\vir{2k+3}$. Therefore, we obtain $M_{k+1}(a_{I_4}) = M_k(a_{I_4}) \cup \left\{2k+1\right\}$, $\rank\mathbb{V}_{0,|I_4|+1}(\vir{2k+1}, W^{I_4}\otimes W_{b_4}) = \rank\mathbb{V}_{0,|I_4|+1}(\vir{2k+3}, W^{I_4}\otimes W_{b_4})$ for $b_4 \ne 2k+1$, and $\rank\mathbb{V}_{0,|I_4|+1}(\vir{2k+3}, W^{I_4}\otimes W_{2k+1}) = 1$.

    Hence, the only problematic term is $b_4 = 2k+1$. In this case, the LHS has no corresponding term, and the contribution to the RHS is $-1$. Therefore, as before, it is enough to show that there exists another pair $(2, 2, 2, b_4)$ such that the corresponding term in the RHS is strictly larger than the corresponding term in the LHS. Since $I_4$ is not a singleton, $2k-1 \in M_k(a_{I_4})$ by \cref{lem:virfus}. It is straightforward to see that $(2, 2, 2, 2k-1)$ satisfies the requirement.

    Now we can assume that such $(b_1', b_2', b_3', b_4')$ does not exist. Hence, we have $\sum_{j \in I_i}(a_j - 1) \le 2(k - 1)$ for each $i$. By the proof of \cref{prop:rkstab}, $M_k(a_{I_i}) = M_{k+1}(a_{I_i})$ and $\rank\mathbb{V}_{0,|I_i|+1}(\vir{2k+1}, W^{I_i}\otimes W_{b_i}) = \rank\mathbb{V}_{0,|I_i|+1}(\vir{2k+3}, W^{I_i}\otimes W_{b_i})$. Hence, it is enough to show that
    \[\mathbb{D}_{0,4}(\vir{2k+3}, \otimes_{i=1}^4 W_{b_i}) \ge \mathbb{D}_{0,4}(\vir{2k+1}, \otimes_{i=1}^4 W_{b_i})\]
    for each $b_i \in M_k(a_{I_i})$. However, $\sum_{i=1}^{4}(b_i - 1) \le 2k+3$, and the case of $(b_1, b_2, b_3, b_4) = (2, 2, c, d)$, $c + d = 2k+3$, is already covered in the previous part. Thus, this follows from \cref{prop:odddeg}.
\end{proof}

\begin{eg}
    Let $(a_1,a_2,a_3,a_4,a_5)=(2,2,4,5,5)$. Then, its critical level is $l=8$. Let
    \[ D=\mathbb{D}_{0,5}(\text{Vir}_{2,13}, W_2^{\otimes 2}\otimes W_4\otimes W_5^{\otimes 2})-\mathbb{D}_{0,5}(\text{Vir}_{2,11}, W_2^{\otimes 2}\otimes W_4\otimes W_5^{\otimes 2}).\]
    Then, $D$ is an F-nef (hence also nef) divisor. Moreover, 
    \[ D\cdot F_{1, 2, \left\{3,4\right\},5}=0, D\cdot F_{\left\{1,2\right\},3,4,5}=2 \]
    hence $D$ is a nonzero, non-ample nef divisor. Note that this is not the pullback of an ample divisor on $\M{0}{4}$. If it were, then since $D\cdot F_{\left\{1,2\right\},3,4,5}=2$, it would have to be a pullback along either $\pi_1$ or $\pi_2$. By symmetry, it would then be the pullback of a divisor along both $\pi_1$ and $\pi_2$, which is impossible.
\end{eg}

Unfortunately, unlike \cref{thm:criteven}, \cref{thm:critodd} is much weaker and does not provide a lot of interesting divisors. Also, \cref{thm:critodd} is not true without the assumption on $k$ since
\[ \mathbb{D}_{0,6}(\text{Vir}_{2,15}, W_5^{\otimes 3}\otimes W_6^{\otimes 3})\cdot F_{\left\{1,2\right\}, \left\{3,6\right\},4,5}<\mathbb{D}_{0,6}(\text{Vir}_{2,13}, W_5^{\otimes 3}\otimes W_6^{\otimes 3})\cdot F_{\left\{1,2\right\}, \left\{3,6\right\},4,5}. \]
Even though, by a brute force experiment using \cite[Difference of Virasoro CBD]{Choigit24}, it seems like a much stronger statement than \cref{thm:critodd} should hold. In particular, \cref{thm:critodd} (2) should be true even when $k = l(a_1, \cdots, a_n) - 4$. As an example,
\begin{align*}
    &\mathbb{D}_{0,5}(\text{Vir}_{2,13}, W_2\otimes W_3\otimes W_5^{\otimes 3})-\mathbb{D}_{0,5}(\text{Vir}_{2,11}, W_2\otimes W_3\otimes W_5^{\otimes 3})\text{ and }\\
    &\mathbb{D}_{0,6}(\text{Vir}_{2,13}, W_2^{\otimes 3}\otimes W_5^{\otimes 3})-\mathbb{D}_{0,6}(\text{Vir}_{2,11}, W_2^{\otimes 3}\otimes W_5^{\otimes 3})
\end{align*}
are nonzero, nonample nef divisors. However, the combinatorial analysis that we need to show \cref{thm:critodd} (2) for $k = l(a_1, \cdots, a_n) - 4$ seems too complicated, so we decided to stop here. Instead, we suggest the following conjecture:

\begin{conj}\label{conj:odddiff}
    For any $r > 0$, there exists $N > 0$ such that for any $n \ge N$,
    \[ \mathbb{D}_{0,n}\left(\vir{2k+3}, \otimes_{i=1}^{n}W_{a_i} \right) - \mathbb{D}_{0,n}\left(\vir{2k+1}, \otimes_{i=1}^{n}W_{a_i} \right) \]
    is F-nef if $2 \le a_1, \cdots, a_n \le 2k$, $\sum_{i=1}^{n}(a_i-1)$ is odd and $l(a_1, \cdots, a_n) - r \le k$.
\end{conj}

The rationale for this conjecture is as follows. Concerning the proof of \cref{thm:critodd}, we noticed that, in most of the terms that occur when intersecting this divisor with an F-curve, the contributions are positive. However, there are some exceptions whose contributions are negative. Nonetheless, if there are sufficiently many modules, then, due to \cref{lem:virfus}, significantly more terms should contribute positively.

We have mentioned that the stabilization of the Virasoro conformal block divisor can be viewed as an analogue of critical level phenomena of coinvariant divisors of affine VOAs. In the case of affine VOAs, the critical level coinvariant divisors are usually nonzero and have interesting properties, such as corresponding to interesting morphisms \cite[Theorem 4.5]{Fak12} and duality \cite{BGM15}. However, Virasoro conformal block divisors typically stabilize earlier than the critical level, and the exact point of stabilization seems to depend significantly on $(a_1, \cdots, a_n)$. \cref{thm:critodd} (1) shows that if $\sum_{i=1}^{n}(a_i - 1)$ is odd, then $\mathbb{D}_{0,n}\left(\vir{2k+1}, \otimes_{i=1}^{n}W_{a_i} \right)$ stabilizes at $k = l(a_1, \cdots, a_n) - 1$, i.e., one below the critical level. This is the best possible since 
\[\mathbb{D}_{0,5}(\text{Vir}_{2,13}, W_2^{\otimes 2}\otimes W_4\otimes W_5^{\otimes 2})\cdot F_{1, 2, \left\{3,4\right\},5}\ne 0.\]
However, if $\sum_{i=1}^{n}(a_i - 1)$ is even, direct experiments suggest that $\mathbb{D}_{0,n}\left(\vir{2k+1}, \otimes_{i=1}^{n}W_{a_i} \right)$ stabilizes much earlier than the critical level. This naturally leads us to the following question.

\begin{qes}\label{qes:crit}
    Can we obtain more refined information on when $\mathbb{D}_{0,n}\left(\vir{2k+1}, \otimes_{i=1}^{n}W_{a_i} \right)$ stabilizes?
\end{qes}

Note that the information on when a Virasoro conformal block divisor stabilizes provides a bound on the number of nonzero differences of Virasoro conformal block divisors. From this viewpoint, the stabilization aligns with the split of the F-conjecture, which implies that there are only finitely many extremal rays of the nef cone.

\section{Basis of \texorpdfstring{$\text{Pic}\left(\M{1}{n}\right)$}{TEXT} and applications}\label{sec:basis}

\subsection{\texorpdfstring{$\vir{5}$}{TEXT}-basis of \texorpdfstring{$\text{Pic}\left(\M{1}{n}\right)$}{TEXT}}\label{subsec:basis}

Let $V := \vir{5}$. This VOA has the central charge of $-\frac{22}{5}$ and is sometimes referred to as the Lee-Yang model. Recall that this VOA has exactly two irreducible modules: $W_1 = V$, the trivial module, and $W_2$. Both modules are self-dual, and their conformal weights are $0$ and $-\frac{1}{5}$, respectively. The fusion rules are given by
\[V \boxtimes V = V, \quad V \boxtimes W_2 = W_2 \boxtimes V = W_2, \quad W_2 \boxtimes W_2 = V \oplus W_2.\]

We will prove the following theorem. Note that we work only with the rational Picard group, not the integral one.

\begin{thm}\label{thm:basis}
If exactly one among $W^1, \dots, W^n$ is $W_2$, then $\mathbb{V}_{1,n}(V, W^\bullet)$ is a trivial line bundle. For any other choices, $\mathbb{D}_{1,n}(V, W^\bullet)$ is nontrivial. Moreover, the set of such nontrivial conformal block divisors is a basis for $\text{Pic}(\overline{\mathrm{M}}_{1,n})$.
\end{thm}

We recall some well-known facts to establish the proof of \cref{thm:basis}.

\begin{lem}\label{lem:linem1n}
    \begin{enumerate}
        \item On $\text{Pic}(\overline{\mathrm{M}}_{1,n})$, we have $12\lambda=\delta_{irr}$ and $12\psi_p=\delta_{irr}+12\sum_{p\in S, |S|\ge 2}\delta_{0,S}$.
        \item The boundary divisors form a basis of $\text{Pic}(\overline{\mathrm{M}}_{1,n})$. In particular, the rank of $\text{Pic}(\overline{\mathrm{M}}_{1,n})$ is $2^n-n$.
        \item Let $\pi_i: \overline{\rm{M}}_{1,n}\to \overline{\rm{M}}_{1,n-1}$ be the projection map that forgets the $i$th point. Then \[\pi_i^\ast \delta_{irr}=\delta_{irr},\ \pi_i^\ast \delta_{0,S}=\delta_{0,S}+\delta_{0,S\cup \left\{i\right\}}\]
        \item The rank of the vector bundles is given by:
        \[ \text{rank }\mathbb{V}_{0,n}(V, W_2^{\otimes n})=F_{n-1},\ \text{rank }\mathbb{V}_{1,n}(V, W_2^{\otimes n})=F_{n-1}+F_{n+1} \]
        where $F_n$ is the Fibonacci sequence, defined by $F_0 = 0, F_1 = 1$, and $F_{n+1} = F_n + F_{n-1}$.
    \end{enumerate}
\end{lem}

\begin{proof}
(1) This follows from Mumford's relation and \cite[Theorem 4(c)]{AC09}.

(2) Follows directly from (1) and \cite[Theorem 4(c)]{AC09}.

(3) Refer to \cite[Lemma 1]{AC09} for the proof.

(4) This is straightforwad by induction using the Factorization Theorem and the fusion rules.
\end{proof}

We begin with the first assertion of \cref{thm:basis}. By the propagation of vacua, it suffices to show that $\mathbb{V}_{1,1}(V, W)$ is a trivial line bundle. The rank of this vector bundle is $1$, as follows from the factorization property and the fusion rules. Therefore, it is enough to compute its first Chern class, which is $0$ by \cref{thm:virdeg} (2).

Next, we prove the rest of the theorem. Since the number of such conformal block divisors is $2^n - n$, it suffices to show that they generate $\text{Pic}(\overline{\mathrm{M}}_{1,n})$, by \cref{lem:linem1n} (2). We will use induction for this purpose. For the base case $n = 1$, \cref{thm:virdeg} implies that $\mathbb{D}_{1,1}(V, V)$ is nonzero, which establishes the result for $n = 1$. For the induction step, we will introduce and rely on two lemmas.

\begin{lem}\label{lem:oper}
    Let $T:\text{Pic}(\overline{\mathrm{M}}_{1,n})\to \Q$ be the linear functional defined by 
    \[ T(\delta_{irr})=0,\ T(\delta_{0,S})=(-1)^{|S|}. \]
    Then, for $n\ge 2$,
    \[ \ker T=\sum_{i=1}^{n}\pi_i^\ast \text{Pic}(\overline{\mathrm{M}}_{1,n-1}). \]
\end{lem}

\begin{proof}
    By \cref{lem:linem1n} (3), it is clear that $\ker T \supseteq \bigoplus_{i=1}^n \pi_i^\ast \text{Pic}(\overline{\mathrm{M}}_{1,n-1})$. The converse is also straightforward, by recursively applying $\pi_i^\ast \delta_{0,S} = \delta_{0,S} + \delta_{0,S \cup \left\{i\right\}}$.
\end{proof}

\begin{lem}\label{lem:Fibo}
    For $n\ge 1$, \[  \sum_{k=0}^n (-1)^{k+1}\binom{n}{k}F_kF_{n-k}=\begin{cases}
        0&\text{ if }n\text{ is odd }\\
        2\cdot 5^{\frac{n}{2}-1}&\text{ if }n\text{ is even },
    \end{cases} \]
    and
    \[  \sum_{k=0}^n (-1)^{k+1}\binom{n}{k}F_kF_{n+2-k}=\begin{cases}
        5^{\frac{n-1}{2}}&\text{ if }n\text{ is odd }\\
        3\cdot 5^{\frac{n}{2}-1}&\text{ if }n\text{ is even }.
    \end{cases}  \]
\end{lem}

\begin{proof}
    Let $a$ and $b$ be roots of the equation $x^2 - x - 1 = 0$. Note that $F_k = \frac{a^k - b^k}{a - b}$ and
    \[F_k F_{n-k} = \frac{1}{(a-b)^2}\left(a^n + b^n - a^k b^{n-k} - a^{n-k} b^k \right).\]
    Using standard binomial coefficient formulas, we have:
    \begin{align*}
    & \sum_{k=0}^n (-1)^{k+1} \binom{n}{k} a^n = \sum_{k=0}^n (-1)^{k+1} \binom{n}{k} b^n = 0, \\
    & \sum_{k=0}^n (-1)^{k+1} \binom{n}{k} a^k b^{n-k} = -(b-a)^n, \quad \sum_{k=0}^n (-1)^{k+1} \binom{n}{k} a^{n-k} b^k = -(a-b)^n.
    \end{align*}
    From these identities, the assertion follows directly by routine computation.
\end{proof}

Now, we will complete the proof. By \cref{lem:oper} and \cref{thm:propvac}, it is enough to show that
\[ T(\mathbb{D}_{1,n}(V,W_2^{\otimes n}))\ne 0 \]
Note that, by \cref{lem:linem1n} (1), if $n\ge 2$,
\[ T(\lambda)=0, T(\psi_p)=\sum_{p\in S, |S|\ge 2}(-1)^{|S|}=\sum_{S\subseteq \left\{1,2,...,n\right\}\setminus \left\{p\right\}, |S|\ge 1}(-1)^{|S|+1}=1. \]
By \cref{thm:c1cdf}, \cref{lem:linem1n} (4) and \cref{lem:Fibo},
\begin{align*}
    T(\mathbb{D}_{1,n}(V,W_2^{\otimes n}))&= -\frac{n}{5}\text{rank }\mathbb{V}_{1,n}(V,W_2^{\otimes n})-\sum_{|S|\ge 2}b_{0,S}(-1)^{|S|} \\
    &=-\frac{n}{5}(F_{n+1}+F_{n-1})+\frac{1}{5}\sum_{|S|\ge 2}(-1)^{|S|}F_{|S|}(F_{n-|S|+2}+F_{n-|S|})\\
    &=-\frac{1}{5}\left(n(F_{n+1}+F_{n-1})-\sum_{k=2}^{n}\binom{n}{k}(-1)^{k}F_{k}(F_{n-k+2}+F_{n-k}) \right)\\
    &= -\frac{1}{5}\sum_{k=0}^n (-1)^{k+1}\binom{n}{k}F_k(F_{n-k}+F_{n-k+2})=-5^{\lfloor \frac{n}{2} \rfloor-1}.
\end{align*}
This completes the proof.

\subsection{Line bundles on Contractions of \texorpdfstring{$\M{1}{n}$}{TEXT} }\label{subsec:cont}

Here, using the basis given by the conformal block divisors, we characterize line bundles on a natural contraction of $\M{1}{n}$ in the sense of \cite{Cho24}. We begin with some basic facts on contractions of $\M{g}{n}$.

Let $[n]=\left\{1,\cdots, n\right\}$ and let $S,T\subseteq [n]$. Assume $2g-2+|S|>0$ and $2g-2+|T|>0$. We consider these as indexing sets of points on a stable curve. Define

\[ \overline{\mathcal{M}}_{g,S,T}:=\begin{cases}
			 \overline{\mathcal{M}}_{g,S}\times_{ \overline{\mathcal{M}}_{g,S\cap T}} \overline{\mathcal{M}}_{g,T}, & \text{ if }g\ge 2\text{ or }g=1, |S\cap T|\ge 1 \text{ or }g=0, |S\cap T|\ge 3\\
             \overline{\mathcal{M}}_{0,S}\times  \overline{\mathcal{M}}_{0,T}, & \text{ if }g=0, |S\cap T|\le 2\text{ and }|S|, |T|\ge 3\\
            \overline{\mathcal{M}}_{1,S} \times_{\overline{\mathcal{M}}_{1,1}} \overline{\mathcal{M}}_{1,T}, & \text{ if }g=1, S\cap T=\emptyset, S,T\ne \emptyset\\
		 \end{cases}  \]
and define $\overline{\rm{M}}_{g,S,T}$ in the same way but using the coarse moduli space instead. The last case is well-defined since any projection map $\pi:\overline{\mathcal{M}}_{1,n}\to \overline{\mathcal{M}}_{1,1}$ that forgets $n-1$ points does not depend on the point they remember, and it corresponds to the $\lambda$ class. Note that this generalizes Knudsen's construction \cite{Knu83}. The case of $g=0$ is already discussed in \cite{Cho24}, following \cite{GF03}, and the case of $g=1$ also appeared in \cite{GKM02}. In particular, we have a natural surjective map $\overline{\mathcal{M}}_{g,n}\to \overline{\mathcal{M}}_{g,S,T}$ in every case. Let $X_{g,S,T}$ be the coarse moduli space of $\overline{\mathcal{M}}_{g,S,T}$. Then we have natural induced maps $\M{g}{n}\to X_{g,S,T}\to \overline{\rm{M}}_{g,S,T}$.

\begin{thm}\label{thm:gencont}
    $\M{g}{n}\to X_{g,S,T}\to \overline{\rm{M}}_{g,S,T}$ is the Stein factorization of $\M{g}{n}\to \overline{\rm{M}}_{g,S,T}$. In particular, $f_{g,S,T}: \M{g}{n}\to X_{g,S,T}$ is a contraction.
\end{thm} 

\begin{proof}
    The proof follows the general outline of the case for $g=0$, which is already established in \cite[Section 3.1]{Cho24}. Note that in this case, $X_{g,S,T}=\overline{\rm{M}}_{g, S,T}$ since the moduli stack is also a scheme. 

    \textbf{Step 1. } $\overline{\mathcal{M}}_{g,S,T}$ is Cohen-Macaulay.
    
    $\overline{\mathcal{M}}_{g,S,T}$ is defined as a fiber product of the form $A \times_B C$, where $A$, $B$, and $C$ are smooth Deligne-Mumford stacks. Since the inclusion $A \times_B C \to A \times C$ is the pullback of a regular embedding $B \to B \times B$, $A \times_B C$ is étale-locally covered by lci schemes; hence, in particular, it is a Cohen-Macaulay stack.    

    \textbf{Step 2. } $\overline{\mathcal{M}}_{g,S,T}$ is normal.

    We use induction on $|T \setminus S|$. If $T \subseteq S$, then there is nothing to prove. Assume the statement holds for $|T \setminus S| = n-1$, and consider the case $|T \setminus S| = n$. Choose $i \in T \setminus S$. We have a pullback diagram.

    \[\begin{tikzcd}[ampersand replacement=\&]
            \overline{\mathcal{M}}_{g,S,T}\arrow[r, "h"] \arrow[d]\&\overline{\mathcal{M}}_{g,T} \arrow[d]\\
            \overline{\mathcal{M}}_{g,S,T-i}\arrow[r]\& \overline{\mathcal{M}}_{g,T-i}.
        \end{tikzcd}.\]

    Note that there exists an open substack $\mathcal{U} \subseteq \overline{\mathcal{M}}_{g,T}$ such that $\pi_i: \overline{\mathcal{M}}_{g,T} \to \overline{\mathcal{M}}_{g,T-i}$ is smooth, and the codimension of its complement is $2$. We can let $\mathcal{U}$ be the locus where the point corresponding to $i$ is a smooth point. Then $h^{-1}(\mathcal{U})$ is smooth over $\overline{\mathcal{M}}_{g,S,T-i}$. By the induction hypothesis, $h^{-1}(\mathcal{U})$ is normal, and its complement has codimension $2$. Therefore, any affine étale cover of $\overline{\mathcal{M}}_{g,S,T}$ is regular in codimension $1$. By Step 1, it is also Cohen-Macaulay. Hence, by Serre's normality criterion, $\overline{\mathcal{M}}_{g,S,T}$ has a normal affine étale cover, hence normal.

    \textbf{Step 3. }$\M{g}{n}\to X_{g,S,T}$ is a contraction.

    Since any projection map $\overline{\mathrm{M}}_{g,n} \to \overline{\mathrm{M}}_{g,n-1}$ is a contraction, we may assume that $S \cup T = [n]$. In this case, it is straightforward to see that $\overline{\mathrm{M}}_{g,n} \to X_{g,S,T}$ is birational by restricting to the smooth locus. For any normal domain $A$ and a finite group $G$ acting on $A$, $A^G$ is also normal. Therefore, by Step 2, $X_{g,S,T}$ is also normal. Then, by the Zariski Main Theorem, $\overline{\mathrm{M}}_{g,n} \to X_{g,S,T}$ is a contraction.

    \textbf{Step 4. } $X_{g,S,T}\to \overline{\rm{M}}_{g, S,T}$ is finite.
    
    Since both are proper varieties, it is enough to show that the map is quasi-finite on closed points. We may assume that the base field $k$ is algebraically closed. By the definition of $\overline{\rm{M}}_{g,S,T}$, a closed point on $\overline{\rm{M}}_{g,S,T}$ is given by a pair $(x,y) \in \overline{\rm{M}}_{g,S} \times \overline{\rm{M}}_{g,T}$, possibly with an additional condition that they coincide after projecting to $\overline{\rm{M}}_{g,S\cap T}$, depending on $g$ and $|S\cap T|$. We will assume we have this additional condition because it is even easier without it. Let $z \in \overline{\rm{M}}_{g,S\cap T}$ be the projection of $x$ and $y$.
    
    Let $p\in \overline{\mathcal{M}}_{g,S}\times_{ \overline{\mathcal{M}}_{g,S\cap T}} \overline{\mathcal{M}}_{g,T}$ be a $k$-point that maps to $(x,y)\in \overline{\mathrm{M}}_{g,S}\times_{ \overline{\mathrm{M}}_{g,S\cap T}} \overline{\mathrm{M}}_{g,T}$. Then, any such $p$ corresponds to a pair $(C_x,C_y,f)$, where $C_x, C_y$ are stable curves corresponding to $x,y$, and $f$ is an isomorphism $\pi_1(C_x)\simeq \pi_2(C_y)$, where $\pi_1, \pi_2$ are maps forgetting points not in $S\cap T$. Since $\pi_1(C_x)$ is also stable, there are only finitely many choices for such $f$. Any inverse image of $(x,y)$ along $X_{g,S,T}\to \overline{\rm{M}}_{g, S,T}$ corresponds to such $p$, so this map is quasi-finite.
\end{proof}

\begin{eg}\label{eg:counter}
    In general, $X_{g,S,T}\to \overline{\rm{M}}_{g, S,T}$ is not an isomorphism, and $\M{g}{n}\to \overline{\rm{M}}_{g, S,T}$ is not a contraction. Let $(C, p_1,\cdots, p_n)$ be a stable curve with a nontrivial automorphism $f$, and assume $f$ does not preserve $q,r\in C$. Let $C_\bullet^1:=(C, p_1,\cdots, p_n,q)$ and $C_\bullet^2:=(C, p_1,\cdots, p_n,r)$. Then, $(C_\bullet^1,C_\bullet^2, \text{id})$ and $(C_\bullet^1,C_\bullet^2, f)$ are two different inverse images of $(C_\bullet^1,C_\bullet^2)$ along $X_{g,[n-1]\cup\left\{q\right\},[n-1]\cup\left\{r\right\}}\to \overline{\rm{M}}_{g, [n-1]\cup\left\{q\right\},[n-1]\cup\left\{r\right\}}$
\end{eg}

We will begin with the characterization of $\psi$ classes.

\begin{thm}\label{thm:charpsi}
    Let $D$ be a divisor on $\M{1}{n}$ such that $D|_{\Delta_{0, \left\{i,n\right\}}}=0$ for $1\le i\le n$ and $D\cdot F=0$ for one of the following F-curves:
    \[ F_5(0,0,I,J)\text{ for } I\sqcup J=[n-1], I,J\ne \emptyset. \]
    Then $D$ is a constant multiple of $\psi_n$.
\end{thm}

\begin{proof}
    Let $T\subset \text{Pic}(\M{1}{n})$ be the subspace that consists of line bundles that are trivial on $\Delta_{0, \left\{i,n\right\}}$ for every $i\in [n-1]$. First, we will prove $\dim T\leq 2$. Let $D\in T$ and
    \[ D=b_{\text{irr}}\delta_{\text{irr}}+\sum_{I\subseteq [n], |I|\ge 2}b_I \delta_{0,I}. \]
    Such an expression exists by \cref{lem:linem1n} (2). Now we will prove that $D$ is determined by $b_{\text{irr}}$ and $b_{[n-1]}$, which implies $\dim T\leq 2$.

    For $i \in [n-1]$, let $f_i:\M{1}{([n-1]\setminus\left\{i\right\})\cup\left\{p\right\}}\to \M{1}{n}$ be the map that corresponds to $\Delta_{0, \left\{i,n\right\}}$. Then, by \cite[Lemma 1]{AC09},
    \[ f_i^\ast \delta_{0,I}=\begin{cases}
        -\psi_p & \text{ if }I=\left\{i,n\right\}\\
        \delta_{0, (I\setminus \left\{i,n\right\})\cup \left\{p\right\}} & \text{ if }\left\{i,n\right\}\subseteq I\\
        \delta_{0,I} & \text{ if }\left\{i,n\right\}\subseteq I^c\\
        0 & \text{ otherwise }.
    \end{cases} \]
    Moreover, by \cref{lem:linem1n} (1), $12\psi_p = \delta_{\text{irr}} + 12\sum_{p \in S, |S| \geq 2}\delta_{0,S}$, hence if $\left\{i,n\right\} \subseteq I$ or $\left\{i,n\right\} \subseteq I^c$, then $b_{0,I}$ is determined by $b_{\text{irr}}$ via the condition $f_i^\ast D = 0$. Since this holds for every $i \in [n-1]$, $b_{\text{irr}}$ determines all coefficients but $b_{[n-1]}$. Hence, $\dim T \leq 2$.

    Let $F$ be one of the F-curves in the statement of the theorem, and $T_0 \subseteq T$ be the subspace consisting of line bundles intersecting trivially with $F$. Note that $\delta_{0, [n-1]} \in T$, but $\delta_{0, [n-1]} \cdot F = 1$, hence $\delta_{0, [n-1]} \not \in T_0$. Therefore, $T_0$ is a proper subspace of $T$, so $\dim T_0 \le 1$. Since $\psi_n \in T_0$, $T_0$ is exactly the subspace spanned by $\psi_n$. This completes the proof.
\end{proof}

Note that we can also prove \cref{thm:charpsi} using \cref{thm:basis}, instead of using boundary divisors as a basis.

Define $F_{S,T}$ as the set of F-curves contracted by $f_{S,T}$. In particular,
\[ F_{[n]\setminus \left\{i\right\}, [n]\setminus \left\{j\right\}}=\left\{F_5(0,0,\left\{i\right\}, \left\{j\right\})\right\} \cup \left\{F_6(1,0,0,0,I,J, \left\{i\right\}, \left\{j\right\})\ |\  I\sqcup J=[n]\setminus \left\{i,j\right\}, I\ne \emptyset  \right\}. \]

\begin{thm}\label{thm:charKnu}
    Let $f:\M{1}{n}\to X_{1, [n]\setminus \left\{i\right\}, [n]\setminus \left\{j\right\}}$ be the contraction in \cref{thm:gencont}. Then for a divisor $D$ on $\M{1}{n}$, the followings are equivalent:
    \begin{enumerate}
    \item $D$ trivially intersects $F_{[n]\setminus \left\{i\right\}, [n]\setminus \left\{j\right\}}$.
    \item $D$ is a pullback of a divisor along $f$.
    \item $D$ is a pullback of a divisor along $\M{1}{n}\to \M{1}{[n]\setminus \left\{i\right\}}\times \M{1}{[n]\setminus \left\{j\right\}}$. 
\end{enumerate}

\end{thm}

\begin{proof}
    First, let $h: \M{0}{n+1} \to \M{1}{n}$ be a map attaching a genus $1$ curve to the $(n+1)$th point. Then, for $I \subseteq [n]$ with $|I| \ge 2$, we have $h^\ast \delta_{0,I} = \delta_{0,I}$. Hence, $h^\ast : \text{Pic}(\M{1}{n}) \to \text{Pic}(\M{0}{n+1})$ is surjective, so the dual map $h_\ast : \NE{\M{0}{n+1}} \to \NE{\M{1}{n}}$ is injective. By \cite[Corollary 4.5]{Cho24},
    \[ F_{\text{Knu}}=\left\{F_6(0,0,0,0,I,J, \left\{i\right\}, \left\{j\right\})\ |\ I\sqcup J=[n]\setminus \left\{i,j\right\}, I,J\ne \emptyset \right\} \]
    is linearly independent. Therefore, 
    \[ h_\ast F_{\text{Knu}}=\left\{F_6(1,0,0,0,I,J, \left\{i\right\}, \left\{j\right\})\ |\  I\sqcup J=[n]\setminus \left\{i,j\right\}, I\ne \emptyset  \right\} \]
    is a linearly independent subset of $\NE{\M{1}{n}}$.

    Now, we will prove that $F_{[n] \setminus {i}, [n] \setminus {j}}$ is a linearly independent subset of $\NE{\M{1}{n}}$. Since $h_\ast F_{\text{Knu}}$ is linearly independent, it suffices to show that there exists a line bundle on $\text{Pic}(\M{1}{n})$ that intersects $h_\ast F_{\text{Knu}}$ trivially but intersects $F_5(0,0,\left\{i\right\}, \left\{j\right\})$ nontrivially. Define $\pi: \M{1}{n} \to \M{1}{2}$ as the map that forgets all points except the $i$th and $j$th. For any ample line bundle $L$ on $\M{1}{2}$, the pullback $\pi^\ast L$ has the desired property: $\pi$ contracts $h_\ast F_{\text{Knu}}$ but does not contract $F_5(0,0,{i},{j})$. Therefore, $F_{[n] \setminus {i}, [n] \setminus {j}}$ is linearly independent.
    
    There are $2^{n-2}$ elements in $F_{[n] \setminus {i}, [n] \setminus {j}}$, and $\text{Im } f^\ast$ intersects trivially with $F_{[n] \setminus {i}, [n] \setminus {j}}$. Hence, the dimension of $\text{Im } f^\ast$ is $\leq 2^n - 2^{n-2} - n$. Moreover, since $f_0 : \M{1}{n} \to \M{1}{n-1} \times \M{1}{n-1}$ factors through $f$, it follows that $\text{Im } f^\ast$ contains the image of
    \[ f_0^\ast:\text{Pic}\left(\M{1}{n-1} \right)\times \text{Pic}\left(\M{1}{n-1} \right)\to \text{Pic}\left(\M{1}{n} \right). \]
    By \cref{thm:basis}, the image of this map is the space spanned by
    \[ \left\{\mathbb{D}_{1,n}(\vir{5}, W^\bullet)\ |\ W^i=V \text{ or } W^j=V, \text{ number of }W_2\text{ among }W^\bullet\text{ is }\ne 1  \right\}. \]
    It is straightforward to verify that exactly $2^n - 2^{n-2} - n$ elements are in this set. Therefore, $\dim \text{Im } f^\ast = 2^n - 2^{n-2} - n$, which implies that $\text{Im } f^\ast$ consists precisely of the line bundles that intersect trivially with $F_{[n] \setminus {i}, [n] \setminus {j}}$.
\end{proof}

\begin{cor}\label{cor:knuex}
    We have the following exact sequence:
    \[
    0 \to \text{Pic}(\M{1}{n-2}) \xrightarrow{(\pi_n^\ast, -\pi_{n-1}^\ast)} \text{Pic}(\M{1}{n-1}) \times \text{Pic}(\M{1}{n-1}) \xrightarrow{\pi_{n-1}^\ast + \pi_n^\ast} \text{Pic}(\M{1}{n}) \xrightarrow{c} \Q^{F_{[n]\setminus \left\{n\right\}, [n]\setminus \left\{n-1\right\}}} \to 0,
    \]
    where $c$ denotes the intersection pairing. The image of $\pi_{n-1}^\ast + \pi_n^\ast$ is $\text{Pic}(X_{1, [n]\setminus \left\{n\right\}, [n]\setminus \left\{n-1\right\}})$.
\end{cor}

\begin{proof}
By the proof of \cref{thm:charKnu}, the sequence
\[
\text{Pic}(\M{1}{n-1}) \times \text{Pic}(\M{1}{n-1}) \xrightarrow{\pi_{n-1}^\ast + \pi_n^\ast} \text{Pic}(\M{1}{n}) \xrightarrow{c} \Q^{F_{[n] \setminus \{n\}, [n] \setminus \{n-1\}}}
\]
is exact. Moreover, the codimension of $\text{Pic}(X_{1, [n] \setminus \{n\}, [n] \setminus \{n-1\}})$ in $\text{Pic}(\M{1}{n})$ equals the cardinality of $F_{[n] \setminus \{n\}, [n] \setminus \{n-1\}}$. Hence, the last term in the sequence is also exact. By definition, the image of $(\pi_n^\ast, -\pi_{n-1}^\ast)$ is contained in the kernel of $\pi_{n-1}^\ast + \pi_n^\ast$. Therefore, by dimension counting, the first term is also exact.
\end{proof}

\begin{thm}\label{thm:chargen}
    Let $f:\M{1}{n}\to X_{1, S, T}$ be a contraction in \cref{thm:gencont}. Then for a divisor $D$ on $\M{1}{n}$, the followings are equivalent:
    \begin{enumerate}
    \item $D$ trivially intersects with the F-curves contracted by $f$.
    \item $D$ is a pullback of a divisor along $f$.
    \item $D$ is a pullback of a divisor along $\M{1}{n}\to \M{1}{S}\times \M{1}{T}$. 
\end{enumerate}
\end{thm}

\begin{proof}
    First, consider the case where $S = T$, so $f$ is the projection map $\pi: \M{1}{n} \to \M{1}{S}$. By \cref{thm:basis},
    \[ I=\left\{\mathbb{D}_{1,n}(\vir{5}, W^\bullet)\ |\ \text{ number of }W_2\text{ among }W^\bullet\text{ is }\ne 1  \right\} \]
    is a basis of $\text{Pic}(\M{1}{n})$. Let $L=\sum_{D\in I}a_D \cdot D \in \text{Pic}(\M{1}{n})$ be a line bundle that contracts $F_{S,S}$. By \cref{thm:propvac}, it suffices to prove that if $a_D \ne 0$ for $D = \mathbb{D}_{1,n}(\vir{5}, W^\bullet)$, then $W^i = V$ for every $i \in S^c$.

    Assume instead that $W^i = W_2$ for some $i \in S^c$. By the definition of $I$, there exists $j \ne i$ such that $W^j = W_2$. Since $F_{[n] \setminus {i}, [n] \setminus {j}} \subseteq F_{S,S}$, \cref{thm:charKnu} implies that $a_D = 0$. This contradiction completes the proof for the case $S = T$.

    Now consider the general case. Since $f$ factors through the projection map $\pi: \M{1}{n} \to \M{1}{S \cup T}$, by the previous case, we may assume $[n] = S \cup T$, so $S^c \cap T^c = \emptyset$. Let $L \in \text{Pic}(\M{1}{n})$ be a line bundle that contracts $F_{S,T}$, and write $L = \sum_{D \in I} a_D D$. By \cref{thm:propvac}, it suffices to prove that if $a_D \ne 0$ for $D = \mathbb{D}_{1,n}(\vir{5}, W^\bullet)$, then $W^i = V$ for every $i \in S^c$ or $W^i = V$ for every $i \in T^c$. This implies that $L$ is contained in the image of
    \[ \text{Pic}\left(\M{1}{S} \right)\times \text{Pic}\left(\M{1}{T} \right)\to \text{Pic}\left(\M{1}{n} \right). \]
    Assume there exist $i \in S^c$ and $j \in T^c$ such that $W^i = W^j = W_2$. Since $S^c \cap T^c = \emptyset$, we have $i \ne j$. Since $F_{[n] \setminus {i}, [n] \setminus {j}} \subseteq F_{S,S}$, \cref{thm:charKnu} implies that $a_D = 0$. This contradiction completes the proof.
\end{proof}

\begin{rmk}
    \cref{thm:chargen} provides a useful criterion for the bigness of a nef divisor on $\M{1}{n}$. By \cite{GKM02}, a nef divisor is not big if and only if it is a pullback of a divisor on $\M{1}{S} \times \M{1}{S^c}$, and this is equivalent to the vanishing of its intersection with certain $F$-curves by \cref{thm:chargen}.
\end{rmk}

From \cref{thm:charpsi} and \cref{thm:chargen}, the following statement directly follows.

\begin{cor}\label{cor:extremal}
    \begin{enumerate}
        \item $\psi_i$ forms an extremal ray of the nef cone of $\M{1}{n}$.
        \item If $D$ is a divisor that forms an extremal ray of the nef cone of $\M{1}{n}$, then the pullback of $D$ along any projection map $\pi:\M{1}{m}\to \M{1}{n}$ is an extremal ray of the nef cone of $\M{1}{m}$. In particular, any pullback of $\psi$ classes forms an extremal ray.
    \end{enumerate}
\end{cor}

\section{Further Discussion}\label{sec:dis}

\subsection{What is special about \texorpdfstring{$\text{Vir}_{2,2k+1}$}{TEXT}?}\label{subsec:vir}

As noted in \cref{rmk:rank}, $\vir{2k+1}$ has some distinguished properties, even among discrete series Virasoro algebras. In particular, on $\M{0}{4}$ or $\M{1}{1}$, its degree is determined by its rank. This special feature is a crucial ingredient of in the proof of many other properties, such as the stabilization of Virasoro conformal block divisors. Moreover, conformal block divisors for $\vir{2k+1}$ are nef, unlike for affine VOAs, for which coinvariant divisors are nef. For affine VOAs, the positivity of divisors in genus $0$ follows from the global generation of the sheaf of coinvariants, whereas the positivity of the conformal blocks divisors for $\vir{2k+1}$ is proved combinatorially and also holds for $g>0$. Therefore, it is natural to ask whether we can develop a more conceptual understanding of  Virasoro conformal blocks, particularly their Chern classes and positivity. With this in mind, we suggest some questions that naturally arise from the results of this paper.

In physics, $\vir{2k+1}$ has some distinguished properties. In such contexts, $\text{Vir}_{p,q}$ corresponds to the $(p,q)$-minimal model, and by coupling it with the appropriate Liouville theory, we obtain a two-dimensional topological gravity, which is a special case of string theory. Moreover, any such theory has a corresponding matrix model that gives the same partition function up to a double scaling limit. If $p = 2$, this theory exhibits interesting properties that are not shared with other minimal models. First, its partition function satisfies the KdV hierarchy (see \cite[Chapter 5.4]{Eyn16} for details), which is why it is sometimes called the KdV minimal model. Second, it is modeled over a single Hermitian matrix, thus admitting a simpler matrix model. We refer to \cite{FGZ95} for a more detailed explanation of two-dimensional gravity associated with minimal models. Given these special properties of KdV minimal models in physics, it is natural to ask what the consequences of these properties are in terms of the sheaf of conformal blocks.

\cref{thm:virdeg} provides a very explicit formula for the degree in terms of the rank, which is $-\frac{r(r-1)}{2}$ (resp. $-r(r-1)$) on $\M{0}{4}$ (resp. $\M{1}{1}$). Note that this corresponds to the formula for the summation of the first $r-1$ negative integers. Since any vector bundle on $\overline{\mathcal{M}}_{0,4}$ or $\overline{\mathcal{M}}_{1,1}$ is a sum of line bundles (the first case is a classical theorem of Grothendieck, and the second case is given by \cite{Me15}), we propose the following conjecture.

\begin{conj}\label{conj:vb04}
    \[ \mathbb{V}_{0,4}(\vir{2k+1}, W^\bullet)=\oplus_{i=0}^{r-1}\mathcal{O}(-i)\text{ and } \mathbb{V}_{1,1}(\vir{2k+1}, W)=\oplus_{i=0}^{r-1}\lambda^{-2i}. \]
    In terms of sheaf of conformal blocks,
    \[ \mathbb{V}_{0,4}(\vir{2k+1}, W^\bullet)^\vee=\oplus_{i=0}^{r-1}\mathcal{O}(i)\text{ and } \mathbb{V}_{1,1}(\vir{2k+1}, W)^\vee=\oplus_{i=0}^{r-1}\lambda^{2i} \]
\end{conj}

In \cref{subsec:crit}, we prove that both $\rank \mathbb{V}_{0,n}(\vir{2k+1}, \otimes_{i=1}^n W_{a_i})$ and $\mathbb{D}_{0,n}(\vir{2k+1}, \otimes_{i=1}^n W_{a_i})$ stabilize as $k \to \infty$. These correspond to the zeroth and first Chern classes, respectively. Hence, we pose the following question.

\begin{qes}\label{qes:vecstab}
    Does $c_r \mathbb{V}_{0,n}(\vir{2k+1}, \otimes_{i=1}^n W_{a_i})$ stabilize as $k \to \infty$ for any $r \geq 2$? Will the sequence of vector bundles $\mathbb{V}_{0,n}(\vir{2k+1}, \otimes_{i=1}^n W_{a_i})$ itself stabilize?
\end{qes}

Note that in the case of affine VOAs, the fiber of the sheaf of coinvariants $\mathbb{V}_{g,n}(L_k(\mathfrak{g}), W^\bullet)$ at a stable curve $(C, p_\bullet)$ has a geometric description. More precisely, 
\[\mathbb{V}_{g,n}(L_k(\mathfrak{g}), W^\bullet) |_{(C, p_\bullet)} = \text{H}^0(\mathcal{B}un_G^{\text{Par}}(C, p_\bullet), \mathcal{L}^k),\]
where $G$ is a simple and simply connected algebraic group, $\mathfrak{g}$ its Lie algebra, $\mathcal{B}un_G^{\text{Par}}(C, p_\bullet)$ is the moduli space of parabolic $G$-bundles on $(C, p_\bullet)$, and $\mathcal{L}$ is the canonical line bundle on $\mathcal{B}un_G^{\text{Par}}(C, p_\bullet)$. The choice of modules influences the definition of the parabolic bundle as well as the selection of the canonical line bundle. This result is proven for smooth curves in \cite{BL94, LS97} and for stable curves in \cite{BF19, BG19}. Notably, \cite{BF19} proved this for families of stable curves, not only for a single fiber. A similar description for the Heisenberg VOA is provided in \cite{Uen95}. In this case, the fibers are given by the theta functions. Hence, the fibers of the sheaf of coinvariants are sometimes called generalized or nonabelian theta functions. Therefore, it is natural to ask whether there exists a `Virasoro theta function' description of the fiber of sheaf of conformal blocks of $\vir{2k+1}$.

\begin{qes}\label{qes:conffib}
    Is there an analogous geometric `generalized theta function' description of the fibers of the sheaf of conformal blocks $\mathbb{V}_{g,n}(\vir{2k+1}, W^\bullet)^\vee$ at a stable curve arising, for instance, from a proof of \cref{conj:vb04}?   
\end{qes}

Note that the first part of the question for general VOA has already been asked in \cite{DGT22b}.

In the case of affine VOAs, we proved that the sheaves of coinvariants are globally generated by naturally defined global sections, which arise from the degree $0$ part of the modules. By \cite{Zhu96}, these correspond to the simple representations of the Zhu algebra, which, in this case, is a quotient of the universal enveloping algebra of $\mathfrak{g}$. The same construction applies to $\vir{2k+1}$. However, by \cite{FZ92}, the Zhu algebra of $\vir{2k+1}$ is commutative, so all its simple representations are $1$-dimensional. Consequently, this yields only one section of the sheaf of coinvariants for $\vir{2k+1}$. This outcome is consistent with the results of this paper, as the sheaf of conformal blocks admits positivity. This naturally raises the question: can we construct global sections of the sheaf of conformal blocks?

\begin{qes}\label{qes:GG}
    Can we construct natural global sections of $\mathbb{V}_{g,n}(\vir{2k+1}, W^\bullet)^\vee$, for instance, using \cref{qes:conffib}? Moreover, can we determine whether $\mathbb{V}_{g,n}(\vir{2k+1}, W^\bullet)^\vee$ is globally generated by such global sections?
\end{qes}

One more important property of critical level conformal block divisors of affine VOAs is the duality \cite{BGM15, BGM16}. Since the stabilization of the Virasoro conformal block divisor is analogous to the critical level, we can ask for the corresponding duality.

\begin{qes}\label{qes:dual}
    Can we formulate critical level duality for Virasoro conformal block divisors?
\end{qes}

\subsection{Inductive System of Line Bundles on \texorpdfstring{$\M{0}{n}$}{TEXT}}\label{subsec:ind}

In this paper, we utilize the factorization (\cref{thm:factor}) and the propagation of vacua (\cref{thm:propvac}) to prove various properties of conformal block divisors. Apart from the explicit representation of such divisors in terms of tautological classes (\cref{thm:c1cdf}), these are the only properties we have used. Hence, it is natural to ask whether a set of line bundles exists that satisfy these properties but do not arise from VOAs.

\begin{defn}\label{defn:ind}
    An \textbf{inductive system of line bundles on }$\M{0}{n}$ is a triple $\mathcal{S}=(S, f, D)$ such that
    \begin{enumerate}
        \item $S$ is a set with a distinguished element $V\in S$.
        \item $f$ is a function $f:S\times S\to \mathbb{R}_{\ge 0}[S]$. Here, $\mathbb{R}_{\ge 0}[S]$ is the monoid generated by $\mathbb{R}_{\ge 0}$ and $S$.
        \item $D=(D_n)$ is a collection of $S_n$-equivariant functions $D_n:S^n\to \text{Pic}(\M{0}{n})_\R$ for $n\ge 3$
    \end{enumerate}
    which satisfies
    \begin{description}
        \item[Multiplication] The $\R$-bilinear extension $\bar{f}:\R[S]\times \R[S]\to \R[S]$ of $f$ makes $\R[S]$ into a commutative $\R$-algebra with identity $V$.
    \end{description}
    Define $f(M_1,\cdots, M_n)(M)$ to be the coefficient of $M$ in $M_1\cdots M_n\in \R[S]$, which is always nonnegative. 
    \begin{description}
        \item[Dual] For any $M\in S$, there exists unique $M'\in S$ such that $f(M,M')(V)\ne 0$, and $f(M,M')(V)=1$.
        \item[Propagation of Vacua] $D_n(M_1,\cdots, M_{n-1}, V)=\pi_n^\ast D_{n-1}(M_1,\cdots, M_{n-1})$
        \item[Factorization] For the gluing map  $\xi:\overline{\mathrm{M}}_{0,n_1+1}\times \overline{\mathrm{M}}_{0,n_2+1}\to \overline{\mathrm{M}}_{0,n}$,
        \begin{align*}
            &\xi^\ast D_n\left(M_1,\cdots, M_n \right) \simeq\\& \bigotimes_{M\in S}\left(\pi_1^\ast D_{n_1+1}\left(M_1,\cdots, M_n, M \right)^{\otimes f(M_{n_1},\cdots, M_n)(M)}\otimes\pi_2^\ast D_{n_2+1}\left(M_{n_1+1},\cdots, M_n, M' \right)^{\otimes f(M_{1},\cdots, M_{n_1})(M')}  \right)
        \end{align*}
    \end{description}
\end{defn}

The following simple proposition gives a useful property that follows from the axioms.

\begin{prop}\label{prop:genfact}
    For an inductive system of line bundle $\mathcal{S}=(S, f, D)$,
    \[ f(M_1,\cdots, M_n)(V)=\sum_{M\in S}f(M_{1},\cdots, M_{n_1})(M')\cdot f(M_{n_1},\cdots, M_n)(M) \]
\end{prop}

\begin{proof}
    The statement directly follows from the multiplication axiom and the dual axiom.
\end{proof}

\begin{defn}
    An inductive system of line bundles $\mathcal{S} = (S, f, D)$ \textbf{admits conformal weights} if there exists a function $\text{cw}: S \to \mathbb{Q}$ such that
    \[ D_n(M_1,\cdots, M_n)=f(M_1,\cdots, M_n)(V)\sum_{i=1}^n \text{cw}(M_i)\psi_i-\sum_{I} b_{I}\delta_{0,I} \]
    where
    \[ b_{I}=\sum_{M\in S} \text{cw}(M) \cdot  f(M_I, M)(V) \cdot f(M_{I^c}, M')(V). \]
    Conversely, a triple $(S, f, \text{cw})$ is \textbf{admissible} if, upon defining $D_n(M_1, \cdots, M_n)$ as above, the resulting $(S, f, D)$ is an inductive system of line bundles. In this case, we refer to $(S, f, D)$ as the \textbf{associated inductive system}.
\end{defn}

\begin{eg}\label{eg:VOA}
    Let $V$ be a $C_2$-cofinite and rational VOA, and let $S$ denote the set of simple modules over $V$. Define $f: S \times S \to \mathbb{R}_{\ge 0}[S]$ as the tensor product map, and let $D_n(M_1, \cdots, M_n) := \mathbb{D}_{0,n}(V, \otimes_{i=1}^n M_i)$. This defines an inductive system of line bundles by \cref{thm:propvac} and \cref{thm:factor}. In this context, $\mathbb{Z}[S]$ is the fusion ring of $V$. This example illustrates that the inductive system of line bundles generalizes the genus $0$ conformal block divisors of a VOA. This inductive system admits conformal weights.
\end{eg}

\begin{eg}
    In \cite{Fe15}, Fedorchuk defined the `divisors from symmetric functions on abelian groups', forming an inductive system of line bundles. Specifically, the set $S$ is an abelian group, and $f$ is defined by addition. This construction can be viewed as a generalization of conformal block divisors arising from a lattice VOA. Thus, considering \cref{eg:VOA} as a fundamental example of an inductive system of line bundles, the notion of divisors from a symmetric function is a precursor to this concept. This inductive system also admits conformal weights, which is the symmetric function.
\end{eg}

\begin{eg}
    The relationship between the inductive system of line bundles and cohomological field theory (CohFT) \cite{KM94} is particularly intriguing. We will use the convention of the Chow field theory \cite[Subsection 1.4]{Pan18}. The inductive system of line bundles exhibits several differences from CohFT.

    \begin{enumerate}
        \item It is restricted to the genus $0$ case.  
        \item It is also restricted to Chow degrees $\le 1$
        \item It includes a positivity requirement given by $f: S \times S \to \R_{\geq 0} [S]$.  
    \end{enumerate}

    The first restriction is not fundamental, as the definition can readily be extended to higher genus. The author did not extend it in this paper because of the lack of new examples of inductive systems of line bundles for higher genera. The second restriction is a distinct feature. It facilitates the construction of new examples more easily. Specifically, for any CohFT, we can naturally associate an inductive system of line bundles by truncating it above chow degree $2$. Hence, the inductive system of line bundles can be seen as an analogue of TQFT to degree $2$, just as TQFT is a truncation of CohFT above degree $1$.

    The third restriction is the most significant difference between the two. This condition requires the CohFT to possess a special basis that satisfies certain positivity properties. In contrast, a general CohFT does not need to fulfill this requirement. This distinction arises because CohFTs are designed to address enumerative problems, while the inductive system of line bundles focuses on the positivity of line bundles. 
\end{eg}

One advantage of \cref{defn:ind} is that the set of stable Virasoro conformal block divisors forms an inductive system of line bundles.

\begin{prop}
    $\mathbb{D}_{0,n}(a_1,\cdots, a_n)$ forms an inductive system of line bundles. Moreover, this admits conformal weights defined by cw$(i)=-\frac{i-1}{2}$.
\end{prop}

\begin{proof}
    The set $S$ is the set of natural numbers $\mathbb{N}$, where each $i \in \mathbb{N}$ is identified with $W_{i+1}$, and $W_1$ is the distinguished element. The fusion rules define the map $f$
    \[W_i \boxtimes W_j = \sum_{r \stackrel{\text{2}}{=}j-i+1}^{i+j-1} W_r.\]
    Since this $f$ represents the stabilized fusion rules of $\vir{2k+1}$, it satisfies both the multiplication and dual axioms. The propagation of vacua follows from \cref{thm:propvac}. Furthermore, as shown in the proof of \cref{prop:rkstab}, $f(W_{a_1}, \cdots, W_{a_n})(V)$ corresponds to the $\rank \mathbb{V}_{0,n}(\vir{2k+1}, W^\bullet)$, where $k$ is above the critical level. Thus, factorization follows from \cref{thm:factor} as shown in the proof of \cref{thm:stnef}. Hence, they form an inductive system of line bundles. The second assertion follows from the proof of \cref{thm:steff}.
\end{proof}

Now, we present a family of examples of inductive systems of line bundles. The underlying idea is as follows: coinvariant divisors of the $L_1(\mathfrak{sl}_2)$ and conformal block divisors of $\vir{5}$ are nef on genus $0$. Both have exactly two simple modules, $V$ and $W$, where $V$ is the trivial module and $W$ is the nontrivial module. The fusion rules are given by 
\[V \boxtimes V = V, \quad V \boxtimes W = W \boxtimes V = W, \quad W \boxtimes W = V + pW,\]
where $p = 0$ and $p = 1$, respectively. 

These examples naturally lead to the question: for each $p \in \Z_{>1}$, does there exist a $C_2$-cofinite and rational VOA $V$ with exactly two simple modules, $V$ and $W$, satisfying the fusion rules as above and are their conformal block divisors nef on genus $0$?  Constructing such a VOA is generally a challenging task, and proving that it is $C_2$-cofinite and rational is even more challenging. Consequently, the author is unaware of any examples with $p = 2$. The advantage of the concept of an inductive system of line bundles is that we do not need to address these difficulties directly. Using \cref{thm:c1cdf}, we can imitate the behavior of their conformal block divisors, making it sufficient to construct an inductive system of line bundles instead.

\begin{thm}\label{thm:indsys}
   Let $p \in \R_{\geq 0}$, $S = \{V, W\}$, and define a function $f: S \times S \to \R_{\geq 0}[S]$ by
    \[f(V, V) = V, \quad f(V, W) = f(W, V) = W, \quad f(W, W) = V + pW.\]
    Additionally, define
    \[D_n(W, \dots, W) = R_n^p \cdot \psi - \sum_{i=2}^{\lfloor \frac{n}{2} \rfloor} R_{i+1}^p \cdot R_{n-i+1}^p \cdot \delta_{0,i},\]
    where
    \[R_n^p := \frac{1}{\sqrt{p^2 + 4}} \left( \left( \frac{p + \sqrt{p^2 + 4}}{2} \right)^{n-1} - \left( \frac{p - \sqrt{p^2 + 4}}{2} \right)^{n-1} \right),\]
    and extend to
    \[D_n(M_1, \dots, M_n) = \pi_I^\ast D_{|I|}(W, \dots, W),\]
    where $I \subseteq [n]$ is the set of indices such that $M_i = W$. Then, $(S, f, D)$ forms an inductive system of line bundles.
\end{thm}

\begin{proof}
    First, we note that
    \[ f(V,\cdots,V,W,\cdots, W)(V) = f(W,\cdots, W)(V) = R_n^p,\  R_{n+1}^p=pR^p_n+R_{n-1}^p \]
    when there are $n$ many $W$'s. This follows directly from the definition of $f$ and by induction.

    All other axioms, except for factorization, are trivial. In particular, $D_n$ is defined to satisfy the propagation of vacua. Moreover, by the definition of $D$, when checking the factorization axiom, it suffices to check the case where $M_1 = \cdots = M_n = W$ due to the definition of $D_n$. Hence, we will compute the restriction of $D_n(W,\cdots, W)$ to $\Delta_{0,I}$.

    Since $\Delta_I \simeq \M{0}{|I|+1} \times \M{0}{|I^c|+1}$ and $\text{Pic}(\M{0}{|I|+1} \times \M{0}{|I^c|+1}) \simeq \text{Pic}(\M{0}{|I|+1}) \times \text{Pic}(\M{0}{|I^c|+1})$ (cf. \cite{Ke92}), it suffices to show that
    \[\iota^\ast D_n(W,\cdots, W) = R_{n-|I|+1}^p D_{|I|+1}(W,\cdots, W) + R_{n-|I|}^p \pi_p^\ast D_{|I|}(W,\cdots, W) \]
    where $\iota: \M{0}{I \cup \{x\}} \to \M{0}{n}$ is the map attaching a stable rational $|I^c|+1$-pointed curve at $x$, and $\pi_p: \M{0}{I \cup \{p\}} \to \M{0}{I}$ is the projection map. By \cite[Lemma 1]{AC09},
    \[ \iota^\ast D_n(W,\cdots, W) = R_n^p \cdot \sum_{i \in I} \psi_i + R_{|I|+1}^p \cdot R_{n-|I|+1}^p \psi_x - \sum_{\substack{J \subseteq I, \\ 2 \le |J| \le |I|-1}} R_{|J|+1}^p \cdot R_{n-|J|+1}^p \delta_{0,J}, \]
    \[  R_{n-|I|+1}^p D_{|I|+1}(W,\cdots, W) = R_{n-|I|+1}^p \cdot R_{|I|+1}^p \cdot \psi - \sum_{\substack{J \subseteq I, \\ 2 \le |J| \le |I|-1}} R_{n-|I|+1}^p \cdot R_{|J|+1}^p \cdot R_{|I|-|J|+2}^p \cdot \delta_{0,J}\]
    \[ R_{n-|I|}^p \pi_p^\ast D_{|I|}(W,\cdots, W) = R_{n-|I|}^p R_{|I|}^p \sum_{i \in I} \psi_i - \sum_{\substack{J \subseteq I, \\ 2 \le |J| \le |I|-1}} R_{n-|I|}^p \cdot R_{|J|+1}^p \cdot R_{|I|-|J|+1}^p \cdot \delta_{0,J}. \]
    Hence, this follows from
    \[ R_m^p = R_{m-i+1}^p R_{i+1}^p + R_{m-i}^p R_{i}^p, \]
    which is straightforward by induction.

\end{proof}

\begin{cor}\label{cor:indpos}
    \[  D_n^p:=R_n^p \cdot \psi - \sum_{i=2}^{\lfloor \frac{n}{2} \rfloor} R_{i+1}^p \cdot R_{n-i+1}^p \cdot \delta_{0,i} \]
    is a nef line bundle for $p\ge 0$. and an ample line bundle for $p>0$.
\end{cor}

\begin{proof}
    Note that $R_{i+1}^p \cdot R_{n-i+1}^p < R_n^p$, so this line bundle is a sum of $\psi - \delta$ and an effective sum of boundaries. By factorization, as in the proof of \cref{lem:VOAGKM}, for an F-curve $F = F(0,0,0,0,I_1,I_2,I_3,I_4)$ of Type 6, we have
    \[D_n^p \cdot F = \prod_{i=1}^4 R_{|I_i|+1}^p \cdot D_4(W,W,W,W),\]
    so this is an F-nef divisor, and it is F-ample for $p > 0$ since $R_n^p > 0$ for $n \geq 1$ in this case.

    We will use \cref{lem:effnef} or \cref{lem:Knef} with induction on $n$. If the base field has characteristic 0, we can apply \cref{lem:Knef}, and then $D_n^p$ satisfies the hypothesis (note that $K_{\M{0}{n}} = \psi - 2\delta$), so it is enough to prove that the boundary restrictions are nef and ample when $p > 0$. This quickly follows the factorization and the induction hypothesis. In general, we can apply \cref{lem:effnef} since the standard form of $\psi - \delta$ gives an effective sum of boundaries.
\end{proof}

Note that if $p = 0$ (resp. $p = 1$), these divisors are proportional to $L_1(\mathfrak{sl}_2)$-coinvariant divisors (resp. $L_1(G_2)$, $L_1(F_4)$-coinvariant divisors or $\vir{5}$-conformal block divisors). Hence, they can be viewed as a nonlinear interpolation and extension of such conformal block/coinvariant divisors. An important advantage of this approach is that it works even for non-integer values of $p$. Thus, this method is both simpler and more general than constructing VOAs with the corresponding fusion rules. In particular, $D_n^p$ cannot come from vertex operator algebras.

\printbibliography

\end{document}